\documentclass[11pt,letterpaper]{article}
\usepackage[T1]{fontenc}
\usepackage{lmodern}
\usepackage[margin=1in]{geometry}
\usepackage{amsmath,amssymb,amsthm}
\usepackage{booktabs,array,tabularx}
\usepackage{microtype,needspace}
\usepackage{xspace}
\usepackage{xcolor}
\definecolor{crossrefblue}{RGB}{0,51,102}
\usepackage[colorlinks=true,linkcolor=crossrefblue,citecolor=black,urlcolor=black]{hyperref}
\usepackage[table]{hypcap}
\hypersetup{
  pdflinkmargin=0.5pt,
  pdftitle={Enumerating Solutions to Type (2,3) Composition Identities of Finite Endofunctions},
  pdfauthor={Max A. Alekseyev}
}

\DeclareRobustCommand{\namedref}[2]{\hyperref[#2]{#1~\ref*{#2}}}

\newtheorem{theorem}{Theorem}[section]
\newtheorem{lemma}[theorem]{Lemma}
\newtheorem{corollary}[theorem]{Corollary}
\theoremstyle{definition}
\newtheorem{remark}[theorem]{Remark}
\DeclareMathOperator{\im}{im}
\DeclareMathOperator{\Fix}{Fix}
\DeclareMathOperator{\Comp}{Comp}
\DeclareMathOperator{\Aut}{Aut}
\newcommand{\fall}[2]{#1^{\underline{#2}}}
\newcommand{\coef}[1]{\mathop{[\![#1]\!]}\xspace}
\newcommand{\NN}{\mathcal N}
\newcommand{\MM}{\mathcal M}
\newcommand{\oeis}[1]{\href{https://oeis.org/#1}{OEIS~\texttt{#1}}}
\newcommand{\oeiss}[1]{\href{https://oeis.org/#1}{\texttt{#1}}}
\newcolumntype{Y}{>{\raggedright\arraybackslash}X}

\title{Enumerating Solutions to Type \((2,3)\)\\Composition Identities of Finite Endofunctions}
\author{Max A. Alekseyev\thanks{The George Washington University,
Washington, DC, United States.
Corresponding author:
\href{mailto:maxal@gwu.edu}{maxal@gwu.edu}.}}
\date{}

\begin{document}
\maketitle

\begin{abstract}
We enumerate ordered pairs \((f,g)\) of endofunctions on a finite
labeled set \(X\) satisfying a composition identity of type \((2,3)\),
with two factors on one side and three on the other.
For all sixteen forms up to exchanging \(f\) and \(g\), we give exact
counting formulas as functions of \(|X|\) in terms of finite sums and
generating functions. Such pairs correspond bijectively to actions
on \(X\) of the semigroup presented by two generators and the
prescribed relation. The proofs use images, fibers, and the
decomposition of functional digraphs into directed cycles with rooted
trees attached. Several constructions apply beyond these sixteen
identities: they give recurrences for homomorphisms between arbitrary
prescribed functional digraphs and formulas for pairs satisfying
\(f^a=g^b\), where \(a,b\) are fixed positive integers. We also derive
an explicit inclusion--exclusion formula over subsets of \(X\) for
the number of endofunctions \(g\) on \(X\) satisfying \(B=gAg\)
for arbitrary prescribed endofunctions \(A,B\) on \(X\). This
generalizes the enumeration of compositional square roots, recovered
when \(A\) is the identity map.
\end{abstract}



\section{Introduction}
\label{sec:scope}

\paragraph*{Notation and background.}
Let \(n\) be a nonnegative integer and put \(X_n:=\{1,\ldots,n\}\).
An \emph{endofunction} on \(X_n\) is a function \(X_n\to X_n\),
which we refer to simply as a map on \(X_n\). Let \(\MM_n\) be the
set of these maps. We express as a function of \(n\) the number of
ordered pairs \((f,g)\in \MM_n^2\) satisfying an identity that equates
compositions of two and three factors from \(\{f,g\}\). We refer to
these as identities of type \((2,3)\).
Identities of type \((1,2)\) admit elementary counting formulas.
In our previous work~\cite{AlekseyevCompanion}, we completed the
enumeration of solutions to identities of type \((1,3)\).
Earlier work also treats identities of type \((2,2)\), including
the commutation relation \(fg=gf\)~\cite{HollowayShattuck2013}.
These results make type \((2,3)\) a natural next family to study
systematically. We give counting formulas for all sixteen normalized
forms.

\begin{table}[!t]
\centering
\small
\setlength{\tabcolsep}{4pt}
\renewcommand{\arraystretch}{1.08}
\begin{tabular}{@{}cccc|cc@{}}
\toprule
Identity & \shortstack{Numerical data\\in the OEIS~\cite{OEIS2026}} & \shortstack{Total\\formula}
& \shortstack{Arithmetic\\complexity}
& \shortstack{Equivalent condition\\in \(\mathfrak S_n\)}
& \shortstack{Number\\of pairs}\\
\midrule
\multicolumn{4}{@{}c|}{\textit{Images, fixed points, and sources (\namedref{Section}{sec:images-sources})}} & \multicolumn{2}{c@{}}{}\\
\(f^2=f^3\) & \(n^n\cdot\)\,\oeiss{A000949} & \eqref{eq:unary23-total} & \(O(n^2)\) & \(f=\mathrm{id}\) & \(n!\)\\
\(fg=g^3\) & \oeiss{A239750} & \eqref{eq:total750}\textsuperscript{*} & \(O(n\log n)\) & \(f=g^2\) & \(n!\)\\
\(fg=fg^2\) & \oeiss{A239785} & \eqref{eq:companion785-sum}\textsuperscript{*} & \(O(n^2)\) & \(g=\mathrm{id}\) & \(n!\)\\
\(fg=f^2g\) & \oeiss{A239781} & \eqref{eq:total781}\textsuperscript{*} & \(O(n^2)\) & \(f=\mathrm{id}\) & \(n!\)\\
\(fg=gfg\) & \oeiss{A239759} & \eqref{eq:total759}\textsuperscript{*} & \(O(n^2)\) & \(g=\mathrm{id}\) & \(n!\)\\
\(f^2=gf^2\) & \oeiss{A239752} & \eqref{eq:total752}\textsuperscript{*} & \(O(n^2)\) & \(g=\mathrm{id}\) & \(n!\)\\
\(f^2=g^2f\) & \oeiss{A239756} & \eqref{eq:E-positive}\textsuperscript{*} & \(O(n^3)\) & \(f=g^2\) & \(n!\)\\
\midrule
\multicolumn{4}{@{}c|}{\textit{Allocations of fibers (\namedref{Section}{sec:fibers})}} & \multicolumn{2}{c@{}}{}\\
\(f^2=fgf\) & \oeiss{A239772} & \eqref{eq:companion772-explicit} & \(O(n^2)\) & \(g=\mathrm{id}\) & \(n!\)\\
\(f^2=f^2g\) & \oeiss{A239775} & \eqref{eq:P}, \eqref{eq:allocation-tree-cycle-total} & \(O(n^5)\) & \(g=\mathrm{id}\) & \(n!\)\\
\(fg=fgf\) & \oeiss{A239783} & \eqref{eq:783total}, \eqref{eq:allocation-tree-cycle-total} & \(O(n^5)\) & \(f=\mathrm{id}\) & \(n!\)\\
\(fg=f^3\) & \oeiss{A239764} & \eqref{eq:764det} & \(O(n^6)\) & \(g=f^2\) & \(n!\)\\
\midrule
\multicolumn{4}{@{}c|}{\textit{Homomorphisms and iterate types (\namedref{Section}{sec:hom})}} & \multicolumn{2}{c@{}}{}\\
\(fg=gf^2\) & \oeiss{A239755} & \eqref{eq:homTotalFrom} & \(O(n^\frac32 3.0^n)\) & \(g^{-1}fg=f^2\) & \(n!q_n\)\\
\(fg=g^2f\) & \oeiss{A239762} & \eqref{eq:homTotalTo} & \(O(n^\frac32 3.0^n)\) & \(fgf^{-1}=g^2\) & \(n!q_n\)\\
\(f^2=g^3\) & \oeiss{A235325} & \eqref{eq:equalpowerhistogram} & \(O(n^\frac32 3.0^n)\) & \(f^2=g^3\) & \(U_n^{\mathrm{bij}}\)\\
\midrule
\multicolumn{4}{@{}c|}{\textit{Colored trees with additional constraints (\namedref{Sections}{sec:preserved}--\ref{sec:remaining})}} & \multicolumn{2}{c@{}}{}\\
\(f^2=fg^2\) & \oeiss{A239766} & \eqref{eq:766GFtotal} & \(O(n^2 4.8^n)\) & \(f=g^2\) & \(n!\)\\
\(f^2=gfg\) & \oeiss{A239749} & \eqref{eq:coupled-kernel-total}, \eqref{eq:coupled-color-total} & \(O(n^4 5.7^n)\) & \((fg)^2=f^3\) & \(U_n^{\mathrm{bij}}\)\\
\bottomrule
\end{tabular}
\caption{All normalized forms of type \((2,3)\), in presentation order.
Asterisks mark references giving explicit finite sums with a fixed
number of summation indices. 
The arithmetic complexities are justified in \namedref{Appendix}{sec:operation-bounds}.
The last two columns require both maps to be permutations, including the unused map in the $f^2=f^3$ row.}
\label{tab:inventory}
\end{table}

We use products of maps to denote composition from right to left: for
\(f,g\in\MM_n\) and \(x\in X_n\), we set \((fg)(x):=f(g(x))\).
Positive powers of a map denote its repeated composition.
To record the counts, for fixed nonempty words \(u,v\) over
\(\{f,g\}\), let \(u(f,g),v(f,g)\) denote the corresponding
compositions and define
\[
 \NN_{u,v}(n):=\#\{(f,g)\in\MM_n^2:u(f,g)=v(f,g)\}.
\]
We define the \emph{arithmetic complexity} of a formula as the number of
arithmetic operations required to evaluate the \(n\)th term.
The definition of \(\NN_{u,v}(n)\) gives the brute-force formula
\[
 \NN_{u,v}(n)=
 \sum_{f,g\in\MM_n}\mathbf1_{\{u(f,g)=v(f,g)\}},
\]
where \(\mathbf1_E\) is the \emph{indicator} of a condition \(E\),
equal to \(1\) when \(E\) holds and \(0\) otherwise.
This gives a baseline arithmetic complexity of \(n^{2n}\), up to
polynomial factors.
Our aim is to obtain more explicit formulas for \(\NN_{u,v}(n)\) with substantially lower arithmetic complexity, ideally closed-form expressions in \(n\).

Under composition, \(\MM_n\) is the full transformation semigroup
on \(X_n\). Actions on \(X_n\) of the semigroup presented by
two generators and the relation \(u=v\) correspond bijectively
to solution pairs \((f,g)\): the generators act by \(f\) and \(g\),
respectively. Thus \(\NN_{u,v}(n)\) is the number of these actions.

The counting arguments use two descriptions of a map. For
\(h\in\MM_n\) and \(y\in h(X_n)\), the \emph{fiber} of \(h\)
over \(y\) is \(h^{-1}(y):=\{x\in X_n:h(x)=y\}\).
We also represent \(h\) by its \emph{functional digraph}
\(G(h)\), with vertex set \(X_n\) and edges \(x\to h(x)\) for
\(x\in X_n\).
Connected components of \(G(h)\) are defined by ignoring edge
directions. Each such component consists of a directed cycle with
rooted trees directed towards it. We write \(\Comp(h)\) for the
set of connected components of \(G(h)\). A \emph{source} of
\(G(h)\) is a vertex with no incoming edge.

We compare functional digraphs by \emph{homomorphisms}, maps between
their vertex sets that preserve directed edges. Functional digraphs
have been used for transformation equations~\cite{Higgins1988},
and their homomorphisms were studied
in~\cite{Novotny1963,Nelson1982}.
The equal-power identity \(f^2=g^3\) also leads to roots of maps
under composition. Square roots were studied by Snowden and
Howie~\cite{SnowdenHowie1982}, Higgins~\cite{Higgins1986,Higgins1988},
and Kozyra~\cite{Kozyra2018}.

\paragraph*{The sixteen identities and counting methods.}

Exchanging \(f\) and \(g\) if necessary, we may assume that the
shorter word starts with \(f\). It is then \(ff\) or \(fg\),
and the longer word is one of eight words of length three. This gives
sixteen normalized forms. The form \(f^2=f^3\) leaves \(g\)
unrestricted, and \(f^2=g^3\) uses the two maps on separate sides;
each other form has a side containing both maps.
Reversing all words is not a counting symmetry for arbitrary maps.

Our previous study~\cite{AlekseyevCompanion} counts solutions of \(fg=fg^2\) by sources
and connected components of functional digraphs, and solutions of \(f^2=fgf\) by the sizes of fibers and
the positions of their image values. We extend these arguments by
allowing the source and component weights to vary independently and
by allowing fiber weights to depend on the sizes of the fibers
assigned to each target fiber. The earlier formulas arise as special
cases of these extensions, which also apply to further identities.
For the remaining cases we use homomorphisms and information about iterates.

\namedref{Table}{tab:inventory} groups the sixteen cases by the methods used
to prove their counting formulas and records their restrictions to
permutations. In the row for \(f^2=f^3\),
the factor \(n^n\) accounts for the unrestricted second map.
The Online Encyclopedia of Integer Sequences (OEIS)~\cite{OEIS2026}
records numerical data and earlier formulas for several of these counts.
We adopt and give combinatorial proofs of Heinz's formulas for
\(fg=g^3\) (\oeis{A239750}) and \(fg=f^2g\) (\oeis{A239781}), and Einstein's sum for \(fg=gfg\) (\oeis{A239759}).
For \(f^2=gf^2\) (\oeis{A239752}), we give a double sum with 
nonnegative terms, improving upon Einstein's triple sum, together with a
generating function that records the number of points outside the second image.
The numbers of maps satisfying \(f^2=f^3\) form the bounded-height
forest sequence \oeis{A000949}; these numbers follow from the general
enumeration of maps with equal prescribed iterates by Krattenthaler and
M\"uller~\cite{KrattenthalerMueller2004}.

Besides completing this enumeration, the paper gives counting
constructions with more general inputs. For allocations of blocks
whose weights depend on the target size and the number and total size
of incoming blocks, decomposition into trees and cycles gives a
determinant generating function (\namedref{Theorem}{thm:allocation-tree-cycle});
the weights for \(f^2=f^2g\) and \(fg=fgf\) give two
specializations. A recurrence on rooted trees counts homomorphisms
between arbitrary prescribed functional digraphs
(\namedref{Theorem}{thm:homDP}), and hence those between a map and its
square needed for \(fg=gf^2\) and \(fg=g^2f\).
Grouping maps by the isomorphism types of their iterates gives counting formulas
for \(f^a=g^b\), where \(a,b\) are fixed positive integers
(\namedref{Theorem}{thm:equalpowerhistogram}). This includes pairs
with equal squares (\oeis{A235327}) or equal cubes (\oeis{A235326}),
as well as the case \(f^2=g^3\) (\oeis{A235325}).

Note that our previous study~\cite{AlekseyevCompanion} also surveys the counting
formulas for identities of type \((1,3)\), with one factor on one
side and three on the other. Its formula for the remaining case
\(f=gfg\) completes the enumeration of all eight normalized forms.
For arbitrary prescribed maps \(A,B\in\MM_n\), we generalize
this problem by deriving an explicit inclusion--exclusion formula
for the number of maps \(g\in\MM_n\) satisfying \(B=gAg\)
(\namedref{Corollary}{cor:coupled-prescribed-sandwich}).
For fixed \(A\in\MM_n\), the associative product
\(r\star s:=rAs\), for \(r,s\in\MM_n\), defines a variant
of the full transformation semigroup~\cite{DolinkaEast2015}.
Solving \(B=gAg\) is therefore equivalent to finding square roots
of \(B\) under this product.
When \(A\) is the identity map, the formula counts compositional
square roots of \(B\). For fixed positive integers \(a,b\),
taking \(A=f^b\) and \(B=f^a\) and summing over \(f\in\MM_n\)
counts pairs satisfying \(f^a=gf^bg\), extending the enumeration
beyond the sixteen identities considered here.

\paragraph*{Restriction to permutations.}

We also count solutions in which both maps are permutations. For maps
from a finite set to itself, injectivity, surjectivity, and bijectivity
are equivalent. Let \(\mathfrak S_n\) be the set of permutations of
\(X_n\), write
\(\mathrm{id}\) for the identity map on \(X_n\), and define
\[
 \NN^{\mathrm{bij}}_{u,v}(n):=
 \#\{(f,g)\in\mathfrak S_n^2:u(f,g)=v(f,g)\}.
\]
The twelve factorial counts in \namedref{Table}{tab:inventory} follow by
cancellation: one permutation is arbitrary and the other is forced.
Later remarks also consider requiring only one map to be bijective.
We obtain the remaining four counts using conjugacy and roots of permutations.

For \(fg=gf^2\), the permutation \(f\) must be conjugate to its
square. Squaring splits every even cycle into two cycles and preserves
the length of every odd cycle. Thus the number of cycles increases
unless all lengths are odd, in which case the cycle type is unchanged.
It follows that \(f\) and \(f^2\) are conjugate exactly when all
cycles of \(f\) have odd length. For such a cycle type, let \(m_\ell\) be the number of cycles of
length \(\ell\), for each positive odd integer \(\ell\), and put
\(z_{\mathbf m}:=\prod_\ell\ell^{m_\ell}m_\ell!\).
There are \(n!/z_{\mathbf m}\) choices of \(f\) and
\(z_{\mathbf m}\) conjugating permutations \(g\), so each type
contributes \(n!\). The same argument applies to \(fg=g^2f\),
with \(g\) conjugate to its square. Both counts are therefore
\(n!q_n\) (\oeis{A088311}), where \(q_n\) is the number of
partitions of \(n\) into odd parts.

The permutation solutions of \(fg=gf^2\) correspond bijectively to
homomorphisms from the Baumslag--Solitar group
\(\langle a,t\mid t^{-1}at=a^2\rangle\) to \(\mathfrak S_n\),
with \(a\mapsto f\) and \(t\mapsto g\).
Its subgroup numbers are given by Gelman's formula, reproved by
Kelley~\cite[Theorem~3.1]{Kelley2020}.

Put \(U_n^{\mathrm{bij}}:=\NN^{\mathrm{bij}}_{ff,ggg}(n)\),
recorded as \oeis{A239837}. Permutation solutions of \(f^2=g^3\)
correspond bijectively to
homomorphisms from the \((2,3)\)-torus knot group to \(\mathfrak S_n\).
Permutation-root formulas are given by
Pavlov~\cite{Pavlov1981} and Lea\~nos, Moreno, and
Rivera-Mart\'inez~\cite{LeanosMorenoRivera2012}.
Baker and Petri~\cite[Proposition~A.1]{BakerPetri2023} use root counts
for exact enumeration and also study the asymptotic behavior of these
representations. The count \(U_n^{\mathrm{bij}}\) is evaluated in
\namedref{Remark}{rem:permutation-equalpowers}.

Finally, \(f^2=gfg\) is equivalent to \((fg)^2=f^3\).
The assignment \((f,g)\mapsto(fg,f)\) maps permutation solutions
of \(f^2=gfg\) bijectively onto pairs \((a,b)\in\mathfrak S_n^2\)
satisfying \(a^2=b^3\), with inverse \((a,b)\mapsto(b,b^{-1}a)\).
Hence this identity also has \(U_n^{\mathrm{bij}}\) permutation solutions.

\section{Images, fixed points, and sources}
\label{sec:images-sources}

\subsection{Prescribed compositions}

For prescribed maps \(a,b,c,d\) on a finite set \(X\), the equation
\(au=b\) restricts each value of the unknown map \(u\) on \(X\)
to a fiber of \(a\), whereas
\(uc=d\) prescribes \(u\) on the image of \(c\).
To express the compatibility condition for such prescribed values,
define the \emph{kernel relation} of a map \(a\) on \(X\) by
\[
 \ker a:=\{(x,y)\in X^2:a(x)=a(y)\}.
\]
Its equivalence classes are the nonempty fibers of \(a\).

\Needspace{12\baselineskip}
\begin{lemma}
\label{lem:prescribed-compositions}
Let \(a,b,c,d\) be prescribed maps on a finite set \(X\) of size \(n\).
Counting maps \(u\) on \(X\), we have
\begin{align}
 \#\{u:au=b\}
   &=\prod_{x\in X}|a^{-1}(b(x))|,
       \label{eq:prescribed-left}\\
 \#\{u:uc=d\}
   &=\begin{cases}
       n^{n-|c(X)|},&\ker c\subseteq\ker d,\\
       0,&\text{otherwise},
      \end{cases}
       \label{eq:prescribed-right}\\
 \#\{u:au=b,\ uc=d\}
   &=\begin{cases}
       \displaystyle\prod_{x\in X\setminus c(X)}|a^{-1}(b(x))|,
           &\ker c\subseteq\ker d\text{ and }ad=bc,\\[4pt]
       0,&\text{otherwise}.
      \end{cases}
       \label{eq:prescribed-two-sided}
\end{align}
The count in
\eqref{eq:prescribed-two-sided} is positive exactly when
\[
 \ker c\subseteq\ker d,\qquad ad=bc,\qquad b(X)\subseteq a(X).
\]
\end{lemma}

\begin{proof}
The equation \(au=b\) permits exactly the choices
\(u(x)\in a^{-1}(b(x))\), independently for each \(x\in X\),
which proves~\eqref{eq:prescribed-left}.
The equation \(uc=d\) prescribes
\(u(c(x)):=d(x)\) on \(c(X)\). This prescription is well-defined
exactly when \(\ker c\subseteq\ker d\), and leaves all other
values free. This proves~\eqref{eq:prescribed-right}.

For both equations, the prescribed values on \(c(X)\) also satisfy
\(au=b\) exactly when
\(a(d(x))=b(c(x))\) for every \(x\in X\), that is, when \(ad=bc\).
The remaining values lie independently in the displayed fibers, proving
\eqref{eq:prescribed-two-sided}. Under the two compatibility conditions,
\(b(c(X))\subseteq a(X)\) already holds. Thus all remaining fibers are
nonempty exactly when \(b(X)\subseteq a(X)\), proving the positivity
criterion.
\end{proof}

When the unknown map is a permutation, write
\(\operatorname{Sym}(X)\) for the permutations of \(X\), and put
\(a_y:=|a^{-1}(y)|\) and \(b_y:=|b^{-1}(y)|\) for \(y\in X\).
Then
\begin{align}
 \#\{u\in\operatorname{Sym}(X):au=b\}
 &=\begin{cases}\displaystyle\prod_{y\in X}a_y!,&a_y=b_y\ (y\in X),\\
 0,&\text{otherwise},\end{cases}
 \label{eq:bijective-prescribed-left}\\
 \#\{u\in\operatorname{Sym}(X):uc=d\}
 &=\begin{cases}(n-|c(X)|)!,&\ker c=\ker d,\\
 0,&\text{otherwise}.\end{cases}
 \label{eq:bijective-prescribed-right}
\end{align}
A permutation satisfying both equations exists exactly when
\(a_y=b_y\) for every \(y\in X\), \(\ker c=\ker d\), and \(ad=bc\).
Under these conditions, the number of such permutations is
\begin{equation}
 \prod_{y\in X}\bigl(a_y-|a^{-1}(y)\cap d(X)|\bigr)!;
 \label{eq:bijective-prescribed-two-sided}
\end{equation}
otherwise it is zero. Indeed, \(u\) must biject each fiber of \(b\)
with the corresponding fiber of \(a\). The prescription
\(u(c(x))=d(x)\) is injective exactly when the two kernel relations agree;
after it is imposed, the remaining points in each corresponding pair
of fibers can be matched arbitrarily.

More generally, imposing \(au(x)=b(x)\) only for \(x\) in a
prescribed set \(D\subseteq X\) gives
\begin{equation}
 \#\{u:(au)|_D=b|_D\}
   =n^{n-|D|}\prod_{x\in D}|a^{-1}(b(x))|.
 \label{eq:prescribed-on-subset}
\end{equation}
The values outside \(D\) are free; on \(D\), the same independent
fiber choices apply. For a permutation \(u\), put
\(m_y:=|\{x\in D:b(x)=y\}|\) for \(y\in X\).
If \(m_y\leq a_y\) for every \(y\in X\), the corresponding count is
\begin{equation}
 (n-|D|)!\prod_{y\in X}\frac{a_y!}{(a_y-m_y)!};
 \label{eq:bijective-prescribed-subset}
\end{equation}
otherwise it is zero. The product is the number of allowed injections
\(D\to X\); each extends to \((n-|D|)!\) permutations of \(X\).

An equality with a common right factor depends only on its image:
for any maps \(a,b,c\) on \(X\),
\[
 ac=bc\quad\Longleftrightarrow\quad a|_{c(X)}=b|_{c(X)}.
\]
For a prescribed map \(h\) on \(X\) and integers \(r\geq s\geq1\), the
kernel condition in~\eqref{eq:prescribed-right} is automatic, giving
\[
 \#\{u:uh^s=h^r\}=n^{n-|h^s(X)|}.
\]
This covers \(f^2=gf^2\) and \(fg=g^3\). More generally, for a fixed
map \(g\) on \(X\), the number of maps \(f\) on \(X\) satisfying
\(fg=g^r\) is independent of the positive integer \(r\).

Filmus gave this prescribed-image argument and the independence from
the positive exponent in a 2014 answer~\cite{Filmus2014}; the observation
is also recorded in \oeis{A239750}.

Applying these observations to the four identities below leaves
restrictions on images and fixed points. For \(h\in\MM_n\), write
\(\Fix(h):=\{x\in X_n:h(x)=x\}\).

\begin{lemma}\label{lem:local}
For fixed \(f\in\MM_n\), or fixed \(g\in\MM_n\) as appropriate,
the following counts over the other map in \(\MM_n\) hold:
\begin{align}
 \#\{g:f^2=gf^2\}&=n^{n-|\im(f^2)|},\label{eq:local752}\\
 \#\{f:fg=gfg\}&=n^{n-|\im g|}|\Fix(g)|^{|\im g|},\label{eq:local759}\\
 \#\{g:fg=f^2g\}&=|\{x\in X_n:f(x)=f^2(x)\}|^n,\label{eq:local781}\\
 \#\{f:fg=g^3\}&=n^{n-|\im g|}.\label{eq:local750}
\end{align}
\end{lemma}

\begin{proof}
In~\eqref{eq:local752}, \(g\) fixes \(\im(f^2)\) pointwise and is free
elsewhere. In~\eqref{eq:local759}, \(f\) sends \(\im g\) into \(\Fix(g)\)
and is free elsewhere. Equation~\eqref{eq:local781} restricts every value of
\(g\), independently, to the displayed set. Finally, \(fg=g^3\) prescribes
\(f(y)=g^2(y)\) on \(\im g\), with no restriction outside that image.
\end{proof}

For the corresponding permutation counts, put \(r:=|g(X_n)|\) in
the rows where \(g\) is fixed. Then
\begin{align}
 \#\{g\in\mathfrak S_n:f^2=gf^2\}
   &=(n-|\im(f^2)|)!,\notag\\
 \#\{f\in\mathfrak S_n:fg=gfg\}
   &=r!(n-r)!\,\mathbf1_{\{g^2=g\}},\notag\\
 \#\{g\in\mathfrak S_n:fg=f^2g\}
   &=n!\,\mathbf1_{\{f^2=f\}},\notag\\
 \#\{f\in\mathfrak S_n:fg=g^3\}
   &=(n-r)!\,\mathbf1_{\{g^2(X_n)=g(X_n)\}}.
 \label{eq:bijective-local-images}
\end{align}
For the first row, \(g\) fixes \(\im(f^2)\) pointwise and
permutes its complement.
For the second, \(f\) must send \(g(X_n)\) injectively into
\(\Fix(g)\subseteq g(X_n)\), so these sets coincide and \(f\)
permutes them. The third follows by cancelling \(g\).
For the last, the prescribed restriction \(g^2|_{g(X_n)}\) is
injective exactly when \(g\) permutes its image.

\subsection{The unary identity}
\label{sec:unary23}

For \(f\in\MM_n\), the equation \(f^2=f^3\) says that every
point of \(X_n\) reaches a fixed point after at most two applications
of \(f\). Thus every component of \(G(f)\) is a rooted tree of
height at most two, with a
loop at its root. These maps are counted by
\oeis{A000949}. Their exponential generating function
is also the specialization \(\alpha=2\), \(\beta=3\), with trivial
group, of Krattenthaler and M\"uller~\cite[Theorem~1]{KrattenthalerMueller2004}.

For a formal power or Laurent series \(F\) in an indeterminate \(z\)
and an integer \(m\), write \(\coef{z^m}F\) for the coefficient
of \(z^m\) in \(F\).
\begin{lemma}
\label{lem:unary23}
For every \(n\geq0\),
\begin{equation}
 \NN_{ff,fff}(n)
 =n^n n!\coef{z^n}\exp\!\left(z\exp(z\exp(z))\right).
 \label{eq:unary23-total}
\end{equation}
\end{lemma}

\begin{proof}
A rooted tree of height at most one is a root with a set of leaves,
so its exponential generating function is \(ze^z\). A rooted tree
of height at most two is a root with a set of such rooted trees, giving
\(z\exp(z\exp(z))\). Taking a set of components gives the exponential
in~\eqref{eq:unary23-total}. The second map \(g\) is arbitrary,
and contributes the factor \(n^n\).
\end{proof}

To evaluate the nested exponentials, let \(C(z)\) be a formal series with
\(C(0)=0\), set \(B(z):=\exp(C(z))\), and, for integers \(m\geq0\),
put \(b_m:=m!\coef{z^m}B(z)\), \(c_m:=m!\coef{z^m}C(z)\).
Differentiation gives the coefficient recurrence
\[
 b_0=1,\qquad
 b_m=\sum_{j=1}^m\binom{m-1}{j-1}c_jb_{m-j}
 \qquad(m\geq1).
\]
Apply it first to \(C(z)=ze^z\), whose coefficients are
\(c_j=j\) for \(j\geq1\), and then to \(C(z)=z\exp(z\exp(z))\), whose coefficients
are \(j b_{j-1}\) from the first application. Multiplication by
\(n^n\) gives~\eqref{eq:unary23-total}.

This is the height-two specialization of Heinz's bounded-height
forest recurrence in \oeis{A210725}. Kruchinin's double sum in
\oeis{A000949} gives another evaluation.

If only \(g\) is required to be a permutation, replace the factor
\(n^n\) in~\eqref{eq:unary23-total} by \(n!\). If only \(f\)
is required to be a permutation, then \(f=\mathrm{id}\) and there
are \(n^n\) choices of the unrestricted \(g\).

\subsection{Sources and connected components}
\label{sec:source-component-totals}

The counts for \(fg=g^3\) and \(fg=fg^2\) can be expressed through
sources and connected components. For \(h\in\MM_n\), the sources
of \(G(h)\) are precisely the points of \(X_n\setminus h(X_n)\). Write \(s(h)\) for their number
and \(c(h)\) for the number of connected components of \(G(h)\).

For \(fg=fg^2\), we first identify which points must have the same
value under \(f\). The following weighted form of this observation
will also be useful later.

\begin{lemma}
\label{lem:weighted-equalities}
Let \(a,b\colon W\to X\) be maps between finite sets, and let
\(Y\) be a finite set. Form an undirected graph on \(X\) by
joining \(a(w)\) to \(b(w)\) for every \(w\in W\), and let
\(\mathcal C(a,b)\) denote the set of its connected components,
including isolated vertices. For weights \(\rho(x,y)\), indexed by
\((x,y)\in X\times Y\) and taking values in a commutative ring,
\begin{equation}
 \sum_{\substack{u\colon X\to Y\\ua=ub}}
       \prod_{x\in X}\rho(x,u(x))
   =\prod_{C\in\mathcal C(a,b)}
       \left(\sum_{y\in Y}\prod_{x\in C}\rho(x,y)\right).
 \label{eq:weighted-equalities}
\end{equation}
Repeated edges and loops do not affect \(\mathcal C(a,b)\).
\end{lemma}

\begin{proof}
The equation \(ua=ub\) states that \(u\) takes the same value
at the endpoints of every edge. It is therefore equivalent to
\(u\) being constant on each component of the graph. Choosing a
common value \(y\in Y\) for a component \(C\in\mathcal C(a,b)\) contributes \(\prod_{x\in C}\rho(x,y)\),
and distinct components can be assigned independently. Summing over
their common values and multiplying proves the formula.
\end{proof}

For \(fg=fg^2\), apply the lemma with
\(u=f\), \(a=g\), \(b=g^2\), and all weights equal to one.
Every source of \(G(g)\) is isolated in the undirected graph
constructed in the lemma. The nonsources in
each component of \(G(g)\) form one component of this graph: they
are closed under forward iteration and remain joined to the cycle in that component of
\(G(g)\). There are consequently \(s(g)+c(g)\) components, giving
\begin{equation}
 \#\{f:fg=fg^2\}=n^{s(g)+c(g)}.
 \label{eq:source-component-fixed}
\end{equation}

If \(u\) is injective, then \(ua=ub\) is equivalent to \(a=b\).
Thus, when \(X=Y\) and all weights are one, the permutation count
in the lemma is \(|X|!\) if \(a=b\), and zero otherwise.
For \(f,g\in\MM_n\) satisfying \(fg=fg^2\), requiring \(f\)
to be a permutation therefore
leaves exactly the idempotent maps \(g\). Requiring \(g\) to be
a permutation instead gives \(f=fg\), so there are \(n^{c(g)}\)
choices of \(f\), one value for each cycle. Summing over such \(g\)
gives \(\prod_{j=0}^{n-1}(n+j)\).

To sum the fixed-map counts, we use the exponential generating function
for rooted labeled trees. Let \(\mathcal T(z)\in\mathbb Q[[z]]\) be the unique
formal series with zero constant term satisfying
\[
 \mathcal T(z)=z\exp(\mathcal T(z)).
\]
The following coefficient form of Lagrange inversion will be used
several times; see~\cite[Sections~2.1 and~3.2]{Gessel2016}.

\begin{lemma}
\label{lem:abel-tree-coefficients}
For indeterminates \(b,c\), let \(V(z)\) be the unique formal series
with zero constant term satisfying \(V=z\exp(bV)\). Then, for
integers \(0\leq r\leq m\),
\[
 \coef{z^m}\frac{V^r\exp(cV)}{1-bV}
 =\frac{(c+bm)^{m-r}}{(m-r)!}.
\]
\end{lemma}

\begin{proof}
The coefficient substitution \(z=v\exp(-bv)\) contributes the
factor \(1-bv\), which cancels the denominator. Hence
\[
 \coef{z^m}\frac{V^r\exp(cV)}{1-bV}
 =\coef{v^{m-r}}\exp((c+bm)v),
\]
as required.
\end{proof}

We now weight sources and connected components independently,
extending~\cite[Sections~3--4]{AlekseyevCompanion}. The finite-sum
formula uses rising factorials: for an indeterminate \(x\) and an
integer \(j\geq1\), write \(x^{\overline j}:=x(x+1)\cdots(x+j-1)\),
with \(x^{\overline0}:=1\).

\begin{lemma}
\label{lem:sources-components}
For indeterminates \(a,b\),
\begin{equation}
 \sum_{m\geq0}\left(\sum_{h\in\MM_m}
               a^{s(h)}b^{c(h)}\right)\frac{z^m}{m!}
 =\left(1-\mathcal T\!\left(z\exp((a-1)z)\right)\right)^{-b}.
 \label{eq:sources-components}
\end{equation}
For \(m\geq0\), its coefficients also satisfy
\[
 \sum_{h\in\MM_m}a^{s(h)}b^{c(h)}
 =\sum_{k=0}^m\binom mk(a-1)^{m-k}
   \sum_{j=0}^k\binom kj(b-1)^{\overline j}k^{m-j}.
\]
For positive integers \(a,b\) and \(m\geq0\),
\begin{equation}
 \sum_{h\in\MM_m}a^{s(h)}b^{c(h)}
 =m!\coef{t^m}
   \frac{(a-1+\exp(t))^{m+b-1}}
        {(a-1+(1-t)\exp(t))^{b-1}}.
 \label{eq:sources-components-elementary}
\end{equation}
\end{lemma}

\begin{proof}
Maps with weight \(b\) per connected component have exponential
generating function \((1-\mathcal T(z))^{-b}\), by the decomposition
into directed cycles with attached rooted trees.
For a positive integer \(a\), give every source one of \(a\) colors,
one distinguished, and delete sources with other colors.
No value of a map is a source, so each deleted point maps to a retained
point. The remaining map is arbitrary, and its functional digraph has the
same number of connected components as before. Conversely, attach to each of its vertices a set of new
predecessors, each with one of \(a-1\) colors; any source of the resulting map among
the retained vertices has the distinguished color.
Thus each retained vertex is replaced by an object with generating
function \(z\exp((a-1)z)\), proving the substitution in
\eqref{eq:sources-components}. The coefficients are polynomials in
\(a\), so the identity holds for an indeterminate \(a\).

For the finite sum, the coefficient substitution \(u=t\exp(-t)\)
gives, for \(k\geq0\),
\[
 k!\coef{u^k}(1-\mathcal T(u))^{-b}
 =k!\coef{t^k}\exp(kt)(1-t)^{1-b}
 =\sum_{j=0}^k\binom kj(b-1)^{\overline j}k^{k-j}.
\]
In~\eqref{eq:sources-components}, the contribution of \(u^k/k!\)
under \(u=z\exp((a-1)z)\) is
\[
 m!\coef{z^m}\frac{z^k\exp(k(a-1)z)}{k!}
 =\binom mk\bigl(k(a-1)\bigr)^{m-k}
 \qquad(0\leq k\leq m).
\]
Multiplying and combining the powers of \(k\) proves the finite sum.

For the elementary coefficient formula, take positive integers
\(a,b\), and put \(U(z):=\mathcal T(z\exp((a-1)z))\) and
\(Y(z):=U(z)+(a-1)z\). Then
\[
 Y=z\bigl(a-1+\exp(Y)\bigr).
\]
Since \(a\ne0\), the coefficient substitution
\(z=t/(a-1+\exp(t))\) is invertible at zero. It gives
\[
 U=\frac{t\exp(t)}{a-1+\exp(t)},\qquad
 \frac{dz}{z}=
 \left(1-\frac{t\exp(t)}{a-1+\exp(t)}\right)\frac{dt}{t}.
\]
Consequently,
\[
 \coef{z^m}(1-U(z))^{-b}
 =\coef{t^m}(a-1+\exp(t))^m
    \left(1-\frac{t\exp(t)}{a-1+\exp(t)}\right)^{1-b},
\]
which proves~\eqref{eq:sources-components-elementary}.
\end{proof}

The counts for \(fg=g^3\) and \(fg=fg^2\) now differ only in
the weight assigned to each connected component. In the coefficient
formulas below, \(n\) is held fixed during coefficient extraction.

\begin{corollary}
\label{cor:source-identities}
For \(n\geq0\),
\begin{align}
 \NN_{fg,ggg}(n)
   &=n!\coef{t^n}\bigl(n-1+\exp(t)\bigr)^n
     =\sum_{k=0}^n\binom nk k^n(n-1)^{n-k},
       \label{eq:total750}\\
 \NN_{fg,fgg}(n)
   &=n!\coef{t^n}\left(
     1+\frac{(n-1+\exp(t))^{2n-1}}
          {(n-1+(1-t)\exp(t))^{n-1}}\right).
       \label{eq:companion785-total}
\end{align}
\end{corollary}

\begin{proof}
For \(n\geq1\) and fixed \(g\), equation~\eqref{eq:local750} gives
\(n^{n-|g(X_n)|}=n^{s(g)}\) choices of \(f\). Substitution
\((a,b)=(n,1)\) in~\eqref{eq:sources-components-elementary}
gives the first coefficient formula. Expanding its \(n\)th power
by the binomial theorem and using
\(n!\coef{t^n}\exp(kt)=k^n\) gives the finite sum.
For \(fg=fg^2\), equation~\eqref{eq:source-component-fixed}
gives \(n^{s(g)+c(g)}\) choices. Substitution \((a,b)=(n,n)\)
gives the quotient in~\eqref{eq:companion785-total}; adding one does
not change its coefficient of positive degree. At \(n=0\), both
expressions in~\eqref{eq:total750} equal one, while the expression
in parentheses in~\eqref{eq:companion785-total} simplifies to
\[
 1+\frac{(1-t)\exp(t)-1}{\exp(t)-1}
 =\frac{(2-t)\exp(t)-2}{\exp(t)-1}
 =1+O(t).
\]
The common factor \(t\) cancels, so this is a formal power series
with constant coefficient one.
\end{proof}

The finite sum in~\eqref{eq:total750} is recorded by Heinz in
\oeis{A239750}. The \((a,b)=(n,n)\) specialization
of~\eqref{eq:sources-components} is the formula
of~\cite[Theorem~2.1]{AlekseyevCompanion}.

Specializing the finite sum in \namedref{Lemma}{lem:sources-components}
at \((a,b)=(n,n)\) gives, for \(n\geq0\),
\begin{equation}
 \NN_{fg,fgg}(n)
 =\sum_{k=0}^n\binom nk(n-1)^{n-k}
    \sum_{j=0}^k\binom kj(n-1)^{\overline j}k^{n-j}.
 \label{eq:companion785-sum}
\end{equation}
Every summand is a nonnegative integer.

\subsection{Fixed points and their preimages}
\label{sec:fixed-point-totals}

For the next two identities, we choose the fixed points and
the points mapped to them. This gives the formulas recorded by Heinz
and Einstein in \oeis{A239781} and \oeis{A239759}, respectively.

\begin{theorem}
\label{thm:fixed-point-totals}
For \(n\geq0\),
\begin{align}
 \NN_{fg,ffg}(n)
  &=\sum_{r=0}^n\binom nr r^n
      \sum_{k=0}^{r}\binom rk
          k^{r-k}(n-k-1)^{n-r},
       \label{eq:total781}\\
 \NN_{fg,gfg}(n)
  &=\sum_{k=0}^n\binom nk k^k
      \sum_{j=0}^{n-k}\binom{n-k}{j}
          k^j(k+j-1)^j
          \bigl((n-k)(k+j)\bigr)^{n-k-j}.
       \label{eq:total759}
\end{align}
\end{theorem}

\begin{proof}
For \(fg=f^2g\), choose the \(r\) points mapped to fixed points
of \(f\), and then the \(k\) fixed points among them. The other
\(r-k\) points have \(k\) possible images each. Every point outside
the chosen \(r\)-element set has \(n-k-1\) possible images,
excluding itself and the fixed points. By~\eqref{eq:local781}, the allowed
values of \(g\) are exactly the chosen \(r\) points. Thus there
are \(r^n\) choices for \(g\), proving~\eqref{eq:total781}.

For \(fg=gfg\), choose the fixed-point set \(F:=\Fix(g)\),
of size \(k\). Put \(S:=\{x\in X_n:f(x)\notin F\}\). The identity
requires \(f(g(X_n))\subseteq F\), so every point of \(S\)
must be a source of \(G(g)\). In particular, \(S\) is disjoint
from \(F\). Suppose \(|S|=n-k-j\), leaving \(j\)
nonfixed points outside \(S\). The values of \(f\) have
\(k^{k+j}(n-k)^{n-k-j}\) choices. The map \(g\) fixes \(F\)
pointwise; each point of \(S\) has \(k+j\) possible images;
and each of the other \(j\) nonfixed points has \(k+j-1\)
possible images, since its own label is excluded. These choices give
exactly the pairs with the chosen \(F,S\). Summing over their choices
proves~\eqref{eq:total759}.
\end{proof}

\subsection{Fixing the second image pointwise}
\label{sec:second-image-total}

For \(f^2=gf^2\), each fixed \(f\) allows
\(n^{n-|f^2(X_n)|}\) choices of \(g\). Partitioning the domain
using the fixed points of \(g\) gives a direct finite sum and a
generating function recording the number of points outside the second image.
Einstein's earlier formula in \oeis{A239752} groups maps by successive
image sizes.

\begin{theorem}
\label{thm:total752}
For an indeterminate \(u\),
\[
 \sum_{m\geq0}\left(\sum_{h\in\MM_m}
       u^{m-|h^2(X_m)|}\right)\frac{z^m}{m!}
 =\left(1-\mathcal T\!\left(
       z\exp\bigl((u-1)z\exp(uz)\bigr)\right)\right)^{-1}.
\]
Consequently, for \(n\geq0\),
\begin{align}
 \NN_{ff,gff}(n)
 &=n!\coef{z^n}
   \left(1-\mathcal T\!\left(
       z\exp\bigl((n-1)z\exp(nz)\bigr)\right)\right)^{-1}
       \label{eq:total752-series}\\
 &=\sum_{k=0}^n\binom nk k^k
     \sum_{j=0}^{n-k}\binom{n-k}{j}
           \bigl(k(n-1)\bigr)^j(nj)^{n-k-j}.
       \label{eq:total752}
\end{align}
\end{theorem}

\begin{proof}
For a pair \((f,g)\in\MM_n^2\) satisfying \(f^2=gf^2\), put
\(F:=\Fix(g)\) and \(H:=f^{-1}(F)\). The identity says
\(f(X_n)\subseteq H\). Partition \(X_n\) into
\[
 K:=F\cap H,\qquad J:=H\setminus F,\qquad L:=X_n\setminus H.
\]
Then
\[
 f(K\cup J)\subseteq K,\qquad f(L)\subseteq J;
\]
the map \(g\) fixes every point of \(K\), moves every point of
\(J\), and is unrestricted on \(L\). Conversely, these conditions
give \(f^2(X_n)\subseteq K\subseteq\Fix(g)\) and recover the
same partition from \(f,g\).

Write \(k:=|K|\), \(j:=|J|\), and \(\ell:=|L|=n-k-j\).
For \(n\geq1\), the set \(K\) is nonempty. Choosing the three
sets and then the values of the maps gives
\[
 \binom nk\binom{n-k}{j}\,k^k
       \bigl(k(n-1)\bigr)^j(nj)^\ell.
\]
Indeed, \(f|_K\) is arbitrary, each point of \(J\) has \(k\)
choices under \(f\) and \(n-1\) under \(g\), and each point
of \(L\) has \(j\) choices under \(f\) and \(n\) under
\(g\). Summing proves~\eqref{eq:total752}.

For the weighted generating function, first take an integer \(u\geq1\).
For each integer \(m\geq0\) and map \(h\in\MM_m\), color the
points of \(X_m\) with \(u\) colors, one distinguished,
requiring every point of \(h^2(X_m)\) to have that color.
There are \(u^{m-|h^2(X_m)|}\) such colorings. Take \(F\) to be
the set with the distinguished color and use the same partition with
\(H:=h^{-1}(F)\). The restriction to \(K\) is an arbitrary map.
Each of its vertices receives a set of rooted stars: the root lies in
\(J\), with \(u-1\) color choices, and its leaves lie in \(L\),
with \(u\) choices each. Such a star has series
\((u-1)z\exp(uz)\). Replacing each vertex of an arbitrary map
by a vertex with a set of these stars gives
\[
 \sum_{k\geq0}\frac{k^k}{k!}
   \left(z\exp\bigl((u-1)z\exp(uz)\bigr)\right)^k
 =\left(1-\mathcal T\!\left(
   z\exp\bigl((u-1)z\exp(uz)\bigr)\right)\right)^{-1}.
\]
Its coefficients are polynomials in \(u\), so the identity holds
for an indeterminate \(u\). Setting \(u=n\) gives
\eqref{eq:total752-series}.
\end{proof}

The rooted-star weights in~\eqref{eq:total752} also enter the
classical inverse relations for Abel polynomials; see~\cite{Sagan1983}.

\subsection{An exact image and inclusion--exclusion}
\label{sec:image}

We count pairs of maps satisfying \(f^2g=g^2\), the orientation
used in \oeis{A239756}. Exchanging the map names gives the normalized
identity \(f^2=g^2f\). Let
\[
 E_n:=\#\{(f,g)\in\MM_n^2:f^2g=g^2\}.
\]

\begin{theorem}\label{thm:E}
For nonnegative integers \(a,t\), put
\[
 P_{a,t}:=\sum_{i=0}^a\binom{a}{i}i^i
                    (a^2-i)^t t^{a-i}.
\]
For \(n\geq0\), put \(m:=n-a\) in the summand indexed by \(a\). Then
\begin{equation}\label{eq:E-positive}
 E_n=\sum_{a=0}^n\binom na
       \sum_{t=0}^{m}\binom mt P_{a,t}(am-t)^{m-t}.
\end{equation}
\end{theorem}

\begin{proof}
Fix \(R\subseteq X_n\) of size \(k\). If \(\im g=R\), every
\(y\in R\) has the form \(g(x)\) for some \(x\in X_n\), so the equation gives
\(g(y)=f^2(y)\). Conversely, \(g|_R=f^2|_R\) implies
\(f^2g=g^2\). For fixed \(A\subseteq R\), with \(|A|=a\), define
\[
 C(n,k,a):=\#\{f\in\MM_n:f^2(R)\subseteq A\}.
\]
Relabeling shows that this depends only on \(n,k,a\).
For fixed \(f\in\MM_n\), inclusion--exclusion over values in
\(R\) missed by \(g\) gives the number of maps \(g\in\MM_n\)
with image exactly \(R\) and \(g|_R=f^2|_R\):
\[
 \sum_{A\subseteq R}(-1)^{k-|A|}
 \mathbf1_{\{f^2(R)\subseteq A\}}|A|^{n-k}.
\]
Indeed, the forced values on \(R\) lie in the allowed set \(A\)
exactly when \(f^2(R)\subseteq A\), and the other \(n-k\)
arguments each have \(|A|\) choices. Summing over \(f\), the
\(\binom nk\) choices of \(R\), and the \(\binom ka\) choices of
\(A\) gives
\begin{equation}\label{eq:E}
 E_n=\sum_{k=0}^n\binom nk
       \sum_{a=0}^k(-1)^{k-a}\binom ka
       a^{n-k}C(n,k,a).
\end{equation}

For each \(a\geq0\), define
\[
 \mathcal C_a(x,y):=
 \sum_{b,c\geq0}C(a+b+c,a+b,a)\frac{x^b}{b!}\frac{y^c}{c!},
\]
where \(x\) records the points of \(R\setminus A\), and \(y\)
the points outside \(R\). Writing \(b:=k-a\) and \(c:=n-k\),
the factor \((-1)^b a^c\) in~\eqref{eq:E} gives \(x=-z\),
\(y=az\), while \(\binom nk\binom ka=n!/(a!b!c!)\). Hence
\[
 \sum_{n\geq0}E_n\frac{z^n}{n!}
   =\sum_{a\geq0}\frac{z^a}{a!}\mathcal C_a(-z,az).
\]
The factor \(z^a\) makes this a well-defined formal series.

To evaluate \(\mathcal C_a\), put \(S:=f^{-1}(A)\). The condition
\(f^2(R)\subseteq A\) is equivalent to \(f(R)\subseteq S\). Set
\[
 i:=|S\cap A|,\qquad
 \beta:=|S\cap(R\setminus A)|,\qquad
 \gamma:=|S\setminus R|.
\]
Choosing these three parts of \(S\) gives
\(\binom{a}{i}\binom{k-a}\beta\binom{n-k}\gamma\).
Points in \(R\cap S\) map into \(A\cap S\), giving
\(i^{i+\beta}\) choices; points in \(R\setminus S\)
map into \(S\setminus A\), giving
\((\beta+\gamma)^{k-i-\beta}\). Points in \(S\setminus R\)
map into \(A\), and points outside \(R\cup S\) map outside \(A\),
giving \(a^\gamma(n-a)^{n-k-\gamma}\). These choices give exactly
the maps with inverse image \(S\), so
\begin{align}
 C(n,k,a)
 ={}&\sum_{i=0}^a\sum_{\beta=0}^{k-a}\sum_{\gamma=0}^{n-k}
 \binom{a}{i}\binom{k-a}\beta\binom{n-k}\gamma\notag\\[-2pt]
 &\quad{}\cdot i^{i+\beta}
 (\beta+\gamma)^{k-i-\beta}
 a^\gamma(n-a)^{n-k-\gamma}.\label{eq:C}
\end{align}

We now sum~\eqref{eq:C} with the factorial weights in
\(\mathcal C_a(x,y)\). Put \(\tau:=\mathcal T(y)\), \(\xi:=xe^\tau\),
\(u:=b-\beta\), and \(v:=c-\gamma\).
For a nonnegative integer \(q\),
\namedref{Lemma}{lem:abel-tree-coefficients}, with its parameters set to
\(b=1,c=q,r=0\), gives
\[
 \sum_{v\geq0}(q+v)^v\frac{y^v}{v!}
   =\frac{\exp(q\mathcal T(y))}{1-\mathcal T(y)}.
\]
Apply this with \(q=\beta+\gamma+u\) to sum over \(v\) in
\eqref{eq:C}. For fixed \(a,i\), after taking out
\(\binom{a}{i}i^i/(1-\tau)\), the remaining series is
\[
 \sum_{\beta,\gamma,u\geq0}
  \frac{(i\xi)^\beta(a\tau)^\gamma\xi^u}
       {\beta!\gamma!u!}
  (\beta+\gamma)^{a-i+u}.
\]
Summing over \(u\) gives \(\exp(\xi(\beta+\gamma))\).
Group by \(t:=\beta+\gamma\), and put
\(Z_i:=(i\xi+a\tau)e^\xi\).
To evaluate the remaining power sum, for nonnegative integers
\(p,q\), let \(S(p,q)\) denote the Stirling number of the second kind, with
\(S(0,0):=1\) and \(S(p,0):=0\) for \(p>0\). Then
\[
\begin{aligned}
 \mathcal C_a(x,y)
 &=\frac{1}{1-\tau}
   \sum_{i=0}^a\binom{a}{i}i^i
   \sum_{t\geq0}t^{a-i}\frac{Z_i^t}{t!}\\
 &=\frac{1}{1-\tau}
   \sum_{i=0}^a\binom{a}{i}i^i e^{Z_i}
   \sum_{j=0}^{a-i}S(a-i,j)Z_i^j.
\end{aligned}
\]

For \(a\geq0\), let \(W\) be the unique formal series with zero
constant term satisfying \(W=z\exp(aW)\). Then
\(\mathcal T(az)=aW\), and the substitution \((x,y)=(-z,az)\)
sends \(\xi\) and \(Z_i\) to \(-W\) and
\((a^2-i)W\exp(-W)\), respectively. The power-sum form above gives
\[
 \mathcal C_a(-z,az)=\frac{1}{1-aW}
 \sum_{i=0}^a\binom{a}{i}i^i
 \sum_{t\geq0}\frac{t^{a-i}(a^2-i)^t}{t!}
       W^t\exp(-tW).
\]
Its contribution to \(E_n\) is
\(\binom na m!\coef{z^m}\mathcal C_a(-z,az)\), where \(m=n-a\).
\namedref{Lemma}{lem:abel-tree-coefficients}, with its parameters set to
\(b=a,c=-t,r=t\), gives
\[
 \coef{z^m}\frac{W^t\exp(-tW)}{1-aW}
 =\frac{(am-t)^{m-t}}{(m-t)!}\qquad(0\leq t\leq m),
\]
and terms with \(t>m\) vanish. Substitution and summation over
\(i\) give~\eqref{eq:E-positive}. For \(a\geq1\), its summands are nonnegative
because \(a^2-i\geq0\) and \(am-t\geq0\). For \(a=0\),
we have \(i=0\); every term with \(t>0\) vanishes, and the
\(t=0\) term is \(0^n\). Thus this index contributes one at
\(n=0\) and zero otherwise.
\end{proof}

Now fix \(f\in\mathfrak S_n\). For a map \(g\in\MM_n\)
satisfying \(f^2g=g^2\), its image \(R:=g(X_n)\) must satisfy
\(f^2(R)=R\), and \(g|_R=f^2|_R\) already maps \(R\) onto itself.
Consequently every extension from \(X_n\setminus R\) into \(R\)
has image exactly \(R\), and inclusion--exclusion is unnecessary.
The count for this fixed \(f\) is
\begin{equation}
 \#\{g:f^2g=g^2\}
 =\sum_{k=0}^n k^{n-k}\coef{t^k}
       \prod_{C\in\Comp(f^2)}(1+t^{|C|}).
 \label{eq:permutation-exact-image}
\end{equation}
The coefficient of \(t^k\) is the number of unions of cycles
of \(G(f^2)\) with \(k\) vertices.

For unrestricted \(f\), extensions into \(R\) can have a smaller
image. Formula~\eqref{eq:E} enforces surjectivity onto \(R\)
without enumerating \(g\).

\section{Allocations of fibers}
\label{sec:fibers}

\subsection{Fiber labels and weighted allocations}

To specify a map \(f\in\MM_n\), we choose its fiber partition of \(X_n\)
and then its image labels. List the blocks as \(B_1,\ldots,B_k\), with
\(s_i:=|B_i|\) for \(i\in X_k\). The image labels
\(y_1,\ldots,y_k\in X_n\) must be distinct and satisfy
\[
 f(B_i)=\{y_i\},\qquad 1\leq i\leq k.
\]
The locations of these labels within the same partition define a map
\(\tau:X_k\to X_k\) by \(y_j\in B_{\tau(j)}\) for \(j\in X_k\).
We call \(B_{\tau(j)}\) the successor of \(B_j\). For each target index
\(i\in X_k\), record the number and total size of the incoming blocks by
\[
 d_i:=|\tau^{-1}(i)|,\qquad
 b_i:=\sum_{\tau(j)=i}s_j.
\]
For nonnegative integers \(s,d\), we use the falling factorial
\(\fall{s}{d}:=s(s-1)\cdots(s-d+1)\), with
\(\fall{s}{0}:=1\) and \(\fall{s}{d}:=0\) for \(d>s\).

Equivalently, let \(\pi:X_n\to X_k\) send a point to its fiber
index, and let \(r:X_k\to X_n\) be the injective map \(r(j):=y_j\).
Then
\[
 f=r\pi,\qquad \tau=\pi r,\qquad
 f^m=r\tau^{m-1}\pi\quad(m\geq1).
\]
\Needspace{12\baselineskip}
\begin{lemma}\label{lem:labels}
For the fixed partition and a fixed \(\tau\), the number of
possible assignments of distinct image labels is
\[
 \prod_{i=1}^k\fall{s_i}{d_i}.
\]
For every resulting map and every \(x\in B_j\),
\begin{equation}\label{eq:iterates}
 f^m(x)=y_{\tau^{m-1}(j)}\qquad(m\geq1).
\end{equation}
In particular, \(|(f^2)^{-1}(y_i)|=b_i\), including the possibility \(b_i=0\).
\end{lemma}

\begin{proof}
The \(d_i\) distinguishable blocks sent by \(\tau\) to \(i\) must receive
distinct labels from \(B_i\). This gives \(\fall{s_i}{d_i}\) assignments,
independently for each \(i\). Since the blocks are disjoint, these assignments
make all \(y_j\) distinct. The first iterate formula is the definition of
\(y_j\); applying \(f\) replaces \(y_j\) by \(y_{\tau(j)}\), proving
\eqref{eq:iterates} by induction. The inverse image of \(y_i\) under \(f^2\) is the
disjoint union of the blocks \(B_j\) with \(\tau(j)=i\).
\end{proof}

In numerical examples we label the underlying set
\(\{0,\ldots,n-1\}\) and write a map as the tuple of its values in this order.
Fiber sizes alone do not determine those of the second iterate:
on \(\{0,1,2\}\), the maps \((0,0,1)\) and
\((0,0,2)\) both have fiber sizes \(2,1\), but their second iterates have
fiber-size lists \((3)\) and \((2,1)\), respectively.

Fix \(f\in\MM_n\). For \(g\in\MM_n\), the identity \(f^2=fgf\)
is equivalent to \((fg)|_R=f|_R\), where \(R:=f(X_n)\).
Formula~\eqref{eq:prescribed-on-subset} therefore gives
\begin{equation}
 \#\{g:f^2=fgf\}
   =n^{n-|R|}\prod_{y\in R}|f^{-1}(f(y))|.
 \label{eq:companion-fgf-fixed}
\end{equation}

To sum over fiber partitions, we use the following multinomial count of
weighted allocations between blocks. For independent indeterminates
\(\mathbf z=(z_1,\ldots,z_n)\) and a vector
\(\mathbf m=(m_1,\ldots,m_n)\) of nonnegative integers, write
\(\coef{\mathbf z^{\mathbf m}}\) for extraction of the coefficient
of \(z_1^{m_1}\cdots z_n^{m_n}\).

\begin{lemma}
\label{lem:weighted-block-allocation}
Fix \(n\geq0\). For \(1\leq i\leq n\) and each vector
\(\mathbf e=(e_1,\ldots,e_n)\) of nonnegative integers with
\(\sum_jje_j\leq n\), let \(W_i(\mathbf e)\) be a nonnegative integer
weight, and define
\begin{equation}
 \Phi_i(\mathbf z):=
 \sum_{\substack{\mathbf e\geq\mathbf0\\\sum_jje_j\leq n}}
 W_i(\mathbf e)\prod_{j=1}^n\frac{z_j^{e_j}}{e_j!}.
 \label{eq:allocation-polynomial}
\end{equation}
Let \(\mathcal B\) be a partition of \(X_n\), with \(m_j\)
blocks of size \(j\). We call \(\mathbf m\) its \emph{profile}.
For a map \(\lambda\colon\mathcal B\to\mathcal B\)
and a target block \(B\in\mathcal B\), put
\[
 e_{B,j}(\lambda):=
   \#\{C\in\mathcal B:|C|=j,\ \lambda(C)=B\}.
\]
Here \(1\leq j\leq n\), and
\(\mathbf e_B(\lambda):=(e_{B,1}(\lambda),\ldots,e_{B,n}(\lambda))\).
Then
\begin{equation}
 \sum_{\lambda\colon\mathcal B\to\mathcal B}
       \prod_{B\in\mathcal B}W_{|B|}(\mathbf e_B(\lambda))
 =\left(\prod_jm_j!\right)
   \coef{\mathbf z^{\mathbf m}}\prod_i\Phi_i(\mathbf z)^{m_i}.
 \label{eq:allocation-fixed-partition}
\end{equation}
Consequently, summing this weight over all partitions \(\mathcal B\)
of \(X_n\) gives
\begin{equation}
 n!\sum_{\substack{\mathbf m\geq\mathbf0\\\sum_jjm_j=n}}
 \coef{\mathbf z^{\mathbf m}}
 \prod_i\left(\frac{\Phi_i(\mathbf z)}{i!}\right)^{m_i}.
 \label{eq:allocation-all-partitions}
\end{equation}
\end{lemma}

\begin{proof}
Fix the numbers \(e_{B,j}\) for every target block \(B\) and source
size \(j\). They arise from an allocation exactly when
\(\sum_Be_{B,j}=m_j\) for each \(j\). The size-\(j\)
source blocks can then be assigned in
\(m_j!/\prod_Be_{B,j}!\) ways. Multiplication over \(j\), with
the factors \(W_{|B|}(\mathbf e_B)\), is exactly the coefficient
on the right of~\eqref{eq:allocation-fixed-partition}.

There are
\[
 \frac{n!}{\prod_jm_j!(j!)^{m_j}}
\]
partitions of profile \(\mathbf m\). Multiplying by this number
cancels \(\prod_jm_j!\), and summing over the profiles gives
\eqref{eq:allocation-all-partitions}. No factorial for permuting target
blocks is introduced: they are already distinguished by their elements.
\end{proof}

Let a target fiber have size \(b\), and let \(e_j\) be the number of
source fibers of size \(j\) assigned to it, for \(1\leq j\leq n\).
Write \(k:=\sum_je_j\) for their number and
\(w:=\sum_jje_j\) for their total size. In the last row below, \(c\)
is the size of the successor of the target fiber. The local weights are
\[
\begin{array}{c|c}
\text{identity}&\text{local weight}\\\hline
f^2=fgf& n^{b-1}\fall{b}{k}b^k\\
f^2=f^2g&\fall{b}{k}w^w\\
fg=fgf&\fall{w}{k}b^w\\
fg=f^3&\fall{b}{k}c^w
\end{array}
\]
For the first row, if \(f\) has \(r\) fibers, the product of
\(n^{b-1}\) over its target fibers is \(n^{n-r}\), accounting for
the free values in~\eqref{eq:companion-fgf-fixed}. For the second row,
the image labels lie in a target fiber of size \(b\), and \(g\)
maps the corresponding inverse image under \(f^2\), of size \(w\),
into itself. In the third row, the image
labels lie in the source union of size \(w\), and \(g\) maps that
union into the size-\(b\) fiber. These two rows exchange the roles
of \(b\) and \(w\) in the label and map choices.
The final row also depends on the successor size \(c\), so we group
target fibers by the pair \((b,c)\). The same allocation argument applies
separately to each prescribed successor size.

Requiring \(f\) to be a permutation restricts every fiber-allocation
sum to singleton fibers. Distinct image labels then make
\(\tau\) a permutation. For a permutation \(g\), three of the
allocation recurrences use the following local weights:
\[
\begin{array}{c|c}
 \text{identity}&\text{local weight for permutation }g\\\hline
 f^2=fgf&\bigl(\fall{b}{k}\bigr)^2\\
 f^2=f^2g&\fall{b}{k}\,w!\\
 fg=fgf&\fall{w}{k}\,b!\,\mathbf1_{\{w=b\}}.
\end{array}
\]
For the first row only, multiply the contribution of a profile with
\(r:=\sum_i m_i\) fibers by \((n-r)!\). Indeed, the \(r\)
prescribed arguments of \(g\) must have distinct images in their
required fibers, and the other values complete a permutation. In the
second row, \(g\) permutes each second-iterate fiber. In the third,
it must biject each source union with its target fiber, which requires
\(w=b\). The allocation of image labels is otherwise unchanged.

The coefficient extraction in \namedref{Lemma}{lem:weighted-block-allocation}
has an integer recurrence for direct evaluation. For a fixed profile \(\mathbf m\), put
\(r:=\sum_jm_j\) and choose a list \(s_1,\ldots,s_r\) containing
\(m_j\) copies of \(j\). Let \(\mathbf a=(a_1,\ldots,a_n)\), where
\(0\leq a_j\leq m_j\), record the numbers of source blocks of each size
still available. Define
\[
 V_{r+1}(\mathbf a):=\mathbf1_{\{\mathbf a=\mathbf0\}},
\]
and, for \(t=r,r-1,\ldots,1\), set
\begin{equation}
 V_t(\mathbf a):=
 \sum_{\mathbf0\leq\mathbf e\leq\mathbf a}
 \left(\prod_j\binom{a_j}{e_j}\right)
 W_{s_t}(\mathbf e)V_{t+1}(\mathbf a-\mathbf e).
 \label{eq:shared-allocation-recurrence}
\end{equation}
The vector inequalities are coordinatewise. The binomial factors choose
the source blocks assigned to the current target, so \(V_1(\mathbf m)\)
is the weighted sum for one partition. Multiplying by the number of
partitions of this profile gives
\[
 \frac{n!}{\prod_jm_j!(j!)^{m_j}}V_1(\mathbf m).
\]

\subsection{A prescribed value after one fiber choice}
\label{sec:companion-fgf-total}

Specializing \namedref{Lemma}{lem:weighted-block-allocation} to weights
depending on a target fiber's size and the number of image labels it
receives recovers the fiber-profile calculation for \(f^2=fgf\)
in~\cite[Theorem~2.2]{AlekseyevCompanion}. Eliminating its auxiliary
constant term gives the following formula.

For fixed \(n\geq0\), define
\[
 R_n(t):=\frac{1-t-\sqrt{(1+t)^2-4t\exp(nt)}}2,
\]
where the square root has constant term one. Equivalently, \(R_n\)
is the unique formal series with zero constant term satisfying
\begin{equation}
 R_n=t\bigl(\exp(nt)-1\bigr)+tR_n+R_n^2.
 \label{eq:companion772-quadratic}
\end{equation}
In particular, \(R_n\) has order at least two in \(t\).

\begin{theorem}
\label{thm:companion772-explicit}
For \(n\geq0\), with \(n\) held fixed during coefficient extraction,
\begin{equation}
 \NN_{ff,fgf}(n)
 =n!\coef{t^n}
   \frac{(1-R_n(t))^2\exp(nR_n(t))}{1-t-2R_n(t)}.
 \label{eq:companion772-explicit}
\end{equation}
\end{theorem}

\begin{proof}
For a profile \(\mathbf m=(m_1,\ldots,m_n)\) with
\(\sum_i im_i=n\), put \(r:=\sum_i m_i\). By
\eqref{eq:companion-fgf-fixed}, a target fiber of size \(i\)
receiving \(k\) image labels has local weight
\(n^{i-1}\fall{i}{k}i^k\). Its allocation polynomial is
\[
 \Phi_i(\mathbf z)=n^{i-1}
                  \left(1+i\sum_jz_j\right)^i,
\]
with terms beyond the required total size omitted. For any formal power
series \(F(u)\),
\[
 \coef{\mathbf z^{\mathbf m}}F\!\left(\sum_jz_j\right)
 =\frac{r!}{\prod_jm_j!}\coef{u^r}F(u),
\]
so \namedref{Lemma}{lem:weighted-block-allocation} gives the profile contribution
\begin{equation}
 \frac{n!\,n^{n-r}r!}{\prod_i m_i!(i!)^{m_i}}
       \coef{u^r}\prod_i(1+iu)^{im_i}.
 \label{eq:companion772-profile}
\end{equation}

Put \(Q_n(y,z):=\sum_{i=1}^n y^i(n+iz)^i/i!\).
In the expansion \(\sum_{r\geq0}z^{-r}Q_n(y,z)^r\), a profile
\(\mathbf m\) with \(\sum_i im_i=n\) contributes
\[
 \frac{r!}{\prod_i m_i!(i!)^{m_i}}
       \coef{z^r}\prod_i(n+iz)^{im_i}
\]
to \(\coef{z^0y^n}\). Extracting degree \(r\) from the last product gives
\(n^{n-r}\coef{u^r}\prod_i(1+iu)^{im_i}\), so summing
\eqref{eq:companion772-profile} gives
\begin{equation}
 \NN_{ff,fgf}(n)=n!\coef{z^0 y^n}
                    \left(1-\frac{Q_n(y,z)}z\right)^{-1}.
 \label{eq:companion772-total}
\end{equation}
The inverse is formed in \(\mathbb Q[z,z^{-1}][[y]]\).
Only \(r\leq n\) and \(i\leq n\) contribute, so extracting
\(\coef{y^n}\) leaves a Laurent polynomial whose constant term in
\(z\) is taken.

Let \(V=V(y,z)\) satisfy \(V=y\exp(zV)\) with zero constant
term in \(y\). \namedref{Lemma}{lem:abel-tree-coefficients} gives
\[
 1+Q_n(y,z)\equiv\frac{\exp(nV)}{1-zV}\pmod{y^{n+1}}.
\]
Thus extending the sum defining \(Q_n\) to all positive indices in
\eqref{eq:companion772-total} does not affect the desired coefficient.
The coefficient substitution \(y=t\exp(-zt)\) replaces \(V\)
by \(t\) and contributes the factor \((1-zt)\exp(nzt)\), giving
\begin{align*}
 \frac{\NN_{ff,fgf}(n)}{n!}
 &=\coef{z^0 t^n}
   \frac{z(1-zt)^2\exp(nzt)}{(z+1)(1-zt)-\exp(nt)}\\
 &=\coef{z^0 t^n}
   \frac{(1-zt)^2\exp(nzt)}
        {1-t-tz-(\exp(nt)-1)/z}.
\end{align*}
Each coefficient of \(t\) is a Laurent polynomial in \(z\),
so the two extractions commute.

Write \(R:=R_n(t)\) and \(B:=1-t-R\).
Equation~\eqref{eq:companion772-quadratic} gives
\(BR=t(\exp(nt)-1)\), and hence
\[
 1-t-tz-\frac{\exp(nt)-1}{z}
 =(B-tz)\left(1-\frac{R}{tz}\right).
\]
The inverse factors have valid geometric expansions because \(B\)
has constant term one and \(R/t\) has zero constant term.
Put \(H(w):=(1-w)^2\exp(nw)=\sum_{j\geq0}h_jw^j\).
Then
\begin{align*}
 \coef{z^0}\frac{H(tz)}{(B-tz)(1-R/(tz))}
 &=\sum_{j,p\geq0}h_jR^{j+p}B^{-p-1}\\
 &=\frac{H(R)}{B-R}
 =\frac{(1-R)^2\exp(nR)}{1-t-2R},
\end{align*}
which proves~\eqref{eq:companion772-explicit}. At \(n=0\), the specified
square root is \(1-t\), so \(R_0(t)=0\) and the right-hand side
is \(\coef{t^0}(1-t)^{-1}=1\).
\end{proof}

A different elimination of the auxiliary constant term is given
in~\cite{AlekseyevCompanion}.

\subsection{Fibers of the second iterate}
\label{sec:profile}

For \(f^2=f^2g\), image labels are allocated in the same way, while
the choices of \(g\) depend on the size of each second-iterate fiber.
Write
\[
 P_n:=\#\{(f,g)\in\MM_n^2:f^2=f^2g\}.
\]
For a fixed \(f\in\MM_n\), the relation says that \(g\in\MM_n\)
preserves every fiber of \(f^2\). Taking \(a=b=f^2\) in
\namedref{Lemma}{lem:prescribed-compositions}, and putting
\(b_y:=|(f^2)^{-1}(y)|\) for \(y\in\im(f^2)\), gives
\begin{equation}\label{eq:fixedP}
 \#\{g:f^2=f^2g\}=\prod_{y\in\im(f^2)}b_y^{b_y}.
\end{equation}
We sum~\eqref{eq:fixedP} by fiber profile.

Define \(w(0):=1\) and \(w(b):=b^b\) for positive integers \(b\).
With the notation of \namedref{Lemma}{lem:labels}, the combined weight for a
fixed block map \(\tau\), including its image-label assignments, is
\begin{equation}\label{eq:blockP}
 \prod_{i=1}^k\fall{s_i}{d_i}\,w(b_i).
\end{equation}
This weight is unchanged under simultaneous relabeling of the blocks.

For a profile \(\mathbf m=(m_1,\ldots,m_n)\), \(m_j\) is the number
of fibers of size \(j\), and \(\sum_jjm_j=n\). For
\(1\leq i\leq n\), define the finite polynomial
\begin{equation}\label{eq:Psi}
 \Psi_i(\mathbf z):=
 \sum_{\substack{e_1,\ldots,e_n\geq0\\e_1+\cdots+e_n\leq i}}
 \fall{i}{e_1+\cdots+e_n}
 w\!\left(\sum_{j=1}^n je_j\right)
 \prod_{j=1}^n\frac{z_j^{e_j}}{e_j!}.
\end{equation}
Here \(e_j\) records how many source blocks of size \(j\) place their
image labels in a particular target block of size \(i\).

\begin{theorem}\label{thm:P}
For \(n\geq0\),
\begin{equation}\label{eq:P}
 P_n=n!\sum_{\substack{m_1,\ldots,m_n\geq0\\\sum_{j=1}^n jm_j=n}}
 \coef{\mathbf z^{\mathbf m}}
 \prod_{i=1}^n\left(\frac{\Psi_i(\mathbf z)}{i!}\right)^{m_i}.
\end{equation}
\end{theorem}

\begin{proof}
By~\eqref{eq:blockP}, the local weights are
\(W_i(\mathbf e):=\fall{i}{\sum_je_j}w(\sum_jje_j)\).
Their allocation polynomial is \(\Psi_i\) after terms of total
source size greater than \(n\) are discarded; those terms cannot
contribute to the required coefficient.
Formula~\eqref{eq:allocation-all-partitions} now gives~\eqref{eq:P}.
\end{proof}

With these weights, recurrence~\eqref{eq:shared-allocation-recurrence}
evaluates each profile's contribution using integer arithmetic.

For the next two formulas, fix \(f\in\MM_n\). Put \(a_y:=|f^{-1}(y)|\) for \(y\in\im f\), and recall
that \(\Comp(f)\) is the set of connected components of \(G(f)\).
The composition and equality lemmas then give the following counts.

\begin{lemma}\label{lem:764783}
For every \(f\in\MM_n\), with \(g\) ranging over \(\MM_n\) in both counts,
\begin{align}
 \#\{g:fg=f^3\}
   &=\prod_{x\in X_n}a_{f^3(x)},\label{eq:764}\\
 \#\{g:fg=fgf\}
   &=\prod_{C\in\Comp(f)}
       \left(\sum_{y\in\im f}a_y^{|C|}\right).\label{eq:783}
\end{align}
\end{lemma}

\begin{proof}
The first formula is~\eqref{eq:prescribed-left} with
\(a=f\), \(b=f^3\), and \(u=g\).
For the second, put \(h:=fg\). The relation is \(h=hf\).
For \(x\in X_n\) and \(y\in\im f\), the constraint \(h(x)=y\) permits \(a_y\) choices for \(g(x)\),
independently for different \(x\). Apply
\namedref{Lemma}{lem:weighted-equalities} with
\(a=\operatorname{id}_{X_n}\), \(b=f\), and
\(\rho(x,y):=a_y\), taking zero weight when \(y\notin\im f\).
The graph in this application has components \(\Comp(f)\), so
its weighted product is~\eqref{eq:783}.
\end{proof}

\subsection{Unions of fibers with a common value}

For \(fg=fgf\), the connected components of \(G(\tau)\) correspond to
those of \(G(f)\). Indeed, all vertices of a fiber are connected through
their common image, and contracting fibers leaves the edges of \(G(\tau)\).
A component \(D\in\Comp(\tau)\) has total size \(\sum_{j\in D}s_j\),
so its factor in~\eqref{eq:783} is
\(\sum_{i=1}^k s_i^{\sum_{j\in D}s_j}\).

To sum~\eqref{eq:783}, group fibers by their common value under \(fg\).
Define
\[
 A_n:=\#\{(f,g)\in\MM_n^2:fg=fgf\}.
\]
For \(\mathbf e=(e_1,\ldots,e_n)\) with nonnegative integer entries, put
\[
 k(\mathbf e):=\sum_{a=1}^n e_a,\qquad
 w(\mathbf e):=\sum_{a=1}^n a e_a.
\]
For \(1\leq b\leq n\), define the finite polynomial
\begin{equation}\label{eq:783Xi}
 \Xi_b(\mathbf z):=
 \sum_{\substack{e_1,\ldots,e_n\geq0\\w(\mathbf e)\leq n}}
 \fall{w(\mathbf e)}{k(\mathbf e)}\,b^{w(\mathbf e)}
 \prod_{a=1}^n\frac{z_a^{e_a}}{e_a!}.
\end{equation}
The vector \(\mathbf e\) records the source fibers assigned to one
value of \(fg\), with \(b\) the size of the fiber over that value.
The empty allocation gives the constant term \(1\).

\begin{theorem}\label{thm:783total}
For \(n\geq0\),
\begin{equation}\label{eq:783total}
 A_n=n!\sum_{\substack{m_1,\ldots,m_n\geq0\\\sum_{a=1}^n a m_a=n}}
 \coef{\mathbf z^{\mathbf m}}
 \prod_{b=1}^n\left(\frac{\Xi_b(\mathbf z)}{b!}\right)^{m_b}.
\end{equation}
\end{theorem}

\begin{proof}
Fix a fiber partition \(B_1,\ldots,B_k\) of \(X_n\), with
\(s_j:=|B_j|\), and write \(f(B_j)=\{y_j\}\) for \(j\in X_k\).
Put \(h:=fg\).
If \(h=hf\), then \(h\) is constant on each \(B_j\), and every value
of \(h\) lies in \(\im f\). Thus there is a unique map
\(\lambda:X_k\to X_k\) such that \(h(B_j)=\{y_{\lambda(j)}\}\).
For a map \(h\) of this form, the relation \(h=hf\) is equivalent to
\begin{equation}\label{eq:783lambda}
 \lambda\tau=\lambda,
\end{equation}
because \(y_j\in B_{\tau(j)}\) and the image labels are distinct.

For a fixed \(\lambda:X_k\to X_k\) and each \(i\in X_k\), define
\begin{equation}\label{eq:783groups}
 \begin{aligned}
 I_i&:=\lambda^{-1}(i),&
 U_i&:=\bigcup_{j\in I_i}B_j,\\
 k_i&:=|I_i|,&
 w_i&:=|U_i|=\sum_{j\in I_i}s_j.
 \end{aligned}
\end{equation}
Condition~\eqref{eq:783lambda} says exactly that \(y_j\in U_i\)
whenever \(j\in I_i\). The \(k_i\) distinct image labels indexed by
\(I_i\) can therefore be chosen in \(\fall{w_i}{k_i}\) ways.
The sets \(U_i\) are disjoint, so these choices are independent and
make all image labels distinct. For each \(x\in U_i\), the condition
\(fg(x)=y_i\) is equivalent to \(g(x)\in B_i\), giving \(s_i^{w_i}\)
choices. Consequently the contribution of \(\lambda\) is
\[
 \prod_{i=1}^k\fall{w_i}{k_i}\,s_i^{w_i}.
\]
Empty groups are allowed and contribute \(1\). The index \(i\) need
not belong to \(I_i\): its image label \(y_i\) is chosen in the group
containing the source index \(i\). Conversely, the choices just
described give \(f(x)\in U_i\) and
\(fg(x)=y_i=fg(f(x))\) for every \(x\in U_i\), so they produce
exactly the required pairs.

For a target fiber of size \(b\), this local weight is
\(\fall{w(\mathbf e)}{k(\mathbf e)}b^{w(\mathbf e)}\), so its
allocation polynomial is \(\Xi_b\). Applying
\namedref{Lemma}{lem:weighted-block-allocation} now gives
\eqref{eq:783total}. Every pair determines its fiber partition,
its distinct image labels, and its map \(\lambda\) uniquely.
\end{proof}

For a permutation \(g\), the composition has the same image as the
first map, so every source group must be nonempty.
There are as many groups as source fibers, so \(\lambda\)
is a permutation and \(\lambda\tau=\lambda\) forces \(\tau\) to
be the identity. Thus \(f\) is idempotent. Conversely, for an
idempotent \(f\), the admissible permutations \(g\) permute its
fibers of equal size and choose a bijection between each pair of matched fibers.
If \(f\) has \(m_b\) fibers of size \(b\), the number of admissible
\(g\) is \(\prod_bm_b!(b!)^{m_b}\). Choosing the labeled fibers and the
fixed point in each gives the restricted total
\begin{equation}
 \#\{(f,g)\in\MM_n\times\operatorname{Sym}(X_n):fg=fgf\}
 =n!\coef{z^n}\prod_{b\geq1}(1-bz^b)^{-1}.
 \label{eq:permutation-783-total}
\end{equation}

For unrestricted \(g\), compare \(\Xi_b\) with \(\Psi_b\) from
\namedref{Section}{sec:profile}: a nonempty source group of \(k\) blocks
and \(w\) points has weight \(\fall{w}{k}b^w\) here and
\(\fall{b}{k}w^w\) there. The image labels are chosen from the
source union and the target fiber, respectively.

The recurrence~\eqref{eq:shared-allocation-recurrence}, with
\(W_b(\mathbf e):=\fall{w(\mathbf e)}{k(\mathbf e)}b^{w(\mathbf e)}\),
evaluates each profile's contribution using nonnegative integers.

For \(n\geq1\), the following two profiles illustrate the formula.
The one-fiber profile \(m_n=1\) has contribution
\(\coef{z_n}\Xi_n=n^{n+1}\), accounting for the \(n\) constant maps
\(f\) and all \(n^n\) maps \(g\). For the all-singleton profile
\(m_1=n\), the specialization
\(\Xi_1(z_1,0,\ldots,0)=1+z_1+\cdots+z_1^n\) gives the contribution
\[
 n!\coef{z_1^n}(1-z_1)^{-n}
   =n!\binom{2n-1}{n}=\fall{(2n-1)}{n}.
\]
Here \(f\) is a permutation, and cancellation of \(f\) gives
\(g=gf\), so \(g\) is constant on every cycle of \(f\). Summing
the resulting \(n^{c(f)}\), where \(c(f)\) is the number of cycles of
\(G(f)\), over \(f\in\mathfrak S_n\) gives the same product.

\subsection{Allocations by trees and cycles}
\label{sec:allocation-polynomial}

The allocations above can also be counted by colored trees and cycles,
with colors recording successive fiber sizes along a directed path. General determinant formulas for colored functional
digraphs are developed
in~\cite{GouldenKulkarni1997}; we derive the weights needed here.

The following lemma assembles the cycles after their attached trees are
counted. All power series are formal;
every decoration has positive size, so each coefficient involves only
finitely many cycles.

\begin{lemma}
\label{lem:cycle-assembly}
Let \(\mathcal I\) be a finite set of colors, and let
\(M=(M_{ij})_{i,j\in\mathcal I}\) be a matrix of exponential generating
functions with zero constant terms; write \(I\) for the identity
matrix of the same size and \(\operatorname{tr}\) for the matrix trace.
Interpret \(M_{ij}\) as the
series for a nonempty labeled decoration on a transition from color
\(i\) to color \(j\). The generating function for a labeled set
of directed cycles of these decorated transitions is
\begin{equation}
 \exp\left(\sum_{\ell\geq1}
       \frac{\operatorname{tr}(M^\ell)}{\ell}\right)
       =\det(I-M)^{-1}.
 \label{eq:cycle-assembly}
\end{equation}
In particular, suppose a color has at most one allowed successor:
there is a partial map \(\Phi\colon\mathcal I\dashrightarrow\mathcal I\), with
\(M_{i,\Phi(i)}=W_i\) for \(i\) in its domain and all other entries zero. Then
\begin{equation}
 \det(I-M)^{-1}
   =\prod_{C\text{ a directed cycle of }\Phi}
        \left(1-\prod_{i\in C}W_i\right)^{-1}.
 \label{eq:functional-cycle-assembly}
\end{equation}
\end{lemma}

\begin{proof}
The trace \(\operatorname{tr}(M^\ell)\) counts decorated cycles
with a distinguished starting transition. Every cycle has exactly
\(\ell\) such choices: its decorations occupy disjoint nonempty
label sets, even if the sequence of colors is periodic. Dividing by
\(\ell\) and then taking a labeled set of cycles gives the
exponential in~\eqref{eq:cycle-assembly}. The determinant identity
follows from
\[
 -\log\det(I-M)=
       \sum_{\ell\geq1}\frac{\operatorname{tr}(M^\ell)}{\ell},
\]
which is valid formally because \(M(0)=0\).

For the partial-map case, every closed walk follows a directed cycle.
A cycle \(C\) of length \(p\), traversed \(r\) times, contributes
\(p(\prod_{i\in C}W_i)^r\) to the trace of \(M^{pr}\).
Its contribution to the logarithm is therefore
\[
 \sum_{r\geq1}\frac1r\left(\prod_{i\in C}W_i\right)^r
       =-\log\left(1-\prod_{i\in C}W_i\right).
\]
Exponentiation proves~\eqref{eq:functional-cycle-assembly}.
\end{proof}

For example, introduce vertex-weight variables \(u_i\) and edge-weight
variables \(k_{ij}\), for \(i,j\in\mathcal I\). Let \(R_i\) be the
series for trees directed toward a root of color \(i\). Then
\begin{equation}
 R_i=u_i\exp\left(\sum_{j\in\mathcal I}k_{ji}R_j\right).
 \label{eq:common-colored-trees}
\end{equation}
Indeed, the root has an unordered set of incoming child trees; their
edge weights are \(k_{ji}\). A cyclic edge has weight \(k_{ij}\),
so \namedref{Lemma}{lem:cycle-assembly} with
\(M_{ij}:=R_i k_{ij}\) gives the generating function for these
labeled functional digraphs.
For a prescribed partial color map \(\Phi\), take
\(k_{ij}:=\mathbf1_{\{j=\Phi(i)\}}\) when \(\Phi(i)\) is defined, and
\(k_{ij}:=0\) otherwise. The tree series and the full digraph series are
\[
 R_i=u_i\exp\left(\sum_{j:\Phi(j)=i}R_j\right),
 \qquad
 \prod_C\left(1-\prod_{i\in C}R_i\right)^{-1}.
\]
The product ranges over the directed cycles \(C\) of \(\Phi\).
Other local restrictions change the tree and cyclic-vertex weights;
the cycle formulas remain the same.

The determinant expression also admits an integer coefficient
recurrence that will be used for several of the allocation formulas.

\begin{lemma}\label{lem:matrix-coefficients}
Fix integers \(D,N\geq1\). Let \(M(z)\) be a \(D\)-by-\(D\)
matrix of formal power series with \(M(0)=0\), and
suppose the matrices \(M_m:=m!\coef{z^m}M(z)\) are integral for
\(1\leq m\leq N\). Put
\[
 R(z):=(I-M(z))^{-1},\qquad Z(z):=\det(I-M(z))^{-1},\qquad
 L(z):=\operatorname{tr}(R(z)M'(z)),
\]
and, for integers \(m\geq0\), write \(R_m:=m!\coef{z^m}R(z)\),
\(\ell_m:=m!\coef{z^m}L(z)\), and
\(z_m:=m!\coef{z^m}Z(z)\). Then
\begin{equation}\label{eq:matrix-coefficients}
\begin{aligned}
 R_0&=I,&
 R_m&=\sum_{j=1}^m\binom mj M_jR_{m-j}
                 &&(1\leq m\leq N),\\
 &&\ell_m&=\sum_{j=0}^m\binom mj
             \operatorname{tr}(R_jM_{m-j+1})
                 &&(0\leq m<N),\\
 z_0&=1,&
 z_m&=\sum_{j=0}^{m-1}\binom{m-1}j\ell_jz_{m-1-j}
                 &&(1\leq m\leq N).
\end{aligned}
\end{equation}
\end{lemma}

\begin{proof}
The first recurrence follows from \(R=I+MR\). Differentiate the
logarithm in~\eqref{eq:cycle-assembly}; cyclicity of the trace gives
\[
 \frac{d}{dz}\operatorname{tr}(M^\ell)
 =\ell\operatorname{tr}(M^{\ell-1}M'),\qquad
 Z'=\operatorname{tr}(RM')Z=LZ.
\]
Binomial convolution of exponential coefficients now gives the other
two recurrences in~\eqref{eq:matrix-coefficients}.
\end{proof}

We now apply the tree and cycle decomposition when a local allocation weight depends
only on the target size, the number of incoming blocks, and their
total size.

Fix an integer \(N\geq1\). For nonnegative integer vectors
\(\mathbf e=(e_1,\ldots,e_N)\) with \(\sum_jje_j\leq N\), suppose the
weights in \namedref{Lemma}{lem:weighted-block-allocation} have the form
\[
 W_b(\mathbf e)=\mathcal W(b,k,w),\qquad
 k:=\sum_j e_j,\quad w:=\sum_j j e_j.
\]
Their values are specified for \(1\leq b\leq N\) and
\(0\leq k,w\leq N\), and extended by zero outside these ranges.
Define \(T_b(z)\) for \(1\leq b\leq N\), together with \(F(z;u,v)\), by
\begin{equation}\label{eq:allocation-tree-series}
\begin{aligned}
 F(z;u,v)&:=\exp\!\left(u\sum_{j=1}^N T_j(z)v^j\right),\\
 T_b(z)&=\frac{z^b}{b!}
       \sum_{k,w=0}^N\mathcal W(b,k,w)\coef{u^k v^w}F(z;u,v).
\end{aligned}
\end{equation}
The factor \(z^b\) determines their coefficients in increasing degree.
All series are truncated after degree \(N\) in \(z\).
Define the \(N\)-by-\(N\) matrix \(B(z)\) by
\begin{equation}\label{eq:allocation-cycle-matrix}
 B_{ab}(z):=\frac{z^b}{b!}
       \sum_{k,w=0}^N\mathcal W(b,k+1,w+a)
                      \coef{u^k v^w}F(z;u,v)
       \qquad(1\leq a,b\leq N).
\end{equation}

\begin{theorem}\label{thm:allocation-tree-cycle}
For \(0\leq m\leq N\), let \(Z_m\) be the weighted total in
\namedref{Lemma}{lem:weighted-block-allocation} on a set of size \(m\),
with the local weights just specified. Then
\begin{equation}\label{eq:allocation-tree-cycle-total}
 \sum_{m=0}^N Z_m\frac{z^m}{m!}
 \equiv \det(I-B(z))^{-1}\pmod{z^{N+1}}.
\end{equation}
\end{theorem}

\begin{proof}
Decompose the allocation map on the blocks into directed cycles with
rooted trees attached. A block of size \(b\) is a labeled set and
contributes \(z^b/b!\). The series \(T_b(z)\) counts weighted rooted
trees whose root block has size \(b\). The variables \(u,v\) in
\(F(z;u,v)\) record the number and total size of the root blocks in
an incoming forest. The second equation in
\eqref{eq:allocation-tree-series} applies the local weight after
attaching that forest to the root.

A cycle vertex of size \(b\) with predecessor size \(a\) has one
additional incoming block, of size \(a\). Its series is therefore
\(B_{ab}(z)\). Extending \(\mathcal W\) by zero outside the supplied
ranges does not affect the totals through size \(N\): the incoming
blocks at a cycle vertex, including its predecessor, have at most
\(N\) points and at most \(N\) blocks in any contributing structure.
\namedref{Lemma}{lem:cycle-assembly} now gives
\eqref{eq:allocation-tree-cycle-total}.
\end{proof}

For integer evaluation, define the following coefficients for \(1\leq b\leq N\)
and integers \(0\leq m,k,w\leq N\):
\[
 t_{b,m}:=m!\coef{z^m}T_b(z),\qquad
 f_{m,k,w}:=m!\coef{z^m u^k v^w}F(z;u,v),\qquad
 B_m:=m!\coef{z^m}B(z).
\]
Set \(f_{0,0,0}=1\) and all other \(f_{0,k,w}=0\);
out-of-range coefficients are zero.
For \(m\geq1\), equations~\eqref{eq:allocation-tree-series}
and~\eqref{eq:allocation-cycle-matrix} give
\begin{equation}\label{eq:allocation-coefficient-recursions}
\begin{aligned}
 t_{b,m}&=\binom mb
       \sum_{k,w}\mathcal W(b,k,w)f_{m-b,k,w},\\
 f_{m,k,w}&=\sum_{r=1}^m\binom{m-1}{r-1}
       \sum_{b=1}^r t_{b,r}f_{m-r,k-1,w-b},\\
 (B_m)_{ab}&=\binom mb
       \sum_{k,w}\mathcal W(b,k+1,w+a)f_{m-b,k,w}.
\end{aligned}
\end{equation}
At degree \(m\), first compute the \(t_{b,m}\), for \(b\leq m\),
and then the \(f_{m,k,w}\). Only coefficients with
\(0\leq k\leq w\leq m\leq N\) need be retained.

\namedref{Theorem}{thm:allocation-tree-cycle} applies to the counts for
\(f^2=f^2g\) and \(fg=fgf\) with
\[
 \mathcal W(b,k,w)=\fall{b}{k}w^w,
 \qquad
 \mathcal W(b,k,w)=\fall{w}{k}b^w,
\]
respectively. The allocation map is \(\tau\) in the first case
and the map \(\lambda\) in~\eqref{eq:783lambda} in the second.
The condition involving \(\tau\) in the second case has already
been included in the factor \(\fall{w}{k}\), as proved in
\namedref{Theorem}{thm:783total}. Thus both are arbitrary allocation maps
with the stated local weights.

For bijective \(g\), use instead
\(\mathcal W(b,k,w)=\fall{b}{k}w!\) and
\(\mathcal W(b,k,w)=\fall{w}{k}b!\,\mathbf1_{\{w=b\}}\),
respectively. For the second identity,
\eqref{eq:permutation-783-total} gives a specialized product formula.
Requiring \(f\) to be bijective restricts the block sizes to one.

For \(f^2=fgf\), the local weight depends on the incoming block
count \(k\), but not on their total size \(w\). All source-size
variables then enter through \(\sum_jz_j\), giving the simplification
in \namedref{Subsection}{sec:companion-fgf-total}. The weights for the two identities
treated here depend on both \(k,w\). For the next identity, they also
depend on the successor size; we use pairs of sizes as colors.

\subsection{Sizes of fibers and their successors}

For \(fg=f^3\), define
\[
 D_n:=\#\{(f,g)\in\MM_n^2:fg=f^3\}.
\]
By~\eqref{eq:iterates}, the first count in \namedref{Lemma}{lem:764783} becomes
\begin{equation}\label{eq:764tau}
 \prod_{j=1}^k s_{\tau^2(j)}^{s_j}.
\end{equation}
Record the size of each fiber and of its successor: for \(1\leq a,b\leq n\),
let \(m_{ab}\) be the number of size-\(a\) fibers with successor size \(b\).
Then \(r_a:=\sum_bm_{ab}\) is the number of size-\(a\) fibers,
and \(\sum_a ar_a=n\).

\begin{theorem}\label{thm:764total}
For \(n\geq0\), let \(M=(m_{ab})_{1\leq a,b\leq n}\) range over
nonnegative integer matrices satisfying \(\sum_{a,b}a m_{ab}=n\).
For each such matrix and \(1\leq b\leq n\), define
\begin{equation}\label{eq:764K}
 K_b(M):=
 \coef{z_1^{m_{1b}}\cdots z_n^{m_{nb}}}
 \prod_{c=1}^n
       \left(1+\sum_{a=1}^n c^a z_a\right)^{b m_{bc}}.
\end{equation}
Then
\begin{equation}\label{eq:764total}
 D_n=n!\sum_M
       \frac{\prod_{b=1}^n K_b(M)}
            {\prod_{a=1}^n(a!)^{r_a}}.
\end{equation}
\end{theorem}

\begin{proof}
Fix a fiber partition with profile \((r_a)\). Assigning to its
size-\(a\) blocks the prescribed successor sizes has
\(r_a!/\prod_bm_{ab}!\) choices. Once these assignments are made, the
\(m_{ab}\) source blocks of size \(a\) must place their image labels in
target blocks of size \(b\).

Consider a particular target block of size \(b\) whose own successor has
size \(c\). If it receives \(e_a\) source blocks of size \(a\), their
distinct image labels can be chosen in \(\fall{b}{\sum_ae_a}\) ways.
Their contribution to~\eqref{eq:764tau} is
\(\prod_a c^{a e_a}\). The polynomial for this target is consequently
\[
 \sum_{\sum_ae_a\leq b}
 \fall{b}{\sum_ae_a}
 \prod_a\frac{(c^a z_a)^{e_a}}{e_a!}
 =\left(1+\sum_a c^a z_a\right)^b.
\]
There are \(m_{bc}\) target blocks of this kind. The multinomial argument
in \namedref{Lemma}{lem:weighted-block-allocation}, now applied separately
for each successor size, gives \(K_b(M)\) by multiplying
their polynomials and extracting the specified coefficient.
Multiplication by \(\prod_a m_{ab}!\) distinguishes the source blocks of
each size, so this product times \(K_b(M)\) counts exactly their image-label
assignments and compatible choices of \(g\).

Taking the product over \(b\), the factorials \(m_{ab}!\) cancel those
from the initial assignment of successor sizes. The contribution of the
fixed partition is therefore
\(\prod_a r_a!\prod_b K_b(M)\). Finally, there are
\[
 \frac{n!}{\prod_a r_a!(a!)^{r_a}}
\]
partitions with this profile. The remaining multiplicity factorials cancel,
giving~\eqref{eq:764total}. Every map has a unique fiber partition, and its
successor-size assignments and image labels determine \(M\) and \(\tau\),
so every pair is counted once.
\end{proof}

For fixed \(f\in\MM_n\), restricting \(g\) to permutations gives
\[
 \#\{g\in\operatorname{Sym}(X_n):fg=f^3\}
 =\begin{cases}
   \displaystyle\prod_{y\in f(X_n)}|f^{-1}(y)|!,
       &|(f^3)^{-1}(y)|=|f^{-1}(y)|\quad(y\in X_n),\\
   0,&\text{otherwise}.
  \end{cases}
\]
This follows from~\eqref{eq:bijective-prescribed-left}. The equality
of fiber sizes forces \(f^2(X_n)=f(X_n)\).
In the notation of \namedref{Lemma}{lem:labels}, \(\tau\) is consequently
a permutation, and the remaining condition is
\(s_{\tau^2(j)}=s_j\) for every \(j\in X_k\). Thus the corresponding
allocation uses only cycles satisfying this size condition, with
weight \(\prod_j s_j s_j!\) including the choices of image labels.

Summing over these cycles gives
\begin{equation}
 \#\{(f,g)\in\MM_n\times\mathfrak S_n:fg=f^3\}
 =n!\coef{z^n}\prod_{b\geq1}(1-bz^b)^{-1}
                \prod_{1\leq a<b}(1-abz^{a+b})^{-1}.
 \label{eq:permutation-764-total}
\end{equation}
Indeed, the sizes around a cycle of \(\tau\) have period at most two.
Combining the weight \(\prod_j s_j s_j!\) with the factorials for
labeled fiber blocks gives a size-\(b\) vertex weight \(bz^b\).
Cycles of a single size give the first product; cycles alternating
between two distinct sizes \(a,b\) give the second, by
\namedref{Lemma}{lem:cycle-assembly}.

For the unrestricted sum, only sizes in \(S:=\{a\in X_n:r_a>0\}\) need be retained. Indeed, if \(r_b=0\),
the polynomial in~\eqref{eq:764K} is \(1\), so \(K_b(M)=0\) whenever
column \(b\) is nonzero. A second necessary condition for a nonzero term is
\begin{equation}\label{eq:764capacity}
 \sum_a m_{ab}\leq b r_b\qquad(b\in S),
\end{equation}
since the \(r_b\) target blocks together have only \(b r_b\) positions
for distinct image labels.

To compute \(K_b(M)\), fix \(b\in S\) and list
\(S=\{a_1,\ldots,a_t\}\). Regard the factor with index \(c\)
in~\eqref{eq:764K} as \(b m_{bc}\) distinct positions, each either unused or
assigned one of the variables recording a source size. Let
\(\boldsymbol\ell=(\ell_c)_{c\in S}\), where \(0\leq\ell_c\leq bm_{bc}\),
record the numbers of remaining positions. Define
\(R_{t+1}(\boldsymbol\ell):=1\) and, for \(i=t,t-1,\ldots,1\), set
\begin{equation}\label{eq:764recurrence}
 R_i(\boldsymbol\ell)
 :=\sum_{\substack{0\leq e_c\leq\ell_c\ (c\in S)\\
                         \sum_{c\in S}e_c=m_{a_i b}}}
       \left(\prod_{c\in S}\binom{\ell_c}{e_c}c^{a_i e_c}\right)
       R_{i+1}(\boldsymbol\ell-\mathbf e).
\end{equation}
Then \(K_b(M)=R_1((b m_{bc})_{c\in S})\). The binomial factors choose the
positions receiving sources of the current size; no extra factorial
is required.

For \(n\geq1\), the all-singleton profile has only \(m_{11}=n\) and
\(K_1=1\), contributing \(n!\): every permutation \(f\) forces
\(g=f^2\). The one-fiber profile has \(m_{nn}=1\) and
\(K_n=\coef{z_n}(1+n^nz_n)^n=n^{n+1}\), accounting for the \(n\) constant
maps \(f\) and the \(n^n\) arbitrary maps \(g\).

The fiber partitions in Formula~\eqref{eq:764total} can also be summed
by trees and cycles to obtain a determinant expression. Fix an integer
\(N\geq1\). For \(1\leq a,b\leq N\), define the formal power series
\(T_{a,b}(x)\), with zero constant term, by
\begin{equation}\label{eq:764trees}
 T_{a,b}(x)
 :=\sum_{c=1}^N\frac{(bx)^c}{c!}
                     \bigl(1+T_{c,a}(x)\bigr)^c.
\end{equation}
These equations determine a unique solution recursively: the coefficient
of degree \(m\) on the right depends only on coefficients of degrees
below \(m\). For \(1\leq a,b,c\leq N\), define the corresponding
cycle weight by
\[
 W(a,b,c):=\frac{x^a c^a}{(a-1)!}
                \bigl(1+T_{a,b}(x)\bigr)^{a-1},
\]
and index an \(N^2\)-by-\(N^2\) matrix \(\mathsf W_N(x)\) by ordered
pairs in \(X_N^2\), with entries
\begin{equation}\label{eq:764matrix}
 (\mathsf W_N)_{(a,b),(d,c)}
 :=\begin{cases}
       W(a,b,c),&d=b,\\
       0,&d\ne b.
    \end{cases}
\end{equation}

\begin{theorem}\label{thm:764trees}
For \(0\leq n\leq N\),
\begin{equation}
 D_n=n!\,\coef{x^n}\det(I-\mathsf W_N(x))^{-1}.
 \label{eq:764det}
\end{equation}
\end{theorem}

\begin{proof}
Use the functional digraph of the block map \(\tau\). A block of size
\(c\) is a nonempty labeled set and contributes \(x^c/c!\) to its
exponential generating function. Interpret \(T_{a,b}\) as the series for
one incoming noncyclic branch attached to a block of size \(a\) whose
successor has size \(b\). The image label of the branch root inside its
parent block is not chosen yet; it will be chosen at that parent together
with its other incoming image labels.

If the branch root has size \(c\), its second successor has size \(b\),
so~\eqref{eq:764tau} supplies \(b^c\). If it has \(d\) incoming child
branches, their unordered set contributes \(T_{c,a}^d/d!\), and their
distinct image labels can be assigned in \(\fall{c}{d}\) ways. Summing
over \(d\) gives \((1+T_{c,a})^c\). Together with \(x^c/c!\), these
factors prove~\eqref{eq:764trees}.

Now consider a cyclic block of size \(a\), with successor size \(b\) and
second-successor size \(c\). Its one cyclic predecessor already needs an
image label inside the block. If \(d\) noncyclic branches are attached,
all their labels and the cyclic predecessor's label must be distinct.
Their combined factor is therefore
\[
 \sum_{d=0}^{a-1}\frac{\fall{a}{d+1}}{d!}T_{a,b}^d
 =a(1+T_{a,b})^{a-1}.
\]
Multiplication by \(x^a/a!\) and \(c^a\) gives \(W(a,b,c)\).
For a loop the cyclic predecessor is the block itself, and still consumes
one image-label position. For a two-cycle the second successor is the
original block, exactly as required by~\eqref{eq:764tau}.

A closed walk in the ordered-pair indices of \(\mathsf W_N\)
specifies a cyclic sequence of block sizes, with the corresponding
decorated cyclic-vertex weights. Every block occupies a nonempty
label set. \namedref{Lemma}{lem:cycle-assembly} therefore gives
\eqref{eq:764det}, including periodic size sequences.

Every block contributes a positive power of \(x\). Consequently a
coefficient of degree \(n\leq N\) uses only sizes at most \(n\)
and cycles of length at most \(n\). Thus restricting sizes to at most
\(N\) does not change these coefficients.
\end{proof}

For coefficients through degree \(N\), this system has \(N^2\)
auxiliary series.
\section{Functional-digraph homomorphisms and iterate types}
\label{sec:hom}

\subsection{Counting homomorphisms by a tree recurrence}

We first count homomorphisms between two prescribed functional digraphs.
For fixed maps \(s:U\to U\) and \(t:V\to V\) on finite sets, define
\[
 H(s,t):=\#\{q:U\to V:qs=tq\}.
\]
The equation \(qs=tq\) means that
\(q\) preserves directed edges from \(G(s)\) to \(G(t)\): for each
\(u\in U\), the edge \(u\to s(u)\) maps to \(q(u)\to t(q(u))\).
Novotn\'y~\cite{Novotny1963} studies such maps; Nelson~\cite{Nelson1982} corrects the
general existence criterion. Higgins~\cite[Theorem~5]{Higgins1988}
uses functional digraphs to solve finite transformation equations of
the form \(qs=tq\). For centralizers, Ara\'ujo and Konieczny~\cite[Propositions~3.1
and~4.1, Corollary~5.1]{AraujoKonieczny2013}
give the edge-preserving interpretation and decomposition into maps between
components. We use these constructions for
\((h^2,h)\) and \((h,h^2)\), where \(h\in\MM_n\).

Holloway and Shattuck~\cite[Propositions~2.6 and~3.7]{HollowayShattuck2013}
give a cycle-type formula when one commuting map is a permutation and
count the general commuting case by successive levels, using cardinalities
of compatible map sets. The recurrence below evaluates \(H(s,t)\)
before summing over isomorphism types.

Fixing the map that is squared gives the two totals
\begin{equation}\label{eq:homtotals}
 \NN_{fg,gff}(n)=\sum_{h\in\MM_n}H(h^2,h),\qquad
 \NN_{fg,ggf}(n)=\sum_{h\in\MM_n}H(h,h^2).
\end{equation}
At \(n=4\), these counts are \(2776\) and \(3568\), respectively.

A vertex \(u\in U\) is \emph{cyclic} if \(s^j(u)=u\) for some positive
integer \(j\). Delete all edges between cyclic vertices of \(G(s)\).
The remaining graph is a forest of trees directed toward the cyclic
vertices. For each \(u\in U\), let \(T_s(u)\) be its incoming subtree
in this forest, and put
\[
 \operatorname{Ch}_s(u):=\{w\in U:w\text{ is noncyclic and }s(w)=u\}.
\]
For \(u\in U\) and \(v\in V\), let \(Q(u,v)\) be the number of
homomorphisms \(T_s(u)\to G(t)\) sending \(u\) to \(v\).
Let \(\mathcal C(s)\) contain each directed cycle of \(G(s)\) once, with an
arbitrary starting point and its forward cyclic order.

\begin{theorem}\label{thm:homDP}
For \(u\in U\) and \(v\in V\), the subtree counts satisfy
\begin{equation}\label{eq:homDP}
 Q(u,v)=\prod_{w\in\operatorname{Ch}_s(u)}
             \left(\sum_{z\in V:t(z)=v}Q(w,z)\right).
\end{equation}
In particular, \(Q(u,v)=1\) when \(u\) is a leaf. The total count is
\begin{equation}\label{eq:homcycle}
 H(s,t)=\prod_{(c_0,\ldots,c_{\ell-1})\in\mathcal C(s)}
       \left(\sum_{v\in V:t^\ell(v)=v}
                   \prod_{i=0}^{\ell-1}Q(c_i,t^i(v))\right).
\end{equation}
\end{theorem}

\begin{proof}
Once \(q(u)=v\), each child \(w\) must map to a preimage of \(v\) under
\(t\). Sum over these possible images, then multiply the independent
counts for the child subtrees. This proves~\eqref{eq:homDP}.

On a cycle of the domain graph, choosing \(q(c_0)=v\) forces
\(q(c_i)=t^i(v)\). These values are consistent precisely when
\(t^\ell(v)=v\). The trees attached to that cycle contribute the displayed
product of subtree counts. Distinct domain components can be mapped
independently; their target images need not be disjoint.
\end{proof}

Edges in a tree of the domain graph may map to cycle edges of \(G(t)\),
so the image of that tree need not lie in a single tree of the target graph.

If either \(s\) or \(t\) is a permutation, the homomorphism
recurrence simplifies to
\[
 H(s,t)=\prod_{C\in\mathcal C(s)}|\Fix(t^{|C|})|.
\]
When \(s\) is a permutation its components are cycles;
when \(t\) is a permutation every required predecessor is unique. Write \(\Aut(s)\) for the automorphism group of \(G(s)\).
A bijective homomorphism is an isomorphism. On sets of the same
cardinality, the number of such \(q\) is therefore \(|\Aut(s)|\) if
\(G(s)\) and \(G(t)\) are isomorphic, and zero otherwise.
In~\eqref{eq:homtotals}, restrict \(h\) to permutations by retaining
only permutation types, or require \(q\) to be bijective by using
this isomorphism count in place of \(H\).
The isomorphism test includes the attached trees.

\subsection{Component types and their squares}

Write \(B_n:=\NN_{fg,gff}(n)\) and \(C_n:=\NN_{fg,ggf}(n)\) for the
two totals in~\eqref{eq:homtotals}. To sum by isomorphism type, we need the
connected-component types of \(G(h)\), for \(h\in\MM_n\), together with
their automorphism orders and the effect of squaring.

For an integer \(N\geq0\), let \(\mathcal K_N\) be the finite set of isomorphism types of nonempty
connected functional digraphs with at most \(N\) vertices; thus
\(\mathcal K_0=\varnothing\). For
\(C\in\mathcal K_N\), write
\[
 w_C:=|C|,\qquad
 \ell_C:=\text{the length of its directed cycle},\qquad
 \alpha_C:=|\Aut(C)|.
\]
Here \(|C|\) is the number of vertices of a representative of \(C\),
whose automorphisms are the permutations preserving its directed edges.
We use any fixed representative map when taking a square or counting
homomorphisms; the choice does not affect these counts.

A rooted-tree type is determined by the multiset of its child-subtree
types. For a rooted-tree type \(T\), write
\(\alpha(T):=|\Aut(T)|\), with the root required to be fixed. If a child
type \(U\) occurs \(m_U\) times, then
\[
 \alpha(T)=\prod_U m_U!\,\alpha(U)^{m_U}.
\]
A connected type is a directed cycle
\(C=(T_0,\ldots,T_{\ell-1})\) of rooted-tree types, considered up to
cyclic rotation. If \(r(C)\) rotations preserve this tuple, then
\begin{equation}\label{eq:homAut}
 \alpha_C=r(C)\prod_{i=0}^{\ell-1}\alpha(T_i).
\end{equation}
Indeed, an automorphism acts on the directed cycle by a rotation that
preserves its attached tree types, and its restrictions to those trees
can then be chosen independently. Arbitrary functional-digraph types
are multisets of connected types. Defrain, Porreca, and
Timofeeva~\cite{DefrainPorrecaTimofeeva2024} give an algorithm that
generates one representative of each isomorphism type.

Let \(\mathfrak T_n\) be the set of functional-digraph types on
\(n\) vertices, and write \(\operatorname{type}(h)\) for the type
of a map \(h\). For \(t\in\mathfrak T_n\), choose a representative
\(h_t\in\MM_n\), and put
\[
 \omega_n(t):=\frac{n!}{|\Aut(t)|}.
\]
This is the number of labeled realizations of \(t\), as in Fuller's
commuting-map tables for \oeis{A181162}. The following
lemma applies it to fixed-map counts and counts grouped by a common iterate.

\begin{lemma}
\label{lem:type-summation}
If a numerical weight \(A(h)\), defined for \(h\in\MM_n\), is
unchanged by relabeling its vertices, then
\begin{equation}
 \sum_{h\in\MM_n}A(h)
       =\sum_{t\in\mathfrak T_n}\omega_n(t)A(h_t).
 \label{eq:type-summation}
\end{equation}
More generally, let \(P,Q\colon\MM_n\to\MM_n\) commute with
relabeling, and, for \(t\in\mathfrak T_n\), define
\[
 H_P(t):=\sum_{\substack{s\in\mathfrak T_n\\
                         \operatorname{type}(P(h_s))=t}}\omega_n(s),
 \qquad
 H_Q(t):=\sum_{\substack{s\in\mathfrak T_n\\
                         \operatorname{type}(Q(h_s))=t}}\omega_n(s).
\]
Then
\begin{equation}
 \#\{(f,g)\in\MM_n^2:P(f)=Q(g)\}
       =\sum_{t\in\mathfrak T_n}\frac{H_P(t)H_Q(t)}{\omega_n(t)}.
 \label{eq:equivariant-histograms}
\end{equation}
\end{lemma}

\begin{proof}
Conjugation by the symmetric group of \(X_n\) gives exactly
\(n!/|\Aut(t)|\) distinct labeled realizations of type \(t\).
Summing their common value of \(A\) proves~\eqref{eq:type-summation}.
Because \(P\) commutes with conjugation, every labeled map of type
\(t\) has the same number of preimages under \(P\), namely
\(H_P(t)/\omega_n(t)\). The corresponding number for \(Q\) is
\(H_Q(t)/\omega_n(t)\). Multiply these numbers for each common
target, and sum over its \(\omega_n(t)\) labeled realizations.
This proves~\eqref{eq:equivariant-histograms}.
\end{proof}

To describe the effect of squaring, for \(C\in\mathcal K_N\), let
\(\sigma(C)\) be the multiset of connected types in the functional
digraph of its square \(C^2\).
It has one member if \(\ell_C\) is odd and two if \(\ell_C\) is even.
Higgins~\cite[Lemmas~3--4]{Higgins1988} describes these component powers.
For a rooted tree \(T\), define \(E(T)\) to be the rooted tree consisting
of the even-depth vertices of \(T\), with each nonroot vertex joined to its
grandparent. Its child-subtree multiset is obtained by applying \(E\)
to every grandchild subtree of \(T\).

Write \(\operatorname{Ch}(T)\) for the multiset of rooted subtrees at
the children of the root of \(T\). For \(C=(T_0,\ldots,T_{\ell-1})\),
label the cycle vertices \(0,\ldots,\ell-1\), with edges from \(i\)
to \(i+1\) and indices taken modulo \(\ell\). The tree attached to
vertex \(i\) in \(C^2\) has child multiset
\begin{equation}\label{eq:homSquareTrees}
 \operatorname{Ch}(E(T_i))
 \ \sqcup\
 \{E(U):U\in\operatorname{Ch}(T_{i-1})\},
\end{equation}
The first part comes from positive even depths at \(i\); the second
comes from odd depths at the preceding cycle vertex. The cycle edges
now go from \(i\) to \(i+2\). Taking the resulting one or two cycles,
with their attached trees, computes \(\sigma(C)\) without choosing
vertex labels.

For \(C,D\in\mathcal K_N\), \namedref{Theorem}{thm:homDP} computes
\(H(C,D)\) from representatives, and
\begin{equation}\label{eq:homSupport}
 H(C,D)>0\quad\Longleftrightarrow\quad \ell_D\mid\ell_C.
\end{equation}
Necessity follows by restricting a homomorphism to the domain cycle.
For sufficiency, map that cycle into the target cycle and each attached
tree backward along the same cycle, using the unique cyclic predecessor
at each step.

Component sizes, cycle lengths, and fiber sizes do not determine these
homomorphism counts. On \(\{0,1,2,3,4\}\),
the maps
\[
 h_1=(0,0,0,1,3),\qquad h_2=(0,0,0,1,2)
\]
are both connected, have a fixed-point cycle, and have fiber sizes
\((3,1,1)\). Nevertheless,
\[
 \bigl(H(h_1^2,h_1),H(h_1,h_1^2)\bigr)=(36,20),\qquad
 \bigl(H(h_2^2,h_2),H(h_2,h_2^2)\bigr)=(81,25).
\]
The tree recurrence and direct enumeration give these values.

\subsection{The two total enumerations}

For a nonnegative integer vector
\(\mathbf m=(m_C)_{C\in\mathcal K_n}\) of component multiplicities,
define the multiplicities of the square, for \(D\in\mathcal K_n\), by
\begin{equation}\label{eq:homSquaredProfile}
 m_D^{[2]}:=\sum_{C\in\mathcal K_n}m_C\,\nu_D(C),\qquad
 \nu_D(C):=\text{the multiplicity of \(D\) in \(\sigma(C)\)}.
\end{equation}
Thus \(\mathbf m\mapsto\mathbf m^{[2]}\) is linear. Squaring preserves the vertex set, so
\(\sum_Dw_Dm_D^{[2]}=\sum_Cw_Cm_C\). For another nonnegative integer
vector \(\mathbf r=(r_S)_{S\in\mathcal K_n}\), define the homomorphism
count and the target automorphism factor by
\begin{equation}\label{eq:homProfile}
 J(\mathbf r,\mathbf m):=
 \prod_{S\in\mathcal K_n}
       \left(\sum_{D\in\mathcal K_n}m_DH(S,D)\right)^{r_S},
 \qquad
 Z(\mathbf m):=\prod_{C\in\mathcal K_n}m_C!\,\alpha_C^{m_C}.
\end{equation}

\begin{theorem}\label{thm:homTotal}
For \(n\geq0\),
\begin{align}
 B_n&=n!\sum_{\substack{\mathbf m\geq0\\
                       \sum_C w_Cm_C=n}}
          \frac{J(\mathbf m^{[2]},\mathbf m)}{Z(\mathbf m)},
          \label{eq:homTotalFrom}\\
 C_n&=n!\sum_{\substack{\mathbf m\geq0\\
                       \sum_C w_Cm_C=n}}
          \frac{J(\mathbf m,\mathbf m^{[2]})}{Z(\mathbf m)}.
          \label{eq:homTotalTo}
\end{align}
\end{theorem}

\begin{proof}
Fix a map \(h\) with component multiplicities \(\mathbf m\).
A connected domain component of type \(S\) must map into one target
component. There are \(m_D\) targets of type \(D\), each admitting
\(H(S,D)\) homomorphisms. Distinct domain components are mapped
independently, even when they have the same type or choose the same target. Thus \(J(\mathbf r,\mathbf m)\) is the
homomorphism count for source and target multiplicities
\(\mathbf r\) and \(\mathbf m\).

The connected components of \(G(h^2)\) have multiplicities \(\mathbf m^{[2]}\).
Consequently the two numerators in~\eqref{eq:homTotalFrom}
and~\eqref{eq:homTotalTo} are exactly \(H(h^2,h)\) and \(H(h,h^2)\).
In particular, the two components obtained by squaring an even-cycle
component may choose their targets independently.

An automorphism of \(h\) permutes the \(m_C\) components of each
type and chooses an automorphism within each component. Hence its
automorphism group has order \(Z(\mathbf m)\), and exactly
\(n!/Z(\mathbf m)\) labeled maps have this type. Multiplying by
this class size and applying \namedref{Lemma}{lem:type-summation} gives
both formulas.
\end{proof}

The component descriptions, squaring rule, and
\namedref{Theorem}{thm:homDP} determine every term. For integer evaluation, form the class size
\(n!/Z(\mathbf m)\) first and multiply by the appropriate value of
\(J\).

More generally, for a positive integer \(j\), let \(\mathbf m^{[j]}\)
record the component multiplicities obtained by taking the \(j\)th power
of each component representative. The same
argument gives, for any fixed positive integers \(r,s\),
\[
 \sum_{h\in\MM_n}H(h^r,h^s)
 =n!\sum_{\substack{\mathbf m\geq0\\\sum_Cw_Cm_C=n}}
       \frac{J(\mathbf m^{[r]},\mathbf m^{[s]})}{Z(\mathbf m)}.
\]

\subsection{A generating function for \texorpdfstring{\(C_n\)}{Cn}}

Fix an integer \(N\geq1\). To express \(C_n\) by rooted trees and cycles,
let \(K\) be the matrix indexed by \(C,D\in\mathcal K_N\), with entries
\[
 K(C,D):=\sum_{S\in\sigma(D)}H(C,S)=H(C,D^2),
\]
where repeated types in the multiset \(\sigma(D)\) are counted with
their multiplicities. Introduce independent variables
\(\mathbf u=(u_C)_{C\in\mathcal K_N}\) to record component types, and define formal power series
with zero constant term, one for each \(C\in\mathcal K_N\), by
\begin{equation}\label{eq:forwardR}
 R_C(\mathbf u):=
 u_C\exp\!\left(\sum_{D\in\mathcal K_N}K(D,C)R_D(\mathbf u)\right).
\end{equation}
The system has a unique solution: the coefficient of total degree \(d\)
on the right uses coefficients of degrees less than \(d\). Define a matrix
\(M_N(x)\), indexed by \(C,D\in\mathcal K_N\), by
\[
 M_N(x)_{C,D}:=
 R_C\bigl((x^{w_E}/\alpha_E)_{E\in\mathcal K_N}\bigr)K(C,D).
\]
Every entry of \(M_N(x)\) has zero constant term.

\begin{theorem}\label{thm:forwarddet}
For \(0\leq n\leq N\),
\begin{equation}\label{eq:forwarddet}
 C_n=n!\coef{x^n}\det(I-M_N(x))^{-1}.
\end{equation}
\end{theorem}

\begin{proof}
Fix \(h\in\MM_n\), and consider \(q\in\MM_n\) with \(qh=h^2q\).
The image under \(q\) of a connected component of \(G(h)\) is
connected in \(G(h^2)\), hence lies in a unique component
of \(G(h)\). Thus \(q\) induces a map on the components of \(G(h)\). If a source component has type \(C\) and its
target has type \(D\), there are exactly \(K(C,D)\) choices for
the restriction of \(q\). These choices are independent for
different source components, since their domains are disjoint.

We therefore count functional digraphs whose vertices represent
connected components of \(G(h)\). A vertex of type \(C\) has weight \(u_C\),
and an edge from a type-\(C\) vertex to a type-\(D\) vertex has
weight \(K(C,D)\). In a rooted tree directed toward a type-\(C\)
root, an incoming child of type \(D\) contributes
\(K(D,C)R_D\). The children form an unordered set, giving
\eqref{eq:forwardR}.

\namedref{Lemma}{lem:cycle-assembly}, applied to the matrix
\(\operatorname{diag}(R_C)K\), assembles these rooted trees into
a set of directed cycles and gives
\begin{equation}\label{eq:forwardcycle}
 \exp\!\left(\sum_{\ell\geq1}\frac1\ell
       \operatorname{tr}
       \left(\bigl(\operatorname{diag}(R_C)K\bigr)^\ell\right)\right)
 =\det(I-\operatorname{diag}(R_C)K)^{-1}.
\end{equation}
The identity is formal, since every \(R_C\) has zero constant term.

For a nonnegative integer vector \(\mathbf m=(m_C)_{C\in\mathcal K_N}\),
the coefficient of \(\prod_Cu_C^{m_C}\) in~\eqref{eq:forwardcycle} is
\[
 \frac{1}{\prod_Cm_C!}
 \prod_C\left(\sum_Dm_DK(C,D)\right)^{m_C}.
\]
Indeed, distinguish the \(m_C\) vertices of each type; each then
chooses its target vertex and edge weight independently. Substituting \(u_C=x^{w_C}/\alpha_C\) supplies
the factors for the labeled components of \(G(h)\), and multiplication
by \(n!\) gives the total count. A type of size greater than \(N\)
cannot occur in a term of degree at most \(N\), which justifies the
finite cutoff in~\eqref{eq:forwarddet}.
\end{proof}

Here \(K(D,C)\) weights an incoming child in the rooted-tree
system, whereas \(K(C,D)\) weights an edge of a directed cycle.
The sum over component multiplicities in~\eqref{eq:homTotalTo}
gives a direct evaluation of the same total.

For \(C,D\in\mathcal K_N\), with cycle lengths \(\ell_C\) and
\(\ell_D\), respectively,
\[
 K(C,D)>0\quad\Longleftrightarrow\quad\ell_D\mid2\ell_C.
\]
Whenever this condition holds, there are at least \(\ell_D\) maps:
send the source cycle into the target cycle and send every attached
tree vertex backward along that cycle. If either component is a pure
cycle, these are all the maps, so \(K(C,D)=\ell_D\).

\subsection{A generating function for \texorpdfstring{\(B_n\)}{Bn}}

In computing \(H(h^2,h)\) for \(h\in\MM_n\), squaring an even-cycle
component gives two source components, and each chooses its target independently. We first apply
the preceding determinant construction to components with odd cycle
length. Euler differential operators then add each even-cycle component
with its two independent target choices.

\begin{lemma}\label{lem:homoddsupport}
For \(h\in\MM_n\), the count \(H(h^2,h)\) is positive if and
only if every cycle length of \(h\) is divisible by an odd cycle length
that also occurs in \(h\). The count \(H(h,h^2)\) is positive for every
map \(h\).
\end{lemma}

\begin{proof}
A cycle of length \(\ell\) splits under squaring into
\(\gcd(\ell,2)\) cycles of length \(\ell/\gcd(\ell,2)\).
Thus a homomorphism from its squared components to \(h\) exists exactly
when some occurring cycle length \(d\) divides
\(\ell/\gcd(\ell,2)\). If \(\ell\) is even and this chosen \(d\)
is still even, apply the same condition to \(d\). Its exponent of
\(2\) strictly decreases at each step, eventually giving an occurring
odd divisor of \(\ell\). Conversely, any occurring odd divisor of
\(\ell\) divides \(\ell/\gcd(\ell,2)\) and supplies a target cycle.
Every squared component can then be mapped independently.
For a homomorphism from \(h\) to \(h^2\), a component of cycle
length \(\ell\) can always map to one of its own squared
components, whose cycle length divides \(\ell\).
\end{proof}

Fix an integer \(N\geq1\). For \(C\in\mathcal K_N\), define the exponent
of \(2\) in its cycle length by
\(\rho_C:=\max\{r\in\mathbb Z_{\geq0}:2^r\mid\ell_C\}\), and, for integers
\(r\geq0\), put
\[
 \mathcal K_{N,r}:=\{C\in\mathcal K_N:\rho_C=r\},\qquad
 r_{\max}:=\lfloor\log_2N\rfloor.
\]
If \(\rho_C=0\), write \(S_C\) for the unique member of
\(\sigma(C)\). If \(\rho_C>0\), write
\(S_{C,0},S_{C,1}\) for its two members, retaining both occurrences
when their types coincide. A possible target of either
\(S_{C,j}\) has smaller exponent of \(2\) in its cycle length.
We therefore add the types in increasing order of \(\rho_C\).

Use independent variables \(\mathbf u=(u_C)_{C\in\mathcal K_N}\),
with weighted degree \(\deg u_C:=w_C\). On the odd-cycle types define
the matrix
\[
 L(C,D):=H(S_C,D)
 \qquad(C,D\in\mathcal K_{N,0}),
\]
and let \((T_C)_{C\in\mathcal K_{N,0}}\) be the unique family of formal
power series with zero constant terms satisfying
\begin{equation}\label{eq:reverseoddtrees}
 T_C:=u_C\exp\!\left(
       \sum_{D\in\mathcal K_{N,0}}L(D,C)T_D\right).
\end{equation}
For a type \(C\in\mathcal K_{N,r}\), \(r\geq1\), define the
Euler differential operators
\begin{equation}\label{eq:reverseEuler}
 E_{C,j}:=
 \sum_{\substack{D\in\mathcal K_N\\\rho_D<r}}
       H(S_{C,j},D)\,u_D\frac{\partial}{\partial u_D}
 \qquad(j=0,1).
\end{equation}
For a linear operator \(\mathcal A\) acting on formal series in
\(\mathbf u\), and such a series \(F\), define the operator exponential
by the following sum whenever each coefficient receives only finitely
many contributions:
\[
 (\exp\mathcal A)F:=\sum_{k\geq0}\frac{\mathcal A^kF}{k!}.
\]
Set
\begin{align}
 F_0(\mathbf u)&:=
    \det\!\left(I-\operatorname{diag}(T_C)L\right)^{-1},
       \label{eq:reverseF0}\\
 F_r(\mathbf u)&:=
    \exp\!\left(
       \sum_{C\in\mathcal K_{N,r}}u_C E_{C,0}E_{C,1}\right)
    F_{r-1}(\mathbf u)
       \qquad(1\leq r\leq r_{\max}).
       \label{eq:reverselevels}
\end{align}
Each application of the exponent in~\eqref{eq:reverselevels} increases
weighted degree, so every coefficient involves only finitely many terms.

\begin{theorem}\label{thm:reverseEuler}
For \(0\leq n\leq N\),
\begin{equation}\label{eq:reverseEulerTotal}
 B_n=n!\coef{x^n}\left.
       F_{r_{\max}}(\mathbf u)\right|_{u_C=x^{w_C}/\alpha_C}.
\end{equation}
The substitutions are made after all operators have been applied.
\end{theorem}

\begin{proof}
First consider maps whose component types all lie in
\(\mathcal K_{N,0}\). Their squared components are connected. Each
original component of type \(C\) chooses a target component of type
\(D\) and one of \(L(C,D)\) homomorphisms. On the collection of
original components these choices form a functional digraph with edge
weight \(L(C,D)\) from type \(C\) to type \(D\).

The rooted-tree equation~\eqref{eq:common-colored-trees} and
\namedref{Lemma}{lem:cycle-assembly}, with matrix \(L\) on the odd-cycle
types, give~\eqref{eq:reverseoddtrees} and~\eqref{eq:reverseF0}.
Consequently, the coefficient of \(\prod_Cu_C^{m_C}\) in \(F_0\)
is the homomorphism count for those component multiplicities divided by
\(\prod_Cm_C!\).

Suppose now that all types with exponent less than \(r\) have been
included. For a fixed multiplicity vector \(\mathbf m\) on those types and a new type
\(C\in\mathcal K_{N,r}\), put, for \(j=0,1\),
\[
 \beta_{C,j}(\mathbf m):=
       \sum_{\rho_D<r}H(S_{C,j},D)m_D.
\]
The two squared components of one occurrence of \(C\) can be mapped
independently, giving
\(\beta_{C,0}(\mathbf m)\beta_{C,1}(\mathbf m)\) choices.
Every target has smaller exponent, since its cycle length divides
\(\ell_C/2\). Thus components newly added with exponent \(r\) do not
serve as targets for one another.
By~\eqref{eq:homSupport}, previously included sources also cannot choose
these new targets, so their previously counted choices remain
unchanged.

On a monomial \(\mathbf u^{\mathbf m}\), the operator \(E_{C,j}\)
acts by multiplication by \(\beta_{C,j}(\mathbf m)\). Hence
\eqref{eq:reverselevels} multiplies that monomial by
\[
 \prod_{C\in\mathcal K_{N,r}}
   \exp\!\left(u_C\,
       \beta_{C,0}(\mathbf m)\beta_{C,1}(\mathbf m)\right).
\]
The coefficient of \(u_C^{m_C}\) is the number of choices for
\(m_C\) new components divided by \(m_C!\). Operators for distinct
types with the same exponent commute because their Euler factors involve
only variables indexed by smaller exponents. Induction now gives the coefficient
interpretation of \(F_{r_{\max}}\).

The two squared components remain distinguished by their source vertex
sets even if their types coincide. They may also choose the same actual
target component. The product of Euler operators gives precisely these
choices; it introduces no factor of \(1/2\) and no requirement of
distinct targets.

Finally, the substitution contributes
\(\prod_C\alpha_C^{-m_C}x^{\sum_Cw_Cm_C}\). Multiplication by
\(n!\) yields the class-size factor
\(n!/\prod_C(\alpha_C^{m_C}m_C!)\) for each multiplicity vector.
Every map on \(n\leq N\) points has all its component types in
\(\mathcal K_N\), completing the proof.
\end{proof}

For \(h\in\mathfrak S_n\), let \(m_\ell\) be its number of cycles of
length \(\ell\), and put \(d_\ell:=\gcd(\ell,2)\) for \(\ell\geq1\). The component formulas reduce to
\[
 H(h^2,h)=
 \prod_{\ell\geq1}
   \left(\sum_{j\mid\ell/d_\ell}j m_j\right)^{d_\ell m_\ell},
 \qquad
 H(h,h^2)=
 \prod_{\ell\geq1}
   \left(\sum_{j\mid2\ell}j m_j\right)^{m_\ell}.
\]
Here the homomorphism \(q\) ranges over \(\MM_n\). Summing these expressions with class
weight \(n!/\prod_\ell(\ell^{m_\ell}m_\ell!)\) gives the corresponding
restricted totals.

If \(q\) is also bijective, it must give an isomorphism. By
\namedref{Section}{sec:scope}, both totals then become
\(n!q_n\), the sequence \oeis{A088311}.

\subsection{Equal powers}
\label{sec:equalpowers}

For \(f^2=g^3\), define
\[
 U_n:=\NN_{ff,ggg}(n),\qquad
 r_j(c):=\#\{u\in\MM_n:u^j=c\}\quad(c\in\MM_n,\ j\geq1).
\]
Partitioning the solution pairs by their common iterate gives
\begin{equation}
 U_n=\sum_{c\in\MM_n}r_2(c)r_3(c).
 \label{eq:equalpowerroots}
\end{equation}
For each prescribed \(c\), the square- and cube-root choices are
independent. Snowden and
Howie~\cite{SnowdenHowie1982} characterize the existence of square roots.
Higgins gives constructions for square roots~\cite{Higgins1986} and
general iterative roots~\cite{Higgins1988}; further square-root
algorithms are given by Kozyra~\cite{Kozyra2018}.
We count all solutions by grouping iterates by their functional-digraph types.

To group the sum by isomorphism type, use the weights \(\omega_n(t)\)
of \namedref{Lemma}{lem:type-summation}. For a positive integer \(j\) and
\(t\in\mathfrak T_n\), define
\begin{equation}
 H_{j,n}(t):=
 \sum_{\substack{s\in\mathfrak T_n\\\operatorname{type}(h_s^j)=t}}\omega_n(s).
 \label{eq:powerhistogram}
\end{equation}

\begin{theorem}
\label{thm:equalpowerhistogram}
For every \(n\geq0\),
\begin{equation}
 U_n=\sum_{t\in\mathfrak T_n}
       \frac{H_{2,n}(t)H_{3,n}(t)}{\omega_n(t)}.
 \label{eq:equalpowerhistogram}
\end{equation}
For fixed positive integers \(a,b\), replacing \(2,3\) by \(a,b\)
counts pairs \((f,g)\in\MM_n^2\) satisfying \(f^a=g^b\).
\end{theorem}

\begin{proof}
Apply \namedref{Lemma}{lem:type-summation} with \(P(u):=u^2\) and
\(Q(u):=u^3\). Its two weighted sums are \(H_{2,n}\) and \(H_{3,n}\),
giving~\eqref{eq:equalpowerhistogram}. The lemma also gives
\(r_j(c)=H_{j,n}(t)/\omega_n(t)\) for every map \(c\) of type
\(t\), so this quotient is an integer. Replacing \(2,3\) by
\(a,b\) gives the general case.
\end{proof}

Numerical data for the specializations \(f^2=g^2\), \(f^3=g^3\),
and \(f^2=g^3\) are recorded in \oeis{A235327}, \oeiss{A235326},
and \oeiss{A235325}, respectively.

The type of an iterate of a disjoint union is the union of the
component iterate types. Thus the squaring rule~\eqref{eq:homSquareTrees}
applies componentwise, and each connected component's cube type need
be computed only once.

\begin{remark}
\label{rem:permutation-equalpowers}
If either map in \(f^2=g^3\) is bijective, so is the other.
Requiring one or both to be bijective therefore gives
\(U_n^{\mathrm{bij}}\) solutions. In the type algorithm it suffices to
retain disjoint unions of cycles, encoded by integer partitions. The
sum over common targets gives a generating function expressed as a
product over cycle lengths. This is the \((2,3)\) case of
the classical permutation-root formulas~\cite{Pavlov1981,LeanosMorenoRivera2012} and of the group-homomorphism
count in Baker and Petri~\cite[Proposition~A.1]{BakerPetri2023}.

For integers \(\ell\geq1\) and \(m\geq0\), let \(a_{\ell,m}\) and
\(b_{\ell,m}\) be the numbers of square and cube roots of a prescribed
permutation of \(X_{\ell m}\) consisting of \(m\) cycles of length \(\ell\). Set
\(a_{\ell,0}:=1\) and \(b_{\ell,0}:=1\), and set both quantities to zero
for negative second indices. The integer recurrences are
\begin{align}
 a_{\ell,m}
   &=\mathbf1_{\ell\ \mathrm{odd}}a_{\ell,m-1}
       +\ell(m-1)a_{\ell,m-2},\notag\\
 b_{\ell,m}
   &=\mathbf1_{3\nmid\ell}b_{\ell,m-1}
       +\ell^2(m-1)(m-2)b_{\ell,m-3}
       \qquad(m\geq1).
 \label{eq:permutation-local-roots}
\end{align}
Indeed, a square root combines target cycles singly, when \(\ell\)
is odd, or in pairs, with \(\ell\) interlacings for a pair.
A cube root combines them singly, when \(3\nmid\ell\), or in
triples, with \(2\ell^2\) interlacings for a triple. Selecting the
block containing one distinguished target cycle gives the recurrences.
Equivalently, taking sets of the same blocks gives
\begin{align*}
 \sum_{m\geq0}a_{\ell,m}\frac{t^m}{m!}
   &=\exp\!\left(\mathbf1_{\ell\ \mathrm{odd}}t
                   +\frac{\ell t^2}{2}\right),\\
 \sum_{m\geq0}b_{\ell,m}\frac{t^m}{m!}
   &=\exp\!\left(\mathbf1_{3\nmid\ell}t
                   +\frac{\ell^2t^3}{3}\right).
\end{align*}
The root choices for different lengths are independent, and a target
permutation of \(X_n\) with cycle multiplicities \((m_\ell)\) has class size
\(n!/\prod_\ell\ell^{m_\ell}m_\ell!\). Hence, as an identity of
formal power series,
\begin{equation}
 \sum_{n\geq0}U_n^{\mathrm{bij}}\frac{z^n}{n!}
   =\prod_{\ell\geq1}
      \left(\sum_{m\geq0}
       \frac{a_{\ell,m}b_{\ell,m}}{\ell^m m!}z^{\ell m}\right).
 \label{eq:permutation-equalpower-product}
\end{equation}

To extract coefficients through an integer cutoff \(N\geq0\), define \(c_{\ell,m}\) for integers
\(\ell\geq1\), \(m\geq0\), and \(F_0(n)\) for \(0\leq n\leq N\), by
\[
 c_{\ell,m}:=\frac{(\ell m)!}{\ell^m m!}
                a_{\ell,m}b_{\ell,m},\qquad
 F_0(n):=\mathbf1_{n=0}.
\]
Then, for \(1\leq\ell\leq N\) and \(0\leq n\leq N\), use
\begin{equation}
 F_\ell(n):=
  \sum_{m=0}^{\lfloor n/\ell\rfloor}
   \binom n{\ell m}c_{\ell,m}F_{\ell-1}(n-\ell m).
 \label{eq:permutation-equalpower-recurrence}
\end{equation}
The resulting values satisfy \(F_N(n)=U_n^{\mathrm{bij}}\).
Theorem~1.2 of Baker and Petri~\cite{BakerPetri2023} further implies
that the proportion of these pairs satisfying
\(f^2=g^3=\mathrm{id}\) tends to one as \(n\to\infty\).
\end{remark}

\section{Fibers preserved by a square}
\label{sec:preserved}

To count pairs \((f,g)\in\MM_n^2\) satisfying \(f^2=fg^2\), first fix
\(g\) and examine the partitions preserved by its square. A partition \(\mathcal P\) of \(X_n\) is \emph{preserved} by
\(h\in\MM_n\) if \(h\) sends every block of \(\mathcal P\) into a
single block of \(\mathcal P\). When this holds, write
\(F_h:\mathcal P\to\mathcal P\) for its induced block map and put
\(d_C:=|F_h^{-1}(C)|\) for \(C\in\mathcal P\).

\begin{theorem}\label{thm:766}
For every fixed map \(g\in\MM_n\),
\begin{equation}\label{eq:766}
 \#\{f:f^2=fg^2\}
 =\sum_{\mathcal P\text{ preserved by }g^2}
       \prod_{C\in\mathcal P}\fall{|C|}{d_C},
\end{equation}
where \(d_C\) is defined using \(F_{g^2}\).
\end{theorem}

\begin{proof}
Use the factorization from \namedref{Section}{sec:fibers}: let
\(\mathcal P\) be the fiber partition of \(f\), let
\(\pi:X_n\to\mathcal P\) send a point to its block, and write \(f=r\pi\).
The map \(r:\mathcal P\to X_n\) is injective because distinct fibers have
distinct image values. Cancelling this injective map from
\(r\pi r\pi=r\pi g^2\) gives \(\pi r\pi=\pi g^2\).
Hence \(\mathcal P\) is preserved by \(g^2\) and \(\pi r=F_{g^2}\).

Conversely, for a partition \(\mathcal P\) preserved by \(g^2\), choose
distinct labels \(r(B)\in F_{g^2}(B)\) for \(B\in\mathcal P\).
The \(d_C\) blocks sent to \(C\) have \(\fall{|C|}{d_C}\) such assignments.
The resulting maps \(f=r\pi\) have fiber partition \(\mathcal P\)
and satisfy the original identity.
\end{proof}

For a permutation \(g\), the map induced by \(g^2\) on a preserved
partition is surjective and hence a permutation. Thus \(d_C=1\) for every block,
and the fixed-\(g\) formula reduces to
\[
 \#\{f:f^2=fg^2\}
 =\sum_{\mathcal P\text{ preserved by }g^2}
      \prod_{C\in\mathcal P}|C|.
\]
The blocks in any cycle of the induced map have equal sizes. If
\(f\) is required to be a permutation, its fiber partition is
discrete and the equation gives \(f=g^2\); in this case
\(g\) is necessarily a permutation too.

Preservation by \(g^2\) does not imply preservation by \(g\). On
\(\{0,1,2\}\), take \(f=(0,0,2)\) and \(g=(0,2,1)\). Then \(f^2=f\)
and \(g^2\) is the identity, but \(g\) splits the fiber \(\{0,1\}\)
between two fibers. To sum~\eqref{eq:766}, we count maps \(g\) whose squares induce
a prescribed block map. With \(f\) fixed, this is an instance of Higgins's equation
\(ax^m b=c\)~\cite{Higgins1988}, with \(a=f\),
\(b=\operatorname{id}_{X_n}\), \(c=f^2\), \(m=2\), and unknown
map \(x=g\). Total enumeration still requires summing over \(f\).

\subsection{Refining the fiber partition}

The total number of solutions is
\[
 T_n:=\#\{(f,g)\in\MM_n^2:f^2=fg^2\}.
\]
For \(f\in\MM_n\) with \(k\) fibers, enumerate them as
\(B_1,\ldots,B_k\), let
\(\pi:X_n\to X_k\) record the fiber index, and write \(f=r\pi\),
where \(r:X_k\to X_n\) is injective. Put \(F:=\pi r\).
The proof of \namedref{Theorem}{thm:766} gives the equivalent conditions
\begin{equation}\label{eq:766quotient}
 \pi g^2=F\pi,\qquad \pi r=F.
\end{equation}
For \(i,j\in X_k\), refine the fibers according to the fiber
containing the image under \(g\):
\[
 C_{ij}:=\{x\in X_n:\pi(x)=i,\ \pi(g(x))=j\},
 \qquad a_{ij}:=|C_{ij}|.
\]
On each cell \(C_{ij}\), the first condition in
\eqref{eq:766quotient} becomes
\begin{equation}\label{eq:766cells}
 g(C_{ij})\subseteq C_{j,F(i)}.
\end{equation}
Conversely, condition~\eqref{eq:766cells} ensures that these are the
specified cells and that the first equation in~\eqref{eq:766quotient} holds. For fixed cells,
the number of choices of \(g\) is therefore
\[
 \prod_{i,j=1}^k a_{j,F(i)}^{a_{ij}}.
\]
In particular, a nonempty cell with an empty target makes
the contribution zero.

For nonnegative integers \(n,k\), let \(\mathcal A_{n,k}\) be the set
of nonnegative integral \(k\)-by-\(k\) matrices with total sum \(n\) and positive row
sums. The set \(\mathcal A_{0,0}\) contains the empty matrix, whereas
\(\mathcal A_{n,0}=\varnothing\) for \(n>0\).
For \(A=(a_{ij})\in\mathcal A_{n,k}\) and \(i,t\in X_k\), put
\[
 s_i:=\sum_{j=1}^k a_{ij},\qquad
 p_{it}(A):=\prod_{j=1}^k a_{jt}^{a_{ij}}.
\]
With independent variables \(u_1,\ldots,u_k\), encode the weighted
choices of distinct fiber images by
\begin{equation}\label{eq:766labelpoly}
 \Lambda_A:=\coef{u_1\cdots u_k}
       \prod_{t=1}^k
           \left(1+\sum_{i=1}^k p_{it}(A)u_i\right)^{s_t}.
\end{equation}

\begin{theorem}\label{thm:766matrix}
For \(n\geq0\),
\begin{equation}\label{eq:766matrix}
 T_n=n!\sum_{k=0}^n\frac1{k!}
       \sum_{A\in\mathcal A_{n,k}}
          \frac{\Lambda_A}{\prod_{i,j=1}^k a_{ij}!}.
\end{equation}
\end{theorem}

\begin{proof}
For the specified cell sizes, there are
\(n!/\prod_{i,j}a_{ij}!\) assignments of the point labels to the
cells. Fix a block map \(F:X_k\to X_k\), and put
\(d_t:=|F^{-1}(t)|\). The choices of \(g\) contribute
\(\prod_i p_{i,F(i)}(A)\). The \(d_t\) distinct labels
\(r(i)\) required to lie in \(B_t\) can be assigned in
\(\fall{s_t}{d_t}\) ways.

Expanding the product in~\eqref{eq:766labelpoly} and selecting each
\(u_i\) exactly once assigns its index \(i\) to one target \(t\).
If these assignments define \(F\), their coefficient is
\[
 \prod_{t=1}^k\fall{s_t}{d_t}\,
 \prod_{i=1}^k p_{i,F(i)}(A).
\]
Consequently, \(\Lambda_A\) is the sum of the required weights
over all \(F\). Every resulting map \(f=r\pi\) has exactly
the prescribed fibers, because \(r\) is injective and all rows are
nonempty. Each pair \((f,g)\) is counted once for each of the
\(k!\) enumerations of these fibers. Division by \(k!\) proves
the formula.
\end{proof}

\subsection{A generating function over induced maps}

Instead of fixing every cell size, we can count the same solutions by
rooted trees and cycles.
Fix \(F:X_k\to X_k\), and define
\(\Phi_F:X_k^2\to X_k^2\) by
\begin{equation}\label{eq:766pairmap}
 \Phi_F(i,j):=(j,F(i)).
\end{equation}
The identity \(\pi g^2=F\pi\) is equivalent to saying that the map
\(X_n\to X_k^2\) given by
\[
 x\longmapsto\bigl(\pi(x),\pi(g(x))\bigr)
\]
preserves directed edges from \(G(g)\) to \(G(\Phi_F)\).
Conversely, the first coordinate of any such edge-preserving map
recovers \(\pi\); its second coordinate is then forced to be
\(\pi g\). This is a bijection even when \(\pi\) is not surjective.

Introduce variables \(\mathbf x=(x_1,\ldots,x_k)\), where \(x_i\)
records vertices \(v\in X_n\) with \(\pi(v)=i\). Let
\(R_{ij}\) be the unique formal series with zero constant term
satisfying
\begin{equation}\label{eq:766trees}
 R_{ij}=x_i\exp\!\left(\sum_{a:F(a)=j}R_{ai}\right),
 \qquad 1\leq i,j\leq k.
\end{equation}
Their coefficients are determined successively by total degree.
Let \(\mathcal C(\Phi_F)\) be the set of directed cycles of
\(G(\Phi_F)\), each taken once, and put
\begin{equation}\label{eq:766cycleGF}
 \mathcal G_F(\mathbf x):=
 \prod_{\gamma\in\mathcal C(\Phi_F)}
       \left(1-\prod_{(i,j)\in\gamma}R_{ij}\right)^{-1}.
\end{equation}
This is also \(\det(I-\operatorname{diag}(R_{ij})P_F)^{-1}\),
where \(P_F\) is indexed by \(X_k^2\), with \((P_F)_{u,v}=1\) when
\(v=\Phi_F(u)\), and all other entries zero.

To obtain fiber partitions, we must require every color to occur and
select distinct image values. The following coefficient rule does both.

\begin{lemma}
\label{lem:color-selection}
Suppose \(\mathcal A(\mathbf x)\), with
\(\mathbf x=(x_1,\ldots,x_k)\), is an exponential generating
function for weighted labeled objects whose vertices have colors in
\(X_k\). For \(\mathbf s\in\mathbb Z_{\geq0}^k\),
\[
 a_{\mathbf s}:=n!\coef{\mathbf x^{\mathbf s}}\mathcal A(\mathbf x),
 \qquad n:=s_1+\cdots+s_k,
\]
is the total weight of objects on \(X_n\) with exactly \(s_i\)
vertices of color \(i\). Fix nonnegative integers
\(d_1,\ldots,d_k\), and put \(Z:=\{i\in X_k:d_i=0\}\).
Introduce one variable \(u_i\) for each \(i\) with \(d_i>0\),
and write \(\mathbf u\) for their collection.
For \(J\subseteq Z\), let \(\mathcal A_J(z,\mathbf u)\)
be obtained by substituting
\[
 x_i:=
 \begin{cases}
 z(1+u_i),&d_i>0,\\
 0,&i\in J,\\
 z,&i\in Z\setminus J.
 \end{cases}
\]
Then
\begin{equation}
 n!\left(\prod_{i=1}^k d_i!\right)
 \coef{z^n\prod_{d_i>0}u_i^{d_i}}
 \sum_{J\subseteq Z}(-1)^{|J|}\mathcal A_J(z,\mathbf u)
 \label{eq:color-selection}
\end{equation}
is the total weight of objects using every color, with an ordered
choice of \(d_i\) distinct vertices of each color \(i\).
The contribution for given color multiplicities \(\mathbf s\) is therefore
\(a_{\mathbf s}\prod_i\fall{s_i}{d_i}\) when every \(s_i>0\),
and zero otherwise.
\end{lemma}

\begin{proof}
For \(d_i>0\), the contribution \(x_i^{s_i}\) becomes
\(z^{s_i}(1+u_i)^{s_i}\). The coefficient of \(u_i^{d_i}\),
multiplied by \(d_i!\), is \(\fall{s_i}{d_i}\), and is zero
unless \(s_i\geq d_i\geq1\). For the other colors, setting
\(x_i=0\) excludes objects using color \(i\). Inclusion--exclusion
over \(J\subseteq Z\) consequently retains precisely the contributions
from objects in which all these colors occur. The factor \(n!\)
converts the exponential coefficient to the total weight on \(X_n\).
\end{proof}

When all \(d_i=0\), this reduces to the surjectivity formula
\begin{equation}
 n!\coef{z^n}\sum_{S\subseteq X_k}(-1)^{k-|S|}
       \mathcal A(z\mathbf1_S),
 \label{eq:color-surjectivity}
\end{equation}
where \(z\mathbf1_S\) means \(x_i=z\) for \(i\in S\)
and \(x_i=0\) otherwise. In degree \(n=k\), a surjective
coloring is a bijection.

To apply this rule to the fixed map \(F:X_k\to X_k\), define its indegrees and source set by
\[
 d_i:=|F^{-1}(i)|\quad(i\in X_k),\qquad
 L_F:=\{i\in X_k:d_i=0\},\qquad \lambda_F:=|L_F|.
\]
For each \(i\in X_k\) with \(d_i>0\), introduce a variable \(z_i\)
recording selected vertices in \(\pi^{-1}(i)\). Write \(\mathbf z\)
for these variables, and use \(x\) to record the total number of vertices.
For \(J\subseteq L_F\), make the substitutions
\begin{equation}\label{eq:766markers}
 x_i=\begin{cases}
       x(1+z_i),&d_i>0,\\
       0,&i\in J,\\
       x,&i\in L_F\setminus J,
     \end{cases}
\end{equation}
and denote the resulting specialization of \(\mathcal G_F\) by
\(\mathcal G_{F,J}(x,\mathbf z)\).
Define
\begin{equation}\label{eq:766Q}
 Q_F(n):=n!\left(\prod_{i=1}^k d_i!\right)
 \coef{x^n\prod_{d_i>0}z_i^{d_i}}
 \sum_{J\subseteq L_F}(-1)^{|J|}
       \mathcal G_{F,J}(x,\mathbf z).
\end{equation}
Let \(\mathcal Q_k\) contain one representative of each
functional-digraph type on \(k\) points, including disconnected types
and the empty type when \(k=0\). For the empty map \(F\),
\(\mathcal G_F=1\) and \(Q_F(n)=\mathbf1_{\{n=0\}}\).
The group \(\Aut(F)\) consists of permutations of \(X_k\)
commuting with \(F\).

\begin{theorem}\label{thm:766GF}
For \(n\geq0\),
\begin{equation}\label{eq:766GFtotal}
 T_n=\sum_{k=0}^n\ \sum_{F\in\mathcal Q_k}
          \frac{Q_F(n)}{|\Aut(F)|}.
\end{equation}
Every displayed quotient is a nonnegative integer. The terms with
\(k+\lambda_F>n\) vanish.
\end{theorem}

\begin{proof}
First count \(g\in\MM_n\) and \(\pi:X_n\to X_k\) satisfying
\(\pi g^2=F\pi\), without requiring \(\pi\) to be surjective. A vertex assigned the pair
\((i,j)\) has weight \(x_i\), and its incoming children are assigned
pairs \((a,i)\) with \(F(a)=j\). A rooted tree is therefore
specified by its root and an unordered set of such incoming trees,
giving~\eqref{eq:766trees}.

Along a directed cycle the assigned pairs follow \(\Phi_F\).
Applying \namedref{Lemma}{lem:cycle-assembly} to this transition gives
\eqref{eq:766cycleGF}, including cycles that traverse the same
cycle of \(\Phi_F\) several times.
Consequently, \(n!\coef{\mathbf x^{\mathbf s}}\mathcal G_F\)
is the number of pairs \((g,\pi)\) with \(|\pi^{-1}(i)|=s_i\)
and \(\sum_i s_i=n\).

Apply \namedref{Lemma}{lem:color-selection} to
\(\mathcal A:=\mathcal G_F\), assigning color \(\pi(v)\) to each vertex
\(v\in X_n\). This gives the substitutions and coefficient
extraction in~\eqref{eq:766markers}--\eqref{eq:766Q}.
For every color \(i\), the \(d_i\) selected distinct vertices
are the labels \(r(j)\) for the distinguishable indices
\(j\in F^{-1}(i)\). Thus \(Q_F(n)\) is the number of triples
\((g,\pi,r)\) in which \(\pi\) is surjective, \(r\) is
injective, and~\eqref{eq:766quotient} holds for the specified
labeled map \(F\).

Every solution \((f,g)\) with \(|\im f|=k\) has exactly
\(k!\) enumerations of its fibers. The type of \(F\) has
\(k!/|\Aut(F)|\) labeled representatives, and relabeling gives
the same value of \(Q_F(n)\) for all of them. This proves
\eqref{eq:766GFtotal}. More directly, \(\Aut(F)\) acts freely
on these triples: a relabeling that fixes
\(\pi\) is the identity because \(\pi\) is surjective.
Thus each quotient in~\eqref{eq:766GFtotal} is an integer.

Finally, each color with \(d_i>0\) requires at least \(d_i\)
points, and each of the \(\lambda_F\) other colors requires one.
Since \(\sum_i d_i=k\), a contribution requires
\(n\geq k+\lambda_F\).
\end{proof}

To evaluate~\eqref{eq:766Q}, group the sources of \(G(F)\) by their
image: for \(v\in X_k\), put
\(L_v:=\{i\in L_F:F(i)=v\}\) and \(\ell_v:=|L_v|\).
Permuting sources within each \(L_v\) preserves \(F\), fixes
every nonsource, and leaves all marker variables unchanged.
Thus \(\mathcal G_{F,J}\) depends on \(J\) only through
\(j_v:=|J\cap L_v|\), and its inclusion--exclusion sum equals
\[
 \sum_{0\leq j_v\leq\ell_v\ (v\in X_k)}
 (-1)^{\sum_vj_v}
 \left(\prod_v\binom{\ell_v}{j_v}\right)
 \mathcal G_{F,J(\mathbf j)}.
\]
Here \(\mathbf j=(j_v)_{v\in X_k}\), and \(J(\mathbf j)\) contains
any \(j_v\) sources from each \(L_v\). Only \(\prod_v(\ell_v+1)\) representative evaluations
are needed.

For \(g\) a permutation, omit the attached trees by replacing
\(R_{ij}\) by \(x_i\) in~\eqref{eq:766cycleGF}. Only permutation
types \(F\) can contribute, since \(g^2\) induces a permutation
on every preserved partition. For these types \(d_i=1\)
and \(L_F=\varnothing\), so the coefficient selection
in~\eqref{eq:766Q} selects one vertex from each fiber and needs no
inclusion--exclusion over empty fibers. The division by
\(|\Aut(F)|\) in~\eqref{eq:766GFtotal} is unchanged.

\subsection{Coefficient recurrences and further identities}
\label{sec:766recurrences}

Fix \(F:X_k\to X_k\), \(J\subseteq L_F\), and a truncation degree
\(n\geq0\), and make
the substitutions~\eqref{eq:766markers} in
\eqref{eq:766trees}--\eqref{eq:766cycleGF}.
Represent a series as
\(A(x)=\sum_{m=0}^n a_mx^m/m!\), with polynomial coefficients
in the \(z_i\), truncated modulo \(z_i^{d_i+1}\) for every \(d_i>0\).
For another such series \(B\) with exponential coefficients \(b_m\),
the product coefficients satisfy
\[
 (AB)_m=\sum_{t=0}^m\binom mt a_t b_{m-t}.
\]
Use the corresponding lower-case letters for the exponential coefficients
of other series. If \(H(0)=0\) and \(E=\exp H\), the coefficients
\(h_m,e_m\) satisfy
\begin{equation}\label{eq:766expRec}
 e_0=1,\qquad
 e_m=\sum_{t=1}^m\binom{m-1}{t-1}h_t e_{m-t}
 \qquad(m\geq1).
\end{equation}
If \(P(0)=0\) and \(U=(1-P)^{-1}\), their coefficients \(p_m,u_m\)
satisfy
\begin{equation}\label{eq:766invRec}
 u_0=1,\qquad
 u_m=\sum_{t=1}^m\binom mt p_tu_{m-t}
 \qquad(m\geq1).
\end{equation}
Write \(r_{ij,m}\) for the exponential coefficient of \(R_{ij}\),
and \(e_{ij,m}\) for that of
\(\exp(\sum_{a:F(a)=j}R_{ai})\), after the substitutions.
A root factor \(x(1+z_i)\) gives
\(r_{ij,0}=0\) and \(r_{ij,m}=m(1+z_i)e_{ij,m-1}\) for \(m\geq1\);
the other root substitutions
are handled in the same way. These recurrences determine every
coefficient required by~\eqref{eq:766Q} with integer arithmetic.

After the substitutions, each vertex contributes at most one factor
\(z_i\). A monomial of \(x\)-degree \(m\) and total degree \(e\)
in \(\mathbf z\) therefore records \(m-e\) unselected vertices.
Since the required total degree in \(\mathbf z\) is \(\sum_i d_i=k\),
terms with \(m-e>n-k\) cannot contribute
through degree \(n\) and may be discarded. This restriction is
preserved under every multiplication used above.

Some restrictions admit simpler counts. Maps \(f\) with one image point
contribute \(n^{n+1}\) to \(T_n\), while those with \(f\) a permutation
contribute \(n!\).

The idempotent restriction also has a direct count. For fixed \(g\in\MM_n\), let the connected components of
\(G(g^2)\) have sizes \(w_1,\ldots,w_c\).
Write \(e_r\) for the elementary symmetric polynomial of degree
\(r\), with \(e_0:=1\). The number of idempotent maps \(f\) with \(r\) image points
satisfying \(f^2=fg^2\) is
\begin{equation}\label{eq:766idempotent}
 e_r(w_1,\ldots,w_c)\,r^{c-r},\qquad 0\leq r\leq c.
\end{equation}
Indeed, \(fg^2=f\) makes \(f\) constant on each of these
components. Each image point is fixed by \(f\), so its component
contains no other image point. Choosing \(r\) components and one
image point in each contributes \(e_r\), and assigning the other
components to these image points contributes \(r^{c-r}\).

The construction extends to fixed powers and includes the earlier
identity \(fg=f^3\).

\begin{corollary}\label{cor:generalquotientpowers}
For fixed positive integers \(a,b\), the number of ordered pairs
\((f,g)\in\MM_n^2\) satisfying \(f^a=fg^b\) has a finite
rooted-tree and cycle expression.
\end{corollary}

\begin{proof}
Write \(f=r\pi\), where \(\pi:X_n\to X_k\) records its \(k\)
indexed fibers and \(r:X_k\to X_n\) is injective, and put
\(F:=\pi r\).
Induction gives \(f^a=rF^{a-1}\pi\), with
\(F^0:=\operatorname{id}_{X_k}\). Cancelling \(r\) shows that
the identity is equivalent to
\[
 \pi g^b=F^{a-1}\pi.
\]
Assign to each point \(v\in X_n\) the tuple
\((\pi(v),\pi(g(v)),\ldots,\pi(g^{b-1}(v)))\in X_k^b\).
The tuple assigned to \(g(v)\) is determined by the transition
\(\Phi:X_k^b\to X_k^b\) given by
\[
 \Phi(i_0,\ldots,i_{b-1})
   :=(i_1,\ldots,i_{b-1},F^{a-1}(i_0)).
\]
Conversely, an assignment following this transition recovers \(\pi\)
from its first coordinate; the successive shifts force the remaining
coordinates and then \(\pi g^b=F^{a-1}\pi\). For \(b=1\), this
transition means simply \(i_0\mapsto F^{a-1}(i_0)\).

For each tuple \(\mathbf i=(i_0,\ldots,i_{b-1})\in X_k^b\), use the
rooted-tree series
\[
 R_{\mathbf i}=x_{i_0}\exp\!\left(
             \sum_{\mathbf j:\Phi(\mathbf j)=\mathbf i}
                    R_{\mathbf j}\right),
\]
and take the same reciprocal product over the cycles of \(\Phi\)
as in~\eqref{eq:766cycleGF}. For \(i\in X_k\), put
\(s_i:=|\pi^{-1}(i)|\) and \(d_i:=|F^{-1}(i)|\). The choices
of \(r\) still have weight \(\prod_i\fall{s_i}{d_i}\). Thus the coefficient extraction,
inclusion--exclusion, and divisor \(|\Aut(F)|\) in
\namedref{Theorem}{thm:766GF} are unchanged; they refer to \(F\).
\end{proof}

The case \((a,b)=(3,1)\) specializes to \(fg=f^3\). When
\(a=b=1\), the transition is the identity and the condition
reduces to \(\pi g=\pi\), so \(g\) preserves each fiber of
\(f\) setwise.
\section{The identity \texorpdfstring{\(f^2=gfg\)}{f squared equals gfg}}
\label{sec:remaining}

\subsection{Two prescribed compositions}
\label{sec:compatiblelifts}

For \((f,g)\in\MM_n^2\) satisfying \(f^2=gfg\), put \(h:=fg\). Then
\[
 fg=h,\qquad gh=f^2,\qquad
 h^2=f(gh)=f^3.
\]
Thus \((h,f)\) satisfies the equal-power equation in
\namedref{Subsection}{sec:equalpowers}. To recover the original count, we weight
each such pair by the number of maps realizing the two prescribed
compositions. For \(f,h\in\MM_n\), define
\begin{equation}
 L(f,h):=\#\{g\in\MM_n:fg=h,\ gh=f^2\}.
 \label{eq:liftweight}
\end{equation}

\begin{lemma}
\label{lem:weightedlifts}
For every \(n\geq0\),
\begin{equation}
 \NN_{ff,gfg}(n)=
 \sum_{\substack{f,h\in\MM_n\\h^2=f^3}}L(f,h).
 \label{eq:weightedlifts}
\end{equation}
For a pair in this sum with \(f\) a permutation, \(h\) is also a
permutation and \(L(f,h)=1\).
\end{lemma}

\begin{proof}
Every solution \((f,g)\) determines \(h:=fg\), with
\(gh=gfg=f^2\) and \(h^2=f(gh)=f^3\). Conversely, the equations in
\eqref{eq:liftweight} give \(gfg=gh=f^2\). These constructions
are inverse, proving~\eqref{eq:weightedlifts}.
If \(f\) is a permutation, then \(h^2=f^3\) is bijective, so \(h\)
is bijective as well. The equation \(fg=h\) forces \(g=f^{-1}h\),
and this map satisfies
\[
 gh=f^{-1}h^2=f^{-1}f^3=f^2.
\]
Hence \(L(f,h)=1\).
\end{proof}

For \(u\in\MM_n\), write \(\operatorname{rank}u:=|u(X_n)|\);
its kernel relation is as defined in \namedref{Section}{sec:images-sources}.
The summands in~\eqref{eq:weightedlifts} have the following fiber formula.

\begin{theorem}
\label{thm:explicitlifts}
Let \(f,h\in\MM_n\) satisfy \(h^2=f^3\). Then
\begin{equation}
 L(f,h)=
 \begin{cases}
 \displaystyle\prod_{x\in X_n\setminus h(X_n)}
                  |f^{-1}(h(x))|,
       &\ker h\subseteq\ker f^2,\\[5pt]
 0,    &\text{otherwise}.
 \end{cases}
 \label{eq:explicitlifts}
\end{equation}
In particular, \(L(f,h)>0\) if and only if
\begin{equation}
 \ker h\subseteq\ker f^2
 \quad\hbox{and}\quad h(X_n)\subseteq f(X_n).
 \label{eq:liftsupport}
\end{equation}
Whenever these conditions hold,
\begin{equation}
 \operatorname{rank}(f^2)\leq\operatorname{rank}h
                         \leq\operatorname{rank}f.
 \label{eq:liftranks}
\end{equation}
\end{theorem}

\begin{proof}
Apply~\eqref{eq:prescribed-two-sided} with
\((a,b,c,d)=(f,h,h,f^2)\) and unknown map \(u=g\).
The equality condition \(ad=bc\) is precisely \(f^3=h^2\).
This gives~\eqref{eq:explicitlifts} and its
positivity criterion~\eqref{eq:liftsupport}. The kernel inclusion
gives the lower rank bound, and the image inclusion gives the upper
bound in~\eqref{eq:liftranks}.
\end{proof}

For arbitrary maps, the weights in~\eqref{eq:weightedlifts} need not
equal one. On \(\{0,1,2\}\), the maps \(f=(0,0,0)\) and
\(h=(0,0,1)\) satisfy
\(h^2=f^3\), but \(L(f,h)=0\), because every \(fg\) is constant.
If \(f=h\) is constant with value \(c\in X_n\), then
\(fg=h\) is automatic and \(gh=f^2\) is equivalent to \(g(c)=c\).
Consequently \(L(f,h)=n^{n-1}\) for \(n\geq1\).
Even when \(L(f,h)>0\), some maps satisfying \(fg=h\) do not satisfy
\(gfg=f^2\): on \(\{0,1\}\), the maps \(f=(0,0)\) and
\(g=(1,0)\) give \((fg)^2=f^3\) but \(gfg\ne f^2\).
Thus the unrestricted totals in~\eqref{eq:equalpowerroots}
and~\eqref{eq:weightedlifts} need not agree, although their permutation
restrictions both equal \(U_n^{\mathrm{bij}}\), evaluated in
\eqref{eq:permutation-equalpower-product}. If only \(g\) is required
to be bijective, use~\eqref{eq:bijective-prescribed-two-sided} for
\(L(f,h)\); see~\eqref{eq:coupled-bijective-g} for the case
\(f^2(X_n)=f(X_n)\).

\begin{lemma}
\label{lem:coupledcentral}
If \(f,g\in\MM_n\) satisfy \(f^2=gfg\), then
\[
 (fg)^2=(gf)^2=f^3,\qquad gf^3=f^3g.
\]
\end{lemma}

\begin{proof}
Multiplication of \(f^2=gfg\) on the left and on the right by
\(f\) gives the first two equalities. For the commutation relation,
\[
 gf^3=gf(gfg)=(gfg)fg=f^3g.
\]
\end{proof}

Commutation with \(f^3\) does not force \(f\) or \(g\) to preserve each
component of \(G(f^3)\). When \(f\) is the identity, any involution
\(g\) is admissible and can exchange two components.
The formulas below account for these exchanges by partitioning the
image and, for the restricted classes, pairing components or cycles.

\subsection{A formula from partitions of the image}
\label{sec:coupled-kernel}

We can also count solutions by partitioning the image according to the
values of the unknown map. Fix a map $f$ on a finite set $X$, and put
\[
 R:=f(X),\qquad Q:=f^2(X),\qquad S:=f^3(X),\qquad
 K(f):=\#\{g\colon X\to X:f^2=gfg\}.
\]
For $r\in R$ and $c\in Q$, define
\[
 s_r:=|f^{-1}(r)|,
 \qquad
 m_c:=\#\{x\in X\setminus R:f^2(x)=c\}.
\]
Let $\mathcal P_f$ be the set of partitions $\pi$ of $R$ such that $f^2$
is constant on each block.  If $B\in\pi$, write
\[
 c_B:=f^2(r)\quad(r\in B),\qquad
 W_B:=\sum_{r\in B}s_r=|f^{-1}(B)|.
\]
Thus $c_B\in S$.  For $r\in R$, let $[r]_\pi$ denote its block.
For each $\pi\in\mathcal P_f$, let $\mathcal I_f(\pi)$ consist of the
injections $\sigma\colon Q\to\pi$ satisfying
\begin{equation}
 [f(c)]_\pi=\sigma\bigl(c_{\sigma(c)}\bigr)
 \qquad(c\in Q).
 \label{eq:coupled-kernel-injection}
\end{equation}
For such an injection and $d\in S$, put
\[
 b_d:=\#\{B\in\pi\setminus\sigma(Q):c_B=d\},\qquad
 A_d:=W_{\sigma(d)}-
       \#\{q\in Q:f(q)\in\sigma(d)\}.
\]
The falling factorial has the convention of \namedref{Section}{sec:fibers}.

\begin{theorem}
\label{thm:coupled-kernel-partitions}
With the notation above,
\begin{equation}
 K(f)=
 \sum_{\pi\in\mathcal P_f}
 \sum_{\sigma\in\mathcal I_f(\pi)}
 \left(\prod_{c\in Q}W_{\sigma(c)}^{m_c}\right)
 \left(\prod_{d\in S}\fall{A_d}{b_d}\right).
 \label{eq:coupled-kernel-partitions}
\end{equation}
\end{theorem}

\begin{proof}
For a solution $g$ on $X$, let $a\colon R\to X$ be its
restriction to $R$. For $r\in R$, the required equation is
\begin{equation}
 a\bigl(f(a(r))\bigr)=f^2(r).
 \label{eq:coupled-image-restriction}
\end{equation}
Consequently, the partition $\pi$ of $R$ into the fibers of $a$
belongs to $\mathcal P_f$.
Moreover, $Q\subseteq a(R)$: for every $x\in X$,
$f^2(x)=g(f(g(x)))$, and $f(g(x))\in R$.
Thus each $c\in Q$ determines a unique block $\sigma(c)$
on which $a$ has value $c$.  Distinct values give distinct blocks, so
$\sigma$ is injective.

On the block $\sigma(c)$, equation~\eqref{eq:coupled-image-restriction}
becomes $a(f(c))=c_{\sigma(c)}$.  This is equivalent to
\eqref{eq:coupled-kernel-injection}.  The values on these blocks are now
fixed. If $B\in\pi\setminus\sigma(Q)$ and $c_B=d$, its common
value $y\in X$ under $a$ must satisfy
\[
 y\in X\setminus Q,\qquad f(y)\in\sigma(d).
\]
There are exactly $A_d$ such values.  Since $\pi$ is the partition into fibers of $a$,
different blocks must receive different values.  The available sets for
distinct $d$ are disjoint, because the blocks $\sigma(d)$ are disjoint.
Thus the free block values can be assigned in
$\prod_{d\in S}\fall{A_d}{b_d}$ ways.

It remains to extend $a$ to $X\setminus R$. For $x\in X\setminus R$,
put $c:=f^2(x)\in Q$. Then
\[
 g(f(g(x)))=c
 \quad\Longleftrightarrow\quad
 f(g(x))\in\sigma(c)
 \quad\Longleftrightarrow\quad
 g(x)\in f^{-1}(\sigma(c)).
\]
There are $W_{\sigma(c)}$ choices, independently for every such $x$.
Their product is $\prod_{c\in Q}W_{\sigma(c)}^{m_c}$.
Each resulting $g$ satisfies the equation and determines $\pi$ and
$\sigma$ uniquely. Summing the weights proves the formula.
\end{proof}

The compatibility condition~\eqref{eq:coupled-kernel-injection} also gives
a direct way to enumerate \(\mathcal I_f(\pi)\): for \(c\in Q\) and \(B\in\pi\), assigning \(\sigma(c)=B\)
forces \(\sigma(c_B)=[f(c)]_\pi\).
Apply this rule repeatedly and discard inconsistent or noninjective
partial assignments.

Choose a representative
\(f_t\in\MM_n\) of each type \(t\in\mathfrak T_n\). Since \(K(f)\) is
unchanged by relabeling, the unrestricted total is
\begin{equation}
 \NN_{ff,gfg}(n)=\sum_{t\in\mathfrak T_n}\omega_n(t)K(f_t).
 \label{eq:coupled-kernel-total}
\end{equation}
\namedref{Theorem}{thm:coupled-exponential} evaluates \(K(f)\) by summing
the weights \(L(f,h)\) with colored generating functions.

\subsection{Counting with colored functional digraphs}
\label{sec:coupled-exponential}

We evaluate~\eqref{eq:weightedlifts} using functional digraphs with
colors in \(X_n\). First allow repeated colors, then use
inclusion--exclusion to require every color in \(X_n\). On \(n\) vertices, such a coloring is bijective.
Each attached tree is indexed by its root color and its parent's color.

For fixed \(f\in\MM_n\), the transition
\((t,x)\mapsto(x,f^3(t))\) on \(X_n^2\) has the form
of~\eqref{eq:766pairmap}. \namedref{Lemmas}{lem:cycle-assembly}
and~\ref{lem:color-selection} again count the cycles and require every
color. The tree restrictions differ: for \(f^2=fg^2\),
the relation \(\pi g^2=F\pi\) gives independent edge conditions.
Here all incoming neighbors of a vertex must have colors with the same
value under \(f^2\), and every leaf contributes the fiber size in
\eqref{eq:explicitlifts}.

The kernel condition is essential. On
\(X:=\{0,1,2,3\}\), the maps
\[
 f=(0,0,1,2),\qquad h=(0,0,0,0)
\]
satisfy \(h^2=f^3\) and \(h(X)\subseteq f(X)\), but
\(f^2=(0,0,0,1)\), so \(\ker h\not\subseteq\ker f^2\).
No map \(g\) on \(X\) satisfies both \(fg=h\) and \(gh=f^2\),
since \(gh\) would be constant.

Fix a map \(f\) on \(X_n\), and put
\[
 R:=f(X_n),\qquad p:=f^2,\qquad q:=f^3,\qquad
 s_y:=|f^{-1}(y)|\quad(y\in X_n).
\]
In particular, \(s_y=0\) outside \(R\). For each \(S\subseteq X_n\),
define the attached-tree series \(T_{x,y}(z)\), for \(x,y\in X_n\), with
zero constant terms. The series \(E_{c,x}(z)\), for \(c,x\in X_n\),
count the allowed sets of children. These series satisfy
\begin{align}
 E_{c,x}(z)
   &:=\exp\left(\sum_{\substack{v\in S\\p(v)=c}}T_{v,x}(z)\right),
       &&c,x\in X_n,\label{eq:coupled-color-E}\\
 T_{x,y}(z)
   &:=z\mathbf1_{x\in S}\left(
       s_y+\mathbf1_{x\in R}
       \sum_{\substack{c\in X_n\\f(c)=y}}(E_{c,x}(z)-1)
       \right),
       &&x,y\in X_n.\label{eq:coupled-color-T}
\end{align}
The dependence on \(f,S\) is suppressed in this notation. These equations
have a unique solution with the stated constant terms: each coefficient
of \(T\) is determined by lower-degree coefficients of \(T\).

Consider the map
\[
 \Phi_f\colon X_n^2\longrightarrow X_n^2,
 \qquad \Phi_f(t,x):=(x,q(t)).
\]
Let \(\mathcal C_f(S)\) be the set of directed cycles of
\(G(\Phi_f)\) contained in \((S\cap R)^2\). Define
\begin{equation}
 \mathcal G_S(z):=
 \prod_{C\in\mathcal C_f(S)}
 \left(1-
   \prod_{(t,x)\in C}zE_{p(t),x}(z)
 \right)^{-1}.
 \label{eq:coupled-color-G}
\end{equation}
Each cycle product has zero constant term,
so every inverse is a well-defined formal series.

\begin{theorem}
\label{thm:coupled-exponential}
For every fixed map \(f\) on \(X_n\),
\begin{equation}
 K(f)=
 \coef{z^n}\sum_{S\subseteq X_n}(-1)^{n-|S|}\mathcal G_S(z).
 \label{eq:coupled-color-total}
\end{equation}
\end{theorem}

\begin{proof}
For a finite labeled set \(Y\), take a map \(h\colon Y\to Y\) and a
coloring \(\varphi\colon Y\to S\). Require
\begin{align}
 \varphi(h^2(v))&=q(\varphi(v)) &&(v\in Y),
       \label{eq:coupled-color-square}\\
 h(u)=h(v)&\ \Longrightarrow\ p(\varphi(u))=p(\varphi(v))
       &&(u,v\in Y),\label{eq:coupled-color-kernel}\\
 \varphi(h(Y))&\subseteq R.
       \label{eq:coupled-color-image}
\end{align}
Give \((h,\varphi)\) the weight
\begin{equation}
 \prod_{v\in Y\setminus h(Y)}s_{\varphi(h(v))}.
 \label{eq:coupled-color-weight}
\end{equation}
Thus every vertex of indegree zero contributes the fiber size at the
color of its parent.

For \(u\in h(Y)\), condition~\eqref{eq:coupled-color-kernel}
gives the common value
\[
 c_u:=p(\varphi(v))\quad\text{for any }v\in h^{-1}(u).
\]
Because \(q=fp\), condition~\eqref{eq:coupled-color-square}
is then equivalent to
\begin{equation}
 f(c_u)=\varphi(h(u))\qquad(u\in h(Y)).
 \label{eq:coupled-color-local}
\end{equation}
The series \(T_{x,y}\) counts an attached tree whose root has color
\(x\) and whose parent, outside the tree, has color \(y\).
If the root has no children, its contribution is \(zs_y\).
Otherwise it belongs to \(h(Y)\), so \(x\in R\). Its children
have colors \(v\) with a common value \(p(v)=c\), where \(f(c)=y\),
and form a nonempty labeled set of trees counted by
\(\sum_{v\in S:p(v)=c}T_{v,x}\). Their contribution is
\(z(E_{c,x}-1)\). Different \(c\) give disjoint possibilities,
because this set of children is nonempty. Requiring \(x\in S\)
gives exactly~\eqref{eq:coupled-color-T}.

Now consider a cyclic vertex with predecessor color \(t\), its own
color \(x\), and successor color \(y\). The predecessor forces the
common value under \(p\) to be \(p(t)\). Hence~\eqref{eq:coupled-color-local}
requires \(y=q(t)\). All off-cycle children must have colors \(v\)
satisfying \(p(v)=p(t)\), and may form an arbitrary labeled set.
The cyclic vertex with its attached trees has weight \(zE_{p(t),x}\),
with transition \((t,x)\mapsto(x,q(t))\). All cyclic colors lie in
\(S\cap R\).

Apply \namedref{Lemma}{lem:cycle-assembly} to the partial map induced by
\(\Phi_f\) on \((S\cap R)^2\), with decoration \(zE_{p(t),x}\).
Its product specialization gives~\eqref{eq:coupled-color-G}.
The lemma also counts repeated traversals of a cycle of \(\Phi_f\),
since decorations have positive size.
In particular, for \(m\geq0\),
\[
 g_m(S):=m!\coef{z^m}\mathcal G_S(z)
\]
is the total weight of the admissible pairs \((h,\varphi)\) on a
fixed labeled \(m\)-element set.

Apply~\eqref{eq:color-surjectivity} from
\namedref{Lemma}{lem:color-selection}, with \(m=n\). The sum
\[
 \sum_{S\subseteq X_n}(-1)^{n-|S|}g_n(S)
\]
retains only contributions from surjective colorings, hence from
bijections since \(|Y|=n\).
For such a coloring, put \(h':=\varphi h\varphi^{-1}\).
Conditions~\eqref{eq:coupled-color-square}--\eqref{eq:coupled-color-image}
become
\[
 (h')^2=f^3,\qquad
 \ker h'\subseteq\ker f^2,\qquad h'(X_n)\subseteq f(X_n).
\]
By \namedref{Theorem}{thm:explicitlifts}, the weight
\eqref{eq:coupled-color-weight} is \(L(f,h')\). For each \(h'\),
the \(n!\) bijections \(\varphi\) determine
\(h=\varphi^{-1}h'\varphi\). Thus inclusion--exclusion gives
\(n!K(f)\), proving~\eqref{eq:coupled-color-total}.

To evaluate the resulting generating function, fix \(S\subseteq X_n\).
For
\(c,x,y\in X_n\) and \(m\geq0\), write
\[
 t_{x,y,m}:=m!\coef{z^m}T_{x,y},\qquad
 a_{c,x,m}:=\sum_{\substack{v\in S\\p(v)=c}}t_{v,x,m},\qquad
 e_{c,x,m}:=m!\coef{z^m}E_{c,x}.
\]
Start with \(t_{x,y,0}=0\), \(e_{c,x,0}=1\), and compute by
increasing \(m\):
\begin{align}
 t_{x,y,1}&=\mathbf1_{x\in S}s_y,\notag\\
 t_{x,y,m}&=m\mathbf1_{x\in S\cap R}
       \sum_{f(c)=y}e_{c,x,m-1} &&(m\geq2),
       \label{eq:coupled-color-t-rec}\\
 e_{c,x,m}&=\sum_{j=1}^{m}\binom{m-1}{j-1}
           a_{c,x,j}e_{c,x,m-j} &&(m\geq1).
       \label{eq:coupled-color-e-rec}
\end{align}
\end{proof}

The same construction counts solutions of \(B=gAg\) for prescribed
\(A,B\in\MM_n\), without a commutation assumption. Under the product
\(u\star v:=uAv\) on \(\MM_n\), studied by Dolinka and
East~\cite{DolinkaEast2015},
this is a square-root equation in a variant of the full transformation
semigroup.

\begin{corollary}
\label{cor:coupled-prescribed-sandwich}
For prescribed maps \(A,B\) on \(X_n\), the formula in
\namedref{Theorem}{thm:coupled-exponential} counts maps \(g\in\MM_n\)
satisfying \(B=gAg\) after the substitutions
\[
 p:=B,\qquad q:=AB,\qquad R:=A(X_n),\qquad
 s_y:=|A^{-1}(y)|\quad(y\in X_n),
\]
and replacement of the condition \(f(c)=y\) by \(A(c)=y\).
For \(f\in\MM_n\) and fixed positive integers \(a,b\), taking
\(A:=f^b\), \(B:=f^a\),
and summing over the functional-digraph types of \(f\) with the weights
\(\omega_n\) gives the number of pairs satisfying \(f^a=gf^bg\).
\end{corollary}

\begin{proof}
Put \(h:=Ag\). Then \(gh=B\) and \(h^2=AB\).
Conversely, a fixed \(h\) satisfying \(h^2=AB\) admits a map
\(g\) with \(Ag=h\) and \(gh=B\) exactly when
\(\ker h\subseteq\ker B\) and \(h(X_n)\subseteq A(X_n)\).
The number of such \(g\) is then
\[
 \prod_{x\in X_n\setminus h(X_n)}|A^{-1}(h(x))|.
\]
This is \namedref{Lemma}{lem:prescribed-compositions} with
\((a,b,c,d)=(A,h,h,B)\).
The proof of \namedref{Theorem}{thm:coupled-exponential} now applies with
the stated substitutions.
Taking \(A:=f^b\) and \(B:=f^a\), and summing over the
functional-digraph types of \(f\), gives the final assertion.
\end{proof}

For prescribed maps \(A,B\) on \(X_n\), we can shorten the
inclusion--exclusion sum by collecting subsets that give the same
generating function. Retain
\(p,q,R,s_y\) from \namedref{Corollary}{cor:coupled-prescribed-sandwich},
and put \(D:=q(X_n)\). For \(c\in X_n\) and
\(0\leq j\leq n\), define
\[
 L_{c,j}:=\{x\in X_n\setminus D:p(x)=c,\ s_x=j\},
 \qquad m_{c,j}:=|L_{c,j}|.
\]
In the sums below, only nonempty sets \(L_{c,j}\) need be indexed.
For \(U\subseteq D\) and integers \(0\leq t_{c,j}\leq m_{c,j}\),
write \(\mathbf t:=(t_{c,j})\). Choose \(S\subseteq X_n\) so that
\(S\cap D=U\) and \(|S\cap L_{c,j}|=t_{c,j}\), and denote
\(\mathcal G_S(z)\) by \(\mathcal G_{U,\mathbf t}(z)\).

\begin{lemma}
\label{lem:coupled-color-grouping}
The series \(\mathcal G_{U,\mathbf t}\) is independent of the
choices within the sets \(L_{c,j}\). With this notation,
\[
 \#\{g:B=gAg\}
 =\coef{z^n}
 \sum_{U\subseteq D}
 \sum_{0\leq t_{c,j}\leq m_{c,j}}
 (-1)^{n-|U|-\sum_{c,j}t_{c,j}}
 \left(\prod_{c,j}\binom{m_{c,j}}{t_{c,j}}\right)
 \mathcal G_{U,\mathbf t}(z).
\]
\end{lemma}

\begin{proof}
Write \(\varepsilon_x:=\mathbf1_{x\in S}\) for \(x\in X_n\),
and, for \(c\in X_n\), put
\[
 N_c:=|U\cap p^{-1}(c)|+\sum_jt_{c,j}
      =|S\cap p^{-1}(c)|.
\]
If \(c\notin p(X_n)\), then \(E_{c,v}=1\). If \(x\notin D\),
there is no \(c\in p(X_n)\) with \(A(c)=x\), because
\(D=A(p(X_n))\). Consequently the defining equations give
\[
 T_{v,x}=z\varepsilon_v s_x\quad(v\in X_n),
 \qquad
 E_{c,x}=\exp(zs_xN_c).
\]
For such an \(x\), they then give
\[
 T_{x,y}=z\varepsilon_x
 \left(s_y+\sum_{A(c)=y}\bigl(\exp(zs_xN_c)-1\bigr)\right).
\]
This also holds when \(s_x=0\), since every term in the sum
then vanishes. Thus the contribution of the colors in \(L_{c,j}\)
to each incoming-tree sum depends only on \(t_{c,j}\).
The remaining tree series, whose root colors lie in \(D\), are
therefore determined by \(U\) and \(\mathbf t\), successively
by degree.

The square of the pair transition is
\((t,x)\mapsto(q(t),q(x))\). Hence both coordinates on a
directed cycle are periodic under \(q\), and in particular lie
in \(D\). The allowed cycles consequently depend only on \(U\).
Their decorations depend only on \(U,\mathbf t\) by the preceding
argument, proving the asserted independence of the series.
There are \(\prod_{c,j}\binom{m_{c,j}}{t_{c,j}}\) subsets
\(S\subseteq X_n\) with these data, all of cardinality \(|U|+\sum_{c,j}t_{c,j}\).
Grouping the inclusion--exclusion terms proves the formula.
\end{proof}

\subsection{Maps whose image is unchanged by a second iterate}
\label{sec:coupled-stable-image}

The formulas simplify when the first two images coincide. Let \(f\) be a map on a finite set \(X\), and suppose
that \(R:=f(X)\) satisfies \(f(R)=R\). Then
\(p:=f|_R\) is a permutation of \(R\), and every point of
\(X\setminus R\) maps directly to a cycle of \(G(f)\). Write \(\operatorname{Sym}(R)\) for the
group of permutations of \(R\). For \(r\in R\), put
\(s_r:=|f^{-1}(r)|\).

\begin{theorem}
\label{thm:coupled-stable-image}
Under these assumptions,
\begin{equation}
 K(f)=\sum_{\substack{y\in\operatorname{Sym}(R)\\y^2=p^3}}
       \prod_{r\in R}s_{y p^{-1}(r)}^{\,s_r-1}.
 \label{eq:coupled-stable-image}
\end{equation}
\end{theorem}

\begin{proof}
The identity \(f^2=gfg\) implies
\[
 R=f^2(X)=gfg(X)\subseteq g(R).
\]
Since \(|g(R)|\leq |R|\), we have \(g(R)=R\), and
\(v:=g|_R\) is a permutation. Restricting the identity to \(R\)
gives \(p^2=vpv\). Thus \(y:=pv\) satisfies \(y^2=p^3\).
Conversely, every such \(y\) determines \(v=p^{-1}y\) and
\(vpv=p^2\).

For \(x\in X\setminus R\) with \(f(x)=r\), the remaining condition is
\[
 v(f(g(x)))=p(r),\qquad
 f(g(x))=v^{-1}p(r)=yp^{-1}(r).
\]
Here \(v^{-1}p=pvp^{-1}=yp^{-1}\) follows from \(vpv=p^2\).
There are \(s_{yp^{-1}(r)}\) choices for \(g(x)\).
Exactly one point of \(f^{-1}(r)\), namely \(p^{-1}(r)\),
lies in \(R\), so there are \(s_r-1\) such independent choices
outside \(R\). This proves~\eqref{eq:coupled-stable-image}.
\end{proof}

The weighted square-root sum in~\eqref{eq:coupled-stable-image}
can be evaluated by treating cycles of \(p^3\) singly or in pairs. Put
\[
 w(u,v):=s_v^{\,s_{p(u)}-1}\qquad(u,v\in R).
\]
The summand associated with a square root is
\(\prod_{u\in R}w(u,y(u))\).
Write a cycle of \(p^3\) as
\(C=(c_0,\ldots,c_{\ell-1})\), with subscripts taken modulo
\(\ell\). Define
\[
 \alpha(C):=
 \begin{cases}
  \displaystyle\prod_{i=0}^{\ell-1}
      w(c_i,c_{i+(\ell+1)/2}),&\ell\text{ odd},\\[4pt]
  0,&\ell\text{ even}.
 \end{cases}
\]
For a distinct cycle \(D=(d_0,\ldots,d_{\ell-1})\) of \(p^3\)
with the same length, define
\begin{equation}
 \beta(C,D):=\sum_{t=0}^{\ell-1}
     \prod_{i=0}^{\ell-1}
       w(c_i,d_{i+t})w(d_{i+t},c_{i+1}).
 \label{eq:coupled-cycle-pair}
\end{equation}
If \(\mathcal A\) is a set of cycles of \(p^3\) of the same length, let
\(Z(\varnothing):=1\) and, for any \(C\in\mathcal A\), set
\begin{equation}
 Z(\mathcal A):=
 \alpha(C)Z(\mathcal A\setminus\{C\})+
 \sum_{D\in\mathcal A\setminus\{C\}}
   \beta(C,D)Z(\mathcal A\setminus\{C,D\}).
 \label{eq:coupled-cycle-matching}
\end{equation}
The sum in~\eqref{eq:coupled-stable-image} is the product of
\(Z(\mathcal A)\) over the different cycle lengths of \(p^3\).
Every square root \(y\) of \(p^3\) commutes with \(p^3\). It either
maps an odd cycle \(C\) to itself by
\(c_i\mapsto c_{i+(\ell+1)/2}\), contributing \(\alpha(C)\), or
alternates between two cycles \(C,D\) of the same length. Summing
over the possible alignments of \(C,D\) gives \(\beta(C,D)\).
Cycles with identical sequences \((s_u,s_{p(u)}-1)\), as \(u\)
runs around each cycle, can be collected by their multiplicities in the
recurrence.

Let \(K_{\mathrm{bij}}(f)\) be the number of permutations \(g\)
of \(X\) satisfying \(f^2=gfg\). Such a solution requires
\(\operatorname{rank}(f^2)=\operatorname{rank}f\), since
\(\operatorname{rank}(gfg)=\operatorname{rank}f\). Hence
\(K_{\mathrm{bij}}(f)=0\) unless \(f^2(X)=f(X)\). When these images coincide,
the proof of \namedref{Theorem}{thm:coupled-stable-image} shows that
\(g\) must map \(f^{-1}(r)\setminus R\) bijectively onto
\(f^{-1}(yp^{-1}(r))\setminus R\) for every \(r\in R\). Consequently
\begin{equation}
 K_{\mathrm{bij}}(f)=
 \left(\prod_{r\in R}(s_r-1)!\right)
 \#\left\{y\in\operatorname{Sym}(R):
       y^2=p^3,\quad s_{yp^{-1}(r)}=s_r\ (r\in R)\right\}.
 \label{eq:coupled-bijective-g}
\end{equation}
The recurrence~\eqref{eq:coupled-cycle-matching} evaluates this count
by taking \(w(u,v):=\mathbf1_{\{s_v=s_{p(u)}\}}\) and multiplying its
value by the factorial product in~\eqref{eq:coupled-bijective-g}.
The map \(f\) need not be a permutation: for example,
\(g=\operatorname{id}_X\) works for every idempotent map \(f\) on \(X\).

For idempotent \(f\), the permutation \(p\) is the identity.
Thus \(y\) is an involution of \(R\), and
\begin{equation}
 K(f)=\sum_{\substack{\sigma\in\operatorname{Sym}(R)\\
                         \sigma^2=\mathrm{id}_R}}
       \prod_{r\in R}s_{\sigma(r)}^{\,s_r-1}.
 \label{eq:coupled-idempotent-fixed}
\end{equation}
To sum over all labeled idempotent \(f\), define
\[
 J_n:=\#\{(f,g)\in\MM_n^2:f^2=f,\ f=gfg\}.
\]

\begin{corollary}
\label{cor:coupled-idempotent-egf}
With \(J_0:=1\),
\begin{equation}
\begin{aligned}
 \sum_{n\geq0}J_n\frac{z^n}{n!}
 &=\exp\left(
   \sum_{a\geq1}a^a\frac{z^a}{a!}
   +\frac12\sum_{a,b\geq1}
        a^b b^a\frac{z^{a+b}}{a!b!}
 \right)\\
 &=\exp\left(
   \frac{\mathcal T(z)\bigl(2+3\mathcal T(z)\bigr)}
        {2\bigl(1-\mathcal T(z)^2\bigr)}
 \right).
\end{aligned}
 \label{eq:coupled-idempotent-egf}
\end{equation}
Define
\[
 D(t):=(1-t)\exp\left(\frac{t(2+3t)}{2(1-t^2)}\right),
 \qquad d_k:=k!\coef{t^k}D(t)\quad(k\geq0).
\]
For every \(n\geq0\),
\[
 J_n=n!\coef{t^n}D(t)\exp(nt)
     =\sum_{k=0}^n\binom nk d_k n^{n-k}.
\]
\end{corollary}

\begin{proof}
For idempotent \(f\), \(G(f)\) consists of disjoint stars with
looped roots. The involution in
\eqref{eq:coupled-idempotent-fixed} either fixes one star or pairs
two stars. A fixed star of size \(a\) has \(a\) choices for its
root and \(a^{a-1}\) choices for the remaining values of \(g\),
giving \(a^a\). A pair of stars of sizes \(a,b\) has \(ab\)
choices of roots and weight \(b^{a-1}a^{b-1}\), giving
\(a^b b^a\). Taking a labeled set of these one-star and two-star
structures proves the first equality in~\eqref{eq:coupled-idempotent-egf}.

By \namedref{Lemma}{lem:abel-tree-coefficients}, the first sum in its
exponent is \(\mathcal T(z)/(1-\mathcal T(z))\). The double sum
counts nonempty maps whose vertices have two colors and whose edges
change color. Every directed cycle of such a functional digraph has even length,
and each connected component has two such colorings. The tree and cycle decomposition gives
\[
 1+\sum_{a,b\geq1}a^b b^a\frac{z^{a+b}}{a!b!}
 =\exp\left(2\sum_{\ell\geq1}
       \frac{\mathcal T(z)^{2\ell}}{2\ell}\right)
 =\frac{1}{1-\mathcal T(z)^2}.
\]
Substituting these two sums and combining their rational expressions
proves the second equality in~\eqref{eq:coupled-idempotent-egf}.

Finally, the substitution \(z=t\exp(-t)\) gives
\(\mathcal T(z)=t\) and \(dz/z=(1-t)\,dt/t\);
expanding \(\exp(nt)\) then gives the coefficient formula and
finite sum.
\end{proof}

\subsection{Maps satisfying \texorpdfstring{\(f^3=f^2\)}{f cubed equals f squared}}
\label{sec:coupled-height-two}

By \namedref{Lemma}{lem:unary23}, a map \(f\in\MM_n\) satisfies
\(f^3=f^2\) precisely when every component of \(G(f)\) is a
rooted tree of height at most two, with a loop at the root.
We now impose \(f^2=gfg\) as well and define
\[
 H_n:=\#\{(f,g)\in\MM_n^2:f^3=f^2,\ f^2=gfg\}.
\]
For such \(f\in\MM_n\), put \(R:=f(X_n)\) and
\(S:=f^2(X_n)=\Fix(f)\). For \(c\in S\), let \(T_c\) be the
vertex set of the component of \(G(f)\) containing \(c\), and put
\[
 P_c:=R\cap T_c,\qquad m_c:=|T_c\setminus R|.
\]
For \(B_c\subseteq P_c\) containing \(c\), define
\[
 a_c:=|f^{-1}(B_c)|,\qquad d_c:=|P_c\setminus B_c|.
\]

\begin{lemma}
\label{lem:coupled-height-two-fixed}
For \(f^3=f^2\),
\begin{equation}
 K(f)=
 \sum_{\substack{\sigma\in\operatorname{Sym}(S)\\
                         \sigma^2=\mathrm{id}_S}}
 \ \sum_{\substack{B_c\subseteq P_c,\ c\in B_c\\c\in S}}
 \prod_{c\in S}
     a_{\sigma(c)}^{\,m_c}
     (a_{\sigma(c)}-1)^{d_c}.
 \label{eq:coupled-height-two-fixed}
\end{equation}
\end{lemma}

\begin{proof}
By \namedref{Lemma}{lem:coupledcentral}, every solution satisfies
\(gf^3=f^3g\). Here \(f^3=f^2\) sends every vertex to the root of its component
in \(G(f)\). Consequently \(g\) sends roots to roots and components
to components.
On the roots the identity becomes \(g^2=\mathrm{id}_S\), so
\(g|_S=\sigma\) is an involution.

Define \(B_c:=\{r\in P_c:g(r)=\sigma(c)\}\), which contains
\(c\). For \(x\in T_c\), the condition \(gfg(x)=c\) is
equivalent to
\[
 g(x)\in f^{-1}(B_{\sigma(c)}).
\]
The values on \(B_c\) are forced. Each point of
\(P_c\setminus B_c\) may take any of the
\(a_{\sigma(c)}-1\) allowed values other than \(\sigma(c)\),
and each point of \(T_c\setminus R\) has
\(a_{\sigma(c)}\) choices. These choices are independent and
also suffice for the identity, proving the formula.
\end{proof}

For \(m\geq1\), fix a labeled \(m\)-element set \(X\). A marked
component is a pair \((f,B)\), where \(f\) is a map on \(X\) satisfying
\(f^3=f^2\), \(G(f)\) is connected, and \(B\subseteq f(X)\)
contains the root of \(G(f)\).
Write
\[
 a:=|f^{-1}(B)|,\qquad d:=|f(X)\setminus B|,
 \qquad s:=m-|f(X)|.
\]
For \(1\leq a\leq m\), define the polynomial in an indeterminate \(x\)
\begin{equation}
 \Delta_{m,a}(x):=
 \sum_{\substack{\text{marked components }(f,B)\\|f^{-1}(B)|=a}}
       (x-1)^d x^s.
 \label{eq:coupled-marked-definition}
\end{equation}
For a positive integer \(x\), each summand is the number of choices in
\namedref{Lemma}{lem:coupled-height-two-fixed} for an allowed target set of size \(x\).
These polynomials have the bivariate generating function
\begin{equation}
 \sum_{m\geq1}\sum_{a=1}^{m}
    \Delta_{m,a}(x)y^a\frac{z^m}{m!}
 =yz\exp\left(yz\bigl(\exp(xyz)+(x-1)\exp(xz)\bigr)\right).
 \label{eq:coupled-marked-egf}
\end{equation}
To obtain this generating function, mark the \(a\) points of
\(f^{-1}(B)\) by \(y\). The root contributes \(yz\). A leaf child of the root contributes
\(xyz\). A nonleaf child in \(B\), together with its nonempty
set of leaf children, contributes \(yz(e^{xyz}-1)\); a nonleaf
child outside \(B\) contributes \(yz(x-1)(e^{xz}-1)\).
Taking a set of the children gives an exponential whose exponent is
\[
 xyz+yz(e^{xyz}-1)+yz(x-1)(e^{xz}-1)
 =yz\bigl(e^{xyz}+(x-1)e^{xz}\bigr),
\]
proving~\eqref{eq:coupled-marked-egf}.

For \(a=m\), the condition \(f^{-1}(B)=X\) forces
\(B=f(X)\). The diagonal polynomials therefore satisfy
\[
 \sum_{k\geq1}\Delta_{k,k}(x)\frac{z^k}{k!}
 =z\exp\left(z\bigl(\exp(xz)+x-1\bigr)\right).
\]
Expanding the factor \(\exp(z\exp(xz))\) gives
\[
 \Delta_{k,k}(x)
 =k!\coef{z^{k-1}}
   \sum_{r\geq0}\frac{z^r}{r!}
      \exp\bigl((x(r+1)-1)z\bigr).
\]
Consequently, for \(k\geq1\),
\begin{equation}
 \Delta_{k,k}(x)
 =k\sum_{r=0}^{k-1}\binom{k-1}{r}
       \bigl(x(r+1)-1\bigr)^{k-1-r}.
 \label{eq:coupled-marked-diagonal}
\end{equation}

For integers \(u\geq0\) and \(j\), let
\(\left\{\begin{smallmatrix}u\\j\end{smallmatrix}\right\}\)
denote a Stirling number of the second kind, with value one for
\((u,j)=(0,0)\) and zero if \(j<0\), \(j>u\), or \(u>0=j\).

\begin{theorem}
\label{thm:coupled-height-two-total}
For \(1\leq a\leq m\), put \(u:=m-a\). Then
\begin{equation}
 \Delta_{m,a}(x)=\binom ma x^u
   \sum_{j=0}^{\min(u,a-1)}
    \left\{\begin{matrix}u\\j\end{matrix}\right\}
    (x-1)^j\fall{a}{j}\Delta_{a-j,a-j}(x).
 \label{eq:coupled-marked-stirling}
\end{equation}
For \(n\geq1\), define
\begin{align}
 \gamma_n:={}&\sum_{a=1}^{n}\Delta_{n,a}(a)\notag\\
 &+\frac12\sum_{m=1}^{n-1}\binom nm
    \sum_{a=1}^{m}\sum_{b=1}^{n-m}
       \Delta_{m,a}(b)\Delta_{n-m,b}(a).
 \label{eq:coupled-height-two-connected}
\end{align}
The generating function and coefficient recurrence are
\begin{align}
 \sum_{n\geq0}H_n\frac{z^n}{n!}
   &=\exp\left(\sum_{n\geq1}\gamma_n\frac{z^n}{n!}\right),\notag\\
 H_0&:=1,\qquad
 H_n=\sum_{m=1}^{n}\binom{n-1}{m-1}\gamma_m H_{n-m}
   \qquad(n\geq1).
 \label{eq:coupled-height-two-egf}
\end{align}
\end{theorem}

\begin{proof}
Factor the generating function in~\eqref{eq:coupled-marked-egf} as
\[
 \left(\sum_{k\geq1}\Delta_{k,k}(x)\frac{(yz)^k}{k!}\right)
 \exp\!\left((x-1)yz\bigl(\exp(xz)-1\bigr)\right).
\]
Expand the outer exponential in powers of \((x-1)yz\), and use
\[
 \frac{(\exp(xz)-1)^j}{j!}
 =\sum_{u\geq0}
   \left\{\begin{matrix}u\\j\end{matrix}\right\}
   \frac{x^u z^u}{u!}.
\]
The coefficient of \(y^a z^m\) has \(k=a-j\) and \(u=m-a\).
Multiplication by \(m!\), with
\[
 \frac{m!}{u!(a-j)!}=\binom ma\fall{a}{j},
\]
gives~\eqref{eq:coupled-marked-stirling}. The range of \(j\) follows
from \(k\geq1\) and the vanishing of the Stirling number for \(j>u\).

The involution of the roots has orbits of size one or two.
A one-component orbit on \(n\) labels contributes
\(\sum_a\Delta_{n,a}(a)\). For a two-component orbit, choose its
two label sets of sizes \(m,n-m\); if their allowed sets have sizes
\(a,b\), its contribution is
\(\Delta_{m,a}(b)\Delta_{n-m,b}(a)\).
Division by two removes the ordering of the components.
This proves~\eqref{eq:coupled-height-two-connected}. Taking a labeled
set of these orbits gives the exponential generating function, and
differentiating it gives the stated coefficient recurrence.

\end{proof}

The first values of \(J_n\) and \(H_n\) are
\[
 \begin{array}{c|r|r}
  n&J_n&H_n\\\hline
  0&1&1\\
  1&1&1\\
  2&6&6\\
  3&49&79\\
  4&554&1754\\
  5&7851&58\,271\\
  6&133\,972&2\,642\,812
 \end{array}
\]
These are the numbers of solutions to \(f^2=gfg\) with
\(f^2=f\) and \(f^3=f^2\), respectively.

\section*{Acknowledgements}

\noindent\textbf{AI Disclosure.}
ChatGPT (OpenAI) assisted with developing and reviewing mathematical
arguments, examining references, and polishing the exposition.
The author directed this work, reviewed and edited the content,
and takes full responsibility for the article.

\bibliographystyle{plainurl}
\begingroup
\small
\let\originalthebibliography\thebibliography
\renewcommand{\thebibliography}[1]{%
  \originalthebibliography{#1}\setlength{\itemsep}{1pt plus .2pt}}
\bibliography{refs_followup_rev078}
\endgroup

\clearpage
\appendix
\makeatletter
\renewcommand{\@seccntformat}[1]{\ifcsname appendixprefix@#1\endcsname\csname appendixprefix@#1\endcsname\else\csname the#1\endcsname\quad\fi}
\newcommand{\appendixprefix@section}{Appendix \thesection:\quad}
\makeatother

\section{Arithmetic operation bounds}
\label{sec:operation-bounds}

This appendix justifies the upper bounds in the ``Arithmetic complexity''
column of \namedref{Table}{tab:inventory}. Addition, subtraction, multiplication,
and exact division each count as one arithmetic operation; operand
bit lengths are stated separately. The bounds are asymptotic as the
relevant size parameter tends to infinity. Empty-set counts are returned
directly in constant time; bit estimates involving \(n\log n\) or
\(N\log N\) concern sizes at least two.

The brute-force indicator sum in \namedref{Section}{sec:scope} has
\(n^{2n}\) terms. For an identity of type \((2,3)\), each indicator
can be evaluated using at most \(5n\) map lookups and \(n\) comparisons.
Thus evaluating all indicators uses \(O(n^{2n+1})\) elementary steps,
while their term-by-term sum uses \(n^{2n}-1\) additions.

\subsection{A common bound for the sixteen counts}

\begin{theorem}
\label{thm:complete23}
Every identity of type \((2,3)\) in two maps is equivalent, after
exchanging their names if necessary, to a form in
\namedref{Table}{tab:inventory}. The formulas for arbitrary maps count all
ordered solutions on \(X_n\). Eleven of these numbers can be computed
with polynomially many arithmetic operations. All sixteen counts can be computed together
with \(O(n^4(17/3)^n)\) integer arithmetic operations on integers of
\(O(n\log n)\) bits.
\end{theorem}

\begin{proof}
The normalization in \namedref{Section}{sec:scope} accounts for
the four possible shorter words
and the eight possible longer words. The formulas in the table are proved
in their indicated sections. The exponential-coefficient recurrence
following \namedref{Lemma}{lem:unary23} has \(m\) summands in degree
\(m\). Applying it twice through degree \(n\), and including
the final factor \(n^n\), evaluates~\eqref{eq:unary23-total} in
\(O(n^2)\) operations. In~\eqref{eq:total750}, the binomial row and powers of
\(n-1\) require \(O(n)\) operations, and repeated squaring computes
the \(n\) powers \(k^n\) in \(O(n\log n)\) operations.
For the double sums~\eqref{eq:companion785-sum}, \eqref{eq:total781},
\eqref{eq:total759}, and~\eqref{eq:total752}, tabulate binomial
coefficients, the powers \(i^j\) for \(0\leq i,j\leq n\), and
the required rising factorials in \(O(n^2)\) operations.
Splitting a power of a product into powers of its factors then leaves
a bounded number of operations per summand. Each of these four sums can be evaluated in
\(O(n^2)\) operations. \namedref{Corollary}{cor:E-cost} gives the
\(O(n^3)\) bound for \(E_n\), and
\namedref{Corollary}{cor:companion772-cost} gives \(O(n^2)\) for \(f^2=fgf\).
\namedref{Corollary}{cor:775783-polynomial} gives \(O(n^5)\) for
\(f^2=f^2g\) and \(fg=fgf\), and
\namedref{Corollary}{cor:764polynomial} gives \(O(n^6)\) for \(fg=f^3\).
The homomorphism and equal-power bounds follow from
\namedref{Corollary}{cor:homCost} and
\namedref{Corollary}{cor:equalpower-cost}. The two remaining bounds follow from
\namedref{Corollary}{cor:766cost} and
\namedref{Corollary}{cor:coupled-exponential-cost}. Each evaluation has arithmetic
cost \(O(n^4(17/3)^n)\); adding the costs of a fixed number of
evaluations preserves this bound.

The bit bounds for the tree and exponential methods are proved with
those results, including \namedref{Corollary}{cor:companion772-cost}.
In the remaining polynomial cases, the displayed factorials, binomial
coefficients, powers, and coefficient recurrences have operands of size
\(n^{O(n)}\). The truncated tree series count labeled objects with at most
\(n\) vertices and weights at most a power of \(n\) per vertex. These values
also have \(O(n\log n)\) bits.
\end{proof}

Coefficient recurrences give several of the polynomial bounds.
The exponential bounds depend on the auxiliary types and coefficient
arrays used in the evaluations.
Let \(\alpha\) denote the rooted-tree growth constant \(2.955765\ldots\).
For a constant \(c>0\), \namedref{Lemma}{lem:functional-type-bound} gives
\(|\mathfrak T_n|\sim c\alpha^n/\sqrt n\), accounting for the
exponential factor in the homomorphism and equal-power bounds.
For \(f^2=gfg\), grouping a map and its allowed-color subset
together gives functional digraphs with two vertex colors.
\namedref{Lemma}{lem:two-colored-functional-types} bounds their number
by \(O((17/3)^n)\), giving the exponential factor in
\namedref{Corollary}{cor:coupled-exponential-cost}.
For \(f^2=fg^2\), the sizes of the coefficient arrays and the
number of inclusion--exclusion subsets both depend on the indegrees
of the induced map. \namedref{Lemmas}{lem:truncated-product-cost}
and~\ref{lem:weighted-functional-types} use this dependence to obtain
\(O(n^2(24/5)^n)\) operations (\namedref{Corollary}{cor:766cost}).
These upper bounds concern the evaluations given here and do not
exclude faster evaluations or simpler formulas.
\subsection{Finite sums and image restrictions}
\label{sec:operation-images}

\begin{corollary}\label{cor:E-cost}
For an integer \(N\geq2\), the numbers \(E_0,\ldots,E_N\) can be
computed with \(O(N^3)\) integer arithmetic operations using
\(O(N^2)\) stored integers of \(O(N\log N)\) bits.
\end{corollary}

\begin{proof}
Precompute binomial coefficients and powers \(i^j\),
\(0\leq i,j\leq N\), in \(O(N^2)\) operations and storage.
For integers \(0\leq i\leq a\leq N\), form
\((a^2-i)^t\) successively for \(t=0,\ldots,N-a\).
This tabulates \(P_{a,t}\) for \(a+t\leq N\)
in \(O(N^3)\) operations and \(O(N^2)\) stored integers.

Set \(E_0:=1\) and initialize the other totals to zero; the
\(a=0\) contribution vanishes for positive sizes. For a contribution
to \(E_n\), with \(1\leq n\leq N\), put \(q:=n-a-t\).
For \(q=0,\ldots,N-1\), maintain the powers
\(w_b:=b^q\) for \(0\leq b\leq\lfloor N^2/4\rfloor\).
For each \(a\geq1,t\geq0\) with \(a+t+q\leq N\), add
\[
 \binom{a+t+q}{a}\binom{t+q}{t}P_{a,t}\,w_{aq+(a-1)t}
\]
to \(E_{a+t+q}\). The required base satisfies
\(0\leq aq+(a-1)t=a(n-a)-t\leq\lfloor N^2/4\rfloor\).
Updating all powers by \(w_b\leftarrow bw_b\) and accumulating
these contributions costs \(O(N^2)\) operations per \(q\).
Thus the prefix takes \(O(N^3)\) operations and \(O(N^2)\)
stored integers. All bases are at most \(N^2\), all exponents
are \(O(N)\), and only polynomially many terms are added,
so every intermediate integer has \(O(N\log N)\) bits.
\end{proof}

For the permutation restriction~\eqref{eq:permutation-exact-image},
successive multiplication by the binomial factors of the displayed
product computes its coefficients in \(O(n^2)\) integer operations.

\begin{corollary}
\label{cor:companion772-cost}
Formula~\eqref{eq:companion772-explicit} computes
\(\NN_{ff,fgf}(n)\) with \(O(n^2)\) integer arithmetic operations
and \(O(n)\) stored integers of \(O(n\log n)\) bits for
\(n\geq2\).
\end{corollary}

\begin{proof}
Fix \(n\geq2\), and put \(r_m:=m!\coef{t^m}R_n(t)\)
for \(0\leq m\leq n\).
Equation~\eqref{eq:companion772-quadratic} gives \(r_0=r_1=0\) and
\begin{equation}
 r_m=m n^{m-1}+m r_{m-1}
       +\sum_{j=2}^{m-2}\binom mj r_jr_{m-j}
       \qquad(2\leq m\leq n).
 \label{eq:companion772-r-recurrence}
\end{equation}
For \(0\leq m\leq n\), let
\[
 e_m:=m!\coef{t^m}\exp(nR_n(t)),\qquad
 d_m:=m!\coef{t^m}(1-t-2R_n(t))^{-1}.
\]
Then \(e_0=d_0=1\), and for \(1\leq m\leq n\),
\begin{align*}
 e_m&=n\sum_{j=2}^m\binom{m-1}{j-1}r_j e_{m-j},\\
 d_m&=m d_{m-1}
       +2\sum_{j=2}^m\binom mj r_jd_{m-j}.
\end{align*}
These recurrences follow respectively by differentiating
\(\exp(nR_n)\) and multiplying the inverse by \(1-t-2R_n\).
They use \(O(n^2)\) integer operations through degree \(n\).
Binomial convolution forms the remaining products in
\eqref{eq:companion772-explicit} within the same bound. Powers of
\(n\) and the needed binomial rows are generated successively;
only a fixed number of arrays of length \(n+1\) are required.

For the bit bound, put \(\rho:=1/(8n)\).
Iteration of~\eqref{eq:companion772-quadratic} from zero has
nonnegative coefficients and is bounded by \(\rho\) at \(t=\rho\),
since
\[
 \exp(n\rho)-1+2\rho
 \leq \frac17+\frac14<1.
\]
Consequently \(R_n(\rho)\leq\rho\), and the values at \(\rho\)
of \(\exp(nR_n)\), \((1-t-2R_n)^{-1}\), and \((1+R_n)^2\)
are bounded by absolute constants. By nonnegativity, their exponential
coefficients in degree \(m\leq n\) are bounded by
\(C m!(8n)^m\) for an absolute constant \(C\).
The same bound, with a possibly larger constant, applies to their
products. The series \((1+R_n)^2\) also bounds the absolute
coefficients arising from \((1-R_n)^2\). Thus all intermediate
integers have \(O(n\log n)\) bits.
\end{proof}

The parameter \(n\) is fixed during each evaluation. Independently
computing the counts through size \(N\) uses \(O(N^3)\)
arithmetic operations.

Truncation at degree \(n\) evaluates the product in
\eqref{eq:permutation-783-total} in polynomially many integer operations.

For the permutation restriction~\eqref{eq:permutation-764-total},
compute the prefix through \(N\) by retaining \(1\leq b\leq N\)
in the first product and \(1\leq a<b\), \(a+b\leq N\), in the second.
Multiply the retained geometric factors in increasing degree. There are \(O(N^2)\) factors, each requiring \(O(N)\)
integer operations, so the prefix takes \(O(N^3)\) operations and
\(O(N)\) stored integers. Their bit lengths are \(O(N\log N)\),
since the coefficients are solution counts divided by the corresponding
factorials.

\subsection{Fiber profiles}
\label{sec:operation-profiles}

\begin{corollary}
\label{cor:shared-allocation-cost}
Suppose the weights \(W_b(\mathbf e)\) in
\namedref{Lemma}{lem:weighted-block-allocation} can each be computed with
polynomially many arithmetic operations in \(n\). Then its total
weighted sum can be computed with \(\exp(O(\sqrt n))\) arithmetic
operations. If the weight computations use \(O(n\log n)\)-bit
integers and, for some fixed constant \(C\) and \(n\geq2\),
\[
 0\leq W_b(\mathbf e)\leq n^{C(b+\sum_jje_j)},
\]
the computation uses integers of \(O(n\log n)\) bits.
\end{corollary}

\begin{proof}
Across all profiles \(\mathbf m\) of total size \(n\), the triples
of nonnegative integer vectors \((\mathbf m,\mathbf a,\mathbf e)\)
of length \(n\), with
\(\mathbf0\leq\mathbf e\leq\mathbf a\leq\mathbf m\), correspond
bijectively to partitions of \(n\) with three colors of parts: their
multiplicity vectors are
\(\mathbf m-\mathbf a\), \(\mathbf a-\mathbf e\), and
\(\mathbf e\). Their number is
\[
 p_3(n):=\coef{x^n}\prod_{j\geq1}(1-x^j)^{-3}.
\]
Each profile has at most \(n\) target blocks, so computing all values in
\eqref{eq:shared-allocation-recurrence} requires at most \(n p_3(n)\)
summands across all profiles. Precomputing factorials and binomial
coefficients, forming each profile multiplier, and evaluating each
summand all take polynomially many arithmetic operations. For every real \(\eta>0\),
\[
 p_3(n)\leq
 \exp\!\left(n\eta+3\sum_{j\geq1}
                \frac{1}{j(\exp(j\eta)-1)}\right)
 \leq\exp\!\left(n\eta+\frac{\pi^2}{2\eta}\right).
\]
Taking \(\eta=n^{-1/2}\) proves the stated operation bound.
The required vectors \(\mathbf m,\mathbf a,\mathbf e\) can also be
generated with polynomially many operations per summand.

For the bit bound, the sizes of the remaining target blocks sum to
at most \(n\), and the total size of the source groups assigned to
them is at most \(n\). A product of local weights is therefore at
most \(n^{2Cn}\). There are at most \(n^n\) assignments of actual
source blocks to targets. Thus every recurrence value is
\(n^{O(n)}\). The factorials, binomial coefficients, partition
multipliers, and sums over all profiles retain the same bound.
\end{proof}

The local weights for \(f^2=f^2g\), \(fg=fgf\), and
\(f^2=fgf\) are computable with polynomially many operations and
are at most \(n^{2(b+\sum_jje_j)}\) for \(n\geq2\).
Thus their common profile recurrence~\eqref{eq:shared-allocation-recurrence}
satisfies the operation and bit bounds of
\namedref{Corollary}{cor:shared-allocation-cost}. Their sharper polynomial evaluations are proved
in \namedref{Corollaries}{cor:companion772-cost} and~\ref{cor:775783-polynomial}.

\subsection{Tree and cycle recurrences}
\label{sec:operation-tree-cycles}

\begin{lemma}\label{lem:matrix-coefficients-cost}
Under the hypotheses of \namedref{Lemma}{lem:matrix-coefficients}, suppose
at most \(s\) entry series of \(M\) are nonzero, where
\(1\leq s\leq D^2\). The coefficients in
\eqref{eq:matrix-coefficients} can be computed in
\(O(ND^2+N^2Ds)\) integer additions and multiplications using
\(O(ND^2)\) stored integers.
\end{lemma}

\begin{proof}
For each degree pair \(1\leq j\leq m\leq N\), multiply each
nonzero entry of \(M_j\) by \(\binom mj\) and add its products with the corresponding row of
\(R_{m-j}\) directly into \(R_m\). This costs \(O(Ds)\)
operations. There are \(O(N^2)\) degree pairs. Initializing the dense coefficient matrices
costs \(O(ND^2)\), and the trace convolutions and final recurrence
cost \(O(N^2s)\) and \(O(N^2)\), respectively. Binomial rows
can be generated successively by Pascal's recurrence in \(O(N^2)\)
operations and \(O(N)\) storage. The matrices \(R_m\) and the
input coefficients therefore require \(O(ND^2)\) stored integers.
\end{proof}

\begin{theorem}\label{thm:allocation-polynomial}
Fix an integer \(N\geq1\), and suppose the weights in
\namedref{Lemma}{lem:weighted-block-allocation} have the form
\[
 W_b(\mathbf e)=\mathcal W(b,k,w),\qquad
 k:=\sum_j e_j,\quad w:=\sum_j j e_j.
\]
Given their values for \(1\leq b\leq N\) and
\(0\leq k,w\leq N\), the weighted totals on sets of sizes
\(0,\ldots,N\) can be computed in \(O(N^5)\) integer additions
and multiplications using \(O(N^3)\) stored integers.
If, for a fixed constant \(C\) and \(N\geq2\),
\[
 0\leq\mathcal W(b,k,w)\leq N^{C(b+w)},
\]
all the integers have \(O(N\log N)\) bits.
\end{theorem}

\begin{proof}
Use the integer coefficients in~\eqref{eq:allocation-coefficient-recursions}
through degree \(N\), in the order described after those recurrences.
The three recurrences in~\eqref{eq:allocation-coefficient-recursions} use
\(O(N^4)\), \(O(N^5)\), and \(O(N^5)\) operations,
respectively. Apply \namedref{Lemma}{lem:matrix-coefficients-cost} to
\eqref{eq:allocation-tree-cycle-total} with \(D=N\) and
\(s=N^2\). This adds \(O(N^5)\) operations. The arrays
\(f_{m,k,w}\), \(B_m\), and the inverse-matrix coefficients use
\(O(N^3)\) stored integers.

The coefficients of \(T_b\) and \(F\) count weighted labeled
trees and forests. In a product of local weights, each block size is
counted once as a vertex size and at most once as an incoming size.
The inverse-matrix coefficients count sequences of the decorated cycle
vertices before closing the cycle; only the first predecessor size is
additional, and it is at most \(N\). Thus products of weights are
\(N^{O(N)}\). The choices of partitions, edges, and vertex orderings
on at most \(N\) labeled points are also \(N^{O(N)}\).
All intermediate sums are nonnegative, and the binomial coefficients
have at most \(O(N)\) bits. Hence every integer has
\(O(N\log N)\) bits.
\end{proof}

\begin{corollary}\label{cor:775783-polynomial}
For a nonnegative integer \(N\), the counts for \(f^2=f^2g\) and
\(fg=fgf\), for all set sizes through \(N\), can be computed with \(O(N^5)\) integer arithmetic
operations and \(O(N^3)\) stored integers of \(O(N\log N)\) bits.
\end{corollary}

\begin{proof}
The two local weights are
\(\mathcal W(b,k,w)=\fall{b}{k}w^w\) and
\(\mathcal W(b,k,w)=\fall{w}{k}b^w\), as established in
\namedref{Subsection}{sec:allocation-polynomial}. Their weight tables can be
formed in \(O(N^3)\) operations after precomputing powers and falling
factorials, and each weight is at most \(N^{2(b+w)}\).
\namedref{Theorem}{thm:allocation-polynomial} applies.
\end{proof}

For bijective \(g\), the weights
\(\mathcal W(b,k,w)=\fall{b}{k}w!\) and
\(\mathcal W(b,k,w)=\fall{w}{k}b!\,\mathbf1_{\{w=b\}}\)
satisfy the same hypotheses. Thus the \(O(N^5)\) evaluation applies
also to these restrictions. For \(fg=fgf\), the product
\eqref{eq:permutation-783-total} gives a simpler specialized evaluation.

For \(f^2=fgf\), the local weight depends on \(k\), but not
on \(w\). The simplification in
\namedref{Subsection}{sec:companion-fgf-total} gives \(O(n^2)\)
operations per term, or \(O(N^3)\) for independently computing
all terms through size \(N\). For the two identities in
\namedref{Corollary}{cor:775783-polynomial}, the local weights depend
on both statistics, and the stated evaluation uses \(O(N^5)\) operations. Their profile recurrences
remain alternative evaluations.

\begin{corollary}\label{cor:764polynomial}
For a nonnegative integer \(N\), the numbers \(D_0,\ldots,D_N\) can be computed with \(O(N^6)\)
integer arithmetic operations and \(O(N^3)\) stored integers, each
of \(O(N\log N)\) bits.
\end{corollary}

\begin{proof}
Store each series through degree \(N\) by its integer exponential
coefficients, so products use binomial convolution. We first compute
the tree series and their powers.
For a nonnegative integer \(r\), let \(U(x)\) be a formal series
with \(U(0)=0\) and integer exponential coefficients, and put
\(P(x):=(1+U(x))^r\). For \(0\leq m\leq N\), write
\(u_m:=m!\coef{x^m}U(x)\) and \(p_m:=m!\coef{x^m}P(x)\). Then \(p_0=1\), and
\((1+U)P'=rU'P\) gives, for \(1\leq m\leq N\), the integer recurrence
\[
 p_m=\sum_{i=1}^m
 \left(r\binom{m-1}{i-1}-\binom{m-1}{i}\right)u_i p_{m-i},
\]
where a binomial coefficient outside its range is zero.
For \(1\leq a,b\leq N\) and \(0\leq m\leq N\), put
\[
 t_{a,b,m}:=m!\coef{x^m}T_{a,b},\qquad
 p_{a,b,m}:=m!\coef{x^m}(1+T_{a,b})^a.
\]
Equation~\eqref{eq:764trees} becomes
\[
 t_{a,b,m}=\sum_{c=1}^m\binom mc b^c p_{c,a,m-c}.
\]
For each degree \(m\), first compute all \(t_{a,b,m}\), and then
update the powers with exponents \(a\) and \(a-1\). These
\(N^2\) series and their powers require \(O(N^4)\) operations.
To factor the transition matrix, for \(1\leq a,b\leq N\), put
\[
 U_{a,b}(x):=\frac{x^a}{(a-1)!}(1+T_{a,b}(x))^{a-1},\qquad
 \mathcal U(x):=\operatorname{diag}(U_{a,b}(x)).
\]
Then \(\mathsf W_N=\mathcal U\mathcal V\), where the rows and columns
of the constant matrix \(\mathcal V\) are indexed by \(X_N^2\), with
\(\mathcal V_{(a,b),(d,c)}=\mathbf1_{\{d=b\}}c^a\).
For \(0\leq m\leq N\), the coefficients
\(u_{a,b,m}:=m!\coef{x^m}U_{a,b}\), also written
\(u_{q,m}\) for \(q=(a,b)\in X_N^2\), require only
\(O(N^3)\) stored integers.

For each \(q\in X_N^2\), compute column \(y^{(q)}(x)\) of
\(R(x):=(I-\mathcal U(x)\mathcal V)^{-1}\) and put
\(h^{(q)}(x):=\mathcal V y^{(q)}(x)\). For \(0\leq m<N\), write
\(y_m^{(q)}:=m!\coef{x^m}y^{(q)}(x)\) and
\(h_m^{(q)}:=m!\coef{x^m}h^{(q)}(x)\).
These vectors start with \(y_0^{(q)}=e_q\),
\(h_0^{(q)}=\mathcal V e_q\), where \(e_q\) is the unit vector
at coordinate \(q\). Process the columns one at a time.
Since \(y^{(q)}=e_q+\mathcal U h^{(q)}\), for \(1\leq m<N\),
\[
 \begin{aligned}
 y_{m,(a,b)}^{(q)}
   &=\sum_{j=a}^{m}\binom mj u_{a,b,j}h_{m-j,(a,b)}^{(q)},\\
 h_{m,(a,b)}^{(q)}
   &=\sum_{c=1}^N c^a y_{m,(b,c)}^{(q)}.
 \end{aligned}
\]
Each recurrence costs \(O(N^4)\) operations per column through
degree \(N-1\). Put \(L:=\operatorname{tr}(R\mathsf W_N')\) and
\(\ell_m:=m!\coef{x^m}L\) for \(0\leq m<N\).
Cyclicity of the trace gives
\(L=\operatorname{tr}(\mathcal V R\mathcal U')\), hence
\[
 \ell_m=\sum_q\sum_{j=0}^{m}\binom mj
             u_{q,j+1}h_{m-j,q}^{(q)}\qquad(0\leq m<N).
\]
Accumulate each column's contribution in \(O(N^2)\) operations,
then use the scalar recurrence in~\eqref{eq:matrix-coefficients}.
There are \(N^2\) columns, giving \(O(N^6)\) operations; the
coefficient tables and one column use \(O(N^3)\) stored integers.
By~\eqref{eq:764det}, the resulting coefficients are
\(D_0,\ldots,D_N\).

The tree and power coefficients count labeled branches and forests
with at most \(m\) points in degree \(m\), with weights bounded by
a power of \(N\) per point; their magnitudes are \(N^{O(m)}\). In
the inverse matrix \((I-\mathsf W_N)^{-1}
=\sum_{\ell\geq0}\mathsf W_N^\ell\), a coefficient of degree
at most \(N\) uses at most \(N\) nonempty decorations. The choices
of intermediate states, each among \(N^2\) pairs, preserve the bound
\(N^{O(N)}\). Multiplication by \(\mathcal V\) adds at most
\(N\) terms with factors bounded by \(N^N\), so the column
coefficients, trace sums, and final coefficients satisfy the same
bound. The signed intermediate sums in the power recurrence also
have this size, since each has at most \(N\) summands. Hence all
integers have \(O(N\log N)\) bits for \(N\geq2\).
\end{proof}

\subsection{Homomorphisms and iterate types}
\label{app:homomorphism-costs}

\begin{lemma}
\label{lem:hom-evaluation}
Let \(U,V\) be finite sets. For fixed maps \(s:U\to U\) and
\(t:V\to V\), the recurrence in \namedref{Theorem}{thm:homDP} computes \(H(s,t)\) using
\(O(|U||V|)\) integer additions and multiplications. In particular,
for a fixed \(g\in\MM_n\), it counts the maps commuting with
\(g\) in \(O(n^2)\) integer operations.
\end{lemma}

\begin{proof}
For each noncyclic vertex \(w\in U\), form the sums in
\eqref{eq:homDP} by adding \(Q(w,z)\), for each \(z\in V\),
to the sum indexed by \(t(z)\). Each pair consisting of a source-tree
edge and a target vertex contributes once. Evaluating~\eqref{eq:homcycle} adds at
most one length-\(\ell\) product per target vertex for a source
cycle of length \(\ell\), and the sum of the source cycle
lengths is at most \(|U|\). This gives the stated bound.
Taking \(s=t=g\) gives the commuting-map count.
\end{proof}

\Needspace{10\baselineskip}
\begin{lemma}
\label{lem:functional-type-bound}
Let \(R(x)\) be the ordinary generating function of unlabelled
rooted trees, let \(\rho\) be its radius of convergence, and set
\(\alpha:=\rho^{-1}=2.955765\ldots\). Then, for some \(c>0\),
\[
 t_n:=|\mathfrak T_n|\sim c\,\frac{\alpha^n}{\sqrt n},
 \qquad
 \sum_{m=0}^nt_m=O\!\left(\frac{\alpha^n}{\sqrt n}\right).
\]
A canonical code and the automorphism order of a map on \(n\)
vertices can be computed with \(O(n^2)\) work. All functional-digraph
types of sizes at most \(n\), with representatives and automorphism orders, can
be generated with \(O(n^{3/2}\alpha^n)\) work.
\end{lemma}

\begin{proof}
A rooted tree consists of its root and a multiset of rooted trees, so
\[
 R(x)=x\exp\left(\sum_{j\geq1}\frac{R(x^j)}j\right).
\]
Let \(C(x)\) be the ordinary generating function for connected
functional-digraph types, and put
\(F(x):=\sum_{n\geq0}t_nx^n\). The cycle and multiset
decompositions give, with \(\varphi\) denoting Euler's totient,
\[
 C(x)=\sum_{d\geq1}\frac{\varphi(d)}d
                  \log\frac1{1-R(x^d)},
 \qquad
 F(x)=\exp\left(\sum_{j\geq1}\frac{C(x^j)}j\right)
     =\prod_{j\geq1}\frac1{1-R(x^j)}.
\]
The last equality uses \(\sum_{d\mid m}\varphi(d)=m\) for
positive integers \(m\).
Otter~\cite[pp.~590--594]{Otter1948} proves that \(\rho\) is the
unique dominant singularity of \(R\), with a square-root expansion
in a slit neighborhood beginning
\[
 R(x)=1-\beta\sqrt{1-x/\rho}+O(1-x/\rho),\qquad \beta>0.
\]
The product \(K(x):=\prod_{j\geq2}(1-R(x^j))^{-1}\) is
analytic for \(|x|<\sqrt\rho\) and \(K(\rho)>0\).
The positive coefficients of \(R\) in every positive degree
exclude other zeros of \(1-R\) on \(|x|=\rho\). Thus \(F(x)=K(x)/(1-R(x))\) has leading term
\(K(\rho)\beta^{-1}(1-x/\rho)^{-1/2}\).
Extracting coefficients from the square-root expansion gives the
asserted asymptotic with \(c=K(\rho)/(\beta\sqrt\pi)\).
Summing the resulting geometric growth gives the cumulative bound.

To construct a canonical code for a functional digraph, encode its
rooted trees by balanced parentheses,
sort child-subtree codes, take the least cyclic rotation of each
component's tree-code list, and sort the component codes, with
delimiters distinguishing components. All codes have length \(O(n)\).
Use insertion sort on each child list. Two child occurrences of
sizes \(a,b\) are compared at most once, at cost
\(O(\min(a,b)+1)=O(ab)\). Summing \(ab\) over all child pairs
and vertices counts pairs of vertices in distinct child subtrees;
each pair has a unique lowest common ancestor. The total is
therefore \(O(n^2)\), including repeated isomorphic children.
Constructing the subtree codes also costs \(O(n^2)\), since the
sum of subtree sizes is at most \(n^2\). Testing all rotations
costs \(O(\ell_C|C|)\) for a component \(C\), where \(\ell_C\)
is its directed-cycle length, hence \(O(n^2)\) in total. Insertion sorting the component codes has
the same quadratic bound, using their disjoint vertex sets.
Adjacent comparisons recover multiplicities of equal children and
components, and rotation tests give the stabilizer orders in
\eqref{eq:homAut}. With factorials precomputed, the automorphism
products and construction of a representative also fit this bound.

Defrain, Porreca, and Timofeeva~\cite[Theorem~19]{DefrainPorrecaTimofeeva2024}
generate the isomorphism types with \(O(n^2)\) delay and
\(O(n)\) words of \(O(\log n)\) bits for generation alone.
Processing each output as above costs \(O(n^2t_n)\) altogether.
Summing over sizes at most \(n\) proves the claimed work bound.
\end{proof}

\Needspace{6\baselineskip}
\begin{corollary}\label{cor:homCost}
The homomorphism totals \(B_n\) and \(C_n\) in~\eqref{eq:homtotals}
can both be computed with \(O(n^{3/2}\alpha^n)\) integer arithmetic operations
on integers of \(O(n\log n)\) bits.
\end{corollary}

\begin{proof}
Generate the types and their class weights by
\namedref{Lemma}{lem:functional-type-bound}. For each representative
\(h\), construct \(h^2\) and apply \namedref{Theorem}{thm:homDP} in
both directions. \namedref{Lemma}{lem:hom-evaluation} gives
\(O(n^2)\) operations per type, hence \(O(n^2t_n)
=O(n^{3/2}\alpha^n)\) in total. The outputs can be processed
successively, retaining only \(O(n^2)\) counting integers.

Each homomorphism count is at most \(n^n\), each class size is at
most \(n!\), and the final totals are at most \(n^{2n}\).
The tree recurrence and the total sums use nonnegative integers;
their partial counts, as well as the factorials used for class
sizes, are \(n^{O(n)}\). Thus the arithmetic operands have
\(O(n\log n)\) bits.
\end{proof}

\Needspace{8\baselineskip}
\begin{corollary}
\label{cor:equalpower-cost}
For fixed positive integers \(a,b\), the number of ordered pairs
\((f,g)\in\MM_n^2\) satisfying \(f^a=g^b\) can be computed by
\namedref{Theorem}{thm:equalpowerhistogram} with \(O(n^{3/2}\alpha^n)\)
integer arithmetic operations on integers of \(O(n\log n)\) bits.
\end{corollary}

\begin{proof}
Generate the types and their class weights by
\namedref{Lemma}{lem:functional-type-bound}, then sort their canonical
codes. There are \(t_n\) codes of length \(O(n)\), and
\(\log t_n=O(n)\), so comparison sorting costs
\(O(nt_n\log t_n)=O(n^2t_n)\) work. For each representative,
compute its two prescribed iterates and their canonical codes in
\(O(n^2)\) work. Binary search locates each code in \(O(n^2)\)
work; add the original representative's labeled class size to the
corresponding sum in~\eqref{eq:powerhistogram}. The final summation
is linear in \(t_n\), giving the stated bound. This procedure
stores \(O(nt_n)\) code symbols and \(O(t_n)\) counting integers.

Each class size and each histogram value is at most \(n^n\), and
their products and the final count are at most \(n^{2n}\).
Thus all integers involved have \(O(n\log n)\) bits.
\end{proof}

\begin{corollary}
\label{cor:permutation-equalpower-cost}
For a nonnegative integer \(N\), the recurrence~\eqref{eq:permutation-equalpower-recurrence}
computes \(U_0^{\mathrm{bij}},\ldots,U_N^{\mathrm{bij}}\)
with \(O(N^2\log N)\) integer arithmetic operations, using
\(O(N)\) stored integers of \(O(N\log N)\) bits.
\end{corollary}

\begin{proof}
Precomputed factorials evaluate
the exact binomial quotients, and summing the numbers of terms over
\(\ell\) and \(n\) gives \(O(N^2\log N)\) integer arithmetic
operations. Keeping two consecutive rows uses \(O(N)\) integers.
The intermediate counting values are numbers of pairs of permutations
on at most \(N\) labels, so the operands have \(O(N\log N)\) bits.
\end{proof}

\subsection{Fibers preserved by powers}
\label{app:fiber-power-costs}

The evaluation of~\eqref{eq:766GFtotal} can be sharpened by retaining
the indegrees of the induced map in both coefficient multiplication
and type enumeration. The first lemma specializes the fast subset
convolution of Bj\"orklund, Husfeldt, Kaski, and
Koivisto~\cite{BjorklundHusfeldtKaskiKoivisto2007} to the truncated
polynomials used here.

\begin{lemma}
\label{lem:truncated-product-cost}
Let \(s\geq1\) and \(d_1,\ldots,d_s\) be positive integers, and put
\(k:=\sum_{i=1}^s d_i\) and \(M:=\prod_{i=1}^s(d_i+1)\).
The product of two polynomials in \(\mathbb Z[z_1,\ldots,z_s]\),
truncated modulo \(z_i^{d_i+1}\) for every \(1\leq i\leq s\),
can be computed with
\(O(k^2M)\) integer arithmetic operations, including exact divisions.
If the input coefficients have at most \(L\) bits, the intermediate
integers have \(O(L+k\log k)\) bits.
\end{lemma}

\begin{proof}
Write \(A=\sum_{\mathbf r}a_{\mathbf r}
\mathbf z^{\mathbf r}\) and
\(B=\sum_{\mathbf r}b_{\mathbf r}
\mathbf z^{\mathbf r}\), where
\(\mathbf r\in\prod_{i=1}^s\{0,\ldots,d_i\}\).
For vectors \(\mathbf r,\mathbf t\) in this set, put
\[
 |\mathbf r|:=\sum_i r_i,\qquad
 \mathbf r!:=\prod_i r_i!,\qquad
 \binom{\mathbf r}{\mathbf t}
   :=\prod_i\binom{r_i}{t_i}.
\]
For \(0\leq j\leq k\), form the binomial sums
\[
 \widehat a_j(\mathbf r):=
 \sum_{\substack{\mathbf t\leq\mathbf r\\
                  |\mathbf t|=j}}
 \binom{\mathbf r}{\mathbf t}
 \mathbf t!a_{\mathbf t},
 \qquad
 \widehat b_j(\mathbf r):=
 \sum_{\substack{\mathbf t\leq\mathbf r\\
                  |\mathbf t|=j}}
 \binom{\mathbf r}{\mathbf t}
 \mathbf t!b_{\mathbf t}.
\]
Multiply with respect to the total degree:
\[
 \widehat c_j(\mathbf r):=
 \sum_{\ell=0}^{j}\widehat a_\ell(\mathbf r)
                       \widehat b_{j-\ell}(\mathbf r).
\]
The required coefficients are
\[
 \coef{\mathbf z^{\mathbf r}}AB
 =\frac1{\mathbf r!}
  \sum_{\mathbf t\leq\mathbf r}
   (-1)^{|\mathbf r|-|\mathbf t|}
   \binom{\mathbf r}{\mathbf t}
   \widehat c_{|\mathbf r|}(\mathbf t).
\]

To verify this identity, write a \(k\)-element set as
\(Y=D_1\sqcup\cdots\sqcup D_s\), with \(|D_i|=d_i\).
Introduce commuting variables \(u_v\), for \(v\in Y\), with
\(u_v^2=0\), and substitute \(z_i=\sum_{v\in D_i}u_v\).
For \(S\subseteq Y\) with \(|S\cap D_i|=r_i\), the coefficient
of \(\prod_{v\in S}u_v\) in the substituted polynomial \(A\) is
\(\mathbf r!a_{\mathbf r}\), and similarly for \(B\).
The binomial sums above sum these coefficients over subsets of
each size. Inversion retains pairs of subsets whose union is \(S\).
When their sizes sum to \(|S|\), these subsets are disjoint, as
required for multiplication with \(u_v^2=0\). The displayed
factorial division therefore recovers the polynomial coefficient
and is exact.

Compute the binomial transforms and their inverses successively in
the \(s\) coordinates. For one total degree, coordinate \(i\)
costs \(O(M(d_i+1))\) operations. Since
\(\sum_i(d_i+1)\leq2k\), all \(k+1\) degrees cost
\(O(k^2M)\). The products over total degrees have the same bound.
Factorial scaling increases absolute values by at most \(k!\),
each binomial transform by at most \(2^k\), and each product
sum has at most \(k+1\) terms. Thus
\((k+1)(k!2^L)^2 2^{3k}\) bounds the intermediate absolute values,
which proves the bit bound.
\end{proof}

For an induced map \(F:X_k\to X_k\), with \(k\geq1\), retain
\(d_i\) and \(\lambda_F\) from \namedref{Section}{sec:preserved}, and put
\[
 M_F:=\prod_{i:d_i>0}(d_i+1).
\]
The coefficient array for the marker variables has \(M_F\) entries
at each fixed degree in \(x\). The next bound combines the array size
with the number of inclusion--exclusion terms and the condition
\(k+\lambda_F\leq n\).

\begin{lemma}
\label{lem:weighted-functional-types}
For every fixed nonnegative integer \(p\),
\[
 \sum_{k=1}^n
 \sum_{\substack{F\in\mathcal Q_k\\k+\lambda_F\leq n}}
      k^p2^{\lambda_F}M_F
 =O\!\left((24/5)^n\right).
\]
The implied constant may depend on \(p\).
\end{lemma}

\begin{proof}
For each functional-digraph type \(F\), write \(|F|\) for its number
of vertices and assign it the weight
\[
 w_F(x):=2^{\lambda_F}M_Fx^{|F|+\lambda_F}.
\]
Weights multiply over connected components. If \(\mathcal C\)
is the set of nonempty connected functional-digraph types, then,
summing over all finite functional-digraph types \(F\) and including
the empty type with weight one,
\[
 \mathcal H(x):=\sum_Fw_F(x)
      =\prod_{C\in\mathcal C}(1-w_C(x))^{-1}.
\]
Each factor records the multiplicity of one connected type.

To bound the connected sum, order the incoming trees at every
vertex and choose a starting point on each directed cycle.
Let \(A(x)\) be the generating function for the resulting plane
rooted trees attached to cycle vertices, with the same vertex weights. A leaf contributes
\(2x^2\), and a nonleaf with \(d\) children contributes
\(x(d+1)A(x)^d\). Consequently,
\[
 A(x)=2x^2+x\sum_{d\geq1}(d+1)A(x)^d
     =2x^2+x\bigl((1-A(x))^{-2}-1\bigr).
\]
A cycle vertex with \(d\) attached children has indegree \(d+1\),
so its series is
\[
 B(x):=x\sum_{d\geq0}(d+2)A(x)^d
       =\frac{x(2-A(x))}{(1-A(x))^2}.
\]
Every connected type has an ordered representation with the same
weight. A nonempty list of these decorated cycle vertices therefore
gives the bound
\[
 \sum_{C\in\mathcal C}w_C(x)\leq\frac{B(x)}{1-B(x)}
\]
wherever the series on the right converge.

Choose the rational radius \(r:=5/24\) so that iteration stays
strictly below \(1/4\), leaving room to absorb \(k^p\).
The right-hand side of the equation defining \(A(x)\) is increasing
in \(A\) for \(0\leq A<1\) and at \(x=r\), \(A=1/4\), equals
\[
 2r^2+r\bigl((1-1/4)^{-2}-1\bigr)
   =\frac{215}{864}<\frac14.
\]
Thus \(A(r)\leq1/4\), and
\[
 B(r)\leq\frac{35}{54},\qquad
 \sum_Cw_C(r)\leq\frac{35}{19}.
\]
Since \(d+1\leq2^d\) for \(d\geq1\), we have
\(M_C\leq2^{|C|}\), and hence
\[
 0\leq w_C(r)\leq(2r)^{|C|+\lambda_C}\leq\frac5{12}.
\]
The inequality \(-\log(1-t)\leq t/(1-t)\), valid for
real \(0\leq t<1\), gives
\[
 \log \mathcal H(r)\leq\frac{12}{7}\sum_Cw_C(r)\leq\frac{60}{19}.
\]
Applying this first to finite products and then passing to the
increasing limit also proves convergence.

The strict test at \(A=1/4\) persists for some \(R>r\), which
may also be chosen so that \(28R/9<1\) and \(2R<1\).
The same iteration gives \(A(R)\leq1/4\) and
\(B(R)\leq28R/9<1\), hence \(\mathcal H(R)<\infty\).
For every fixed \(p\), with \(k:=|F|\), the quantities
\(k^p(r/R)^{k+\lambda_F}\) are bounded independently of \(F\).
It follows that \(\sum_F|F|^p w_F(r)<\infty\), with the empty
type omitted. Finally, \(|F|+\lambda_F\leq n\) implies
\(r^{|F|+\lambda_F}\geq r^n\), so
\[
 \sum_{\substack{F\ne\varnothing\\|F|+\lambda_F\leq n}}
 |F|^p2^{\lambda_F}M_F
 \leq r^{-n}\sum_{F\ne\varnothing}|F|^p w_F(r)
 =O\!\left((24/5)^n\right).
\]
\end{proof}

\Needspace{6\baselineskip}
\begin{corollary}\label{cor:766cost}
The numbers \(T_0,\ldots,T_n\) can be computed using
\(O(n^2(24/5)^n)\) integer arithmetic operations on integers of
\(O(n\log n)\) bits.
\end{corollary}

\begin{proof}
Fix \(1\leq k\leq n\), a map \(F\in\mathcal Q_k\), and
\(J\subseteq L_F\). In
\eqref{eq:766trees}--\eqref{eq:766cycleGF}, after making
\eqref{eq:766markers}, work modulo \(x^{n+1}\) and
\(z_i^{d_i+1}\) for every \(d_i>0\).
\namedref{Lemma}{lem:truncated-product-cost} multiplies the coefficient
arrays with \(O(k^2M_F)\) operations.

Precompute the binomial coefficients through degree \(n\) in
\(O(n^2)\) work. The recurrences~\eqref{eq:766expRec}--\eqref{eq:766invRec}
determine each of the \(k^2\) rooted-tree series through degree
\(n\) using \(O(n^2)\) array products. Grouping the incoming-state
sums uses \(O(nk^2M_F)\) additional operations. The disjoint cycles
of \(\Phi_F\) contain at most \(k^2\) states in total, so their
products, reciprocals, and the final product require \(O(k^2)\)
length-\(n\) series operations. The cost for fixed \(F,J\) is
therefore \(O(n^2k^4M_F)\).

Grouping sources with the same image uses
\(C_F:=\prod_{v\in X_k}(\ell_v+1)\leq2^{\lambda_F}\) evaluations,
with \(\ell_v\) as defined after \namedref{Theorem}{thm:766GF}.
Only types with \(k+\lambda_F\leq n\) contribute.
\namedref{Lemma}{lem:weighted-functional-types} with \(p=4\)
therefore bounds the total evaluation cost by
\[
 O\!\left(n^2
 \sum_{k=1}^n
 \sum_{\substack{F\in\mathcal Q_k\\k+\lambda_F\leq n}}
       k^4 C_FM_F\right)
 =O\!\left(n^2(24/5)^n\right).
\]
Generating the types, their representatives, and their automorphism
orders costs \(O(n^{3/2}\alpha^n)\) by
\namedref{Lemma}{lem:functional-type-bound}, which is absorbed in this bound.
The same arrays give \(Q_F(m)\) for every \(0\leq m\leq n\), and the
final divisions by \(|\Aut(F)|\) are exact by \namedref{Theorem}{thm:766GF}.

Before inclusion--exclusion, the coefficients count labeled trees,
forests, or maps on at most \(n\) vertices, with at most \(n^2\)
choices for the ordered pair assigned to each vertex, and choices
of selected vertices. Their magnitudes are \(n^{O(n)}\).
\namedref{Lemma}{lem:truncated-product-cost} preserves the
\(O(n\log n)\) bit bound during coefficient multiplication.
The binomial multiplicities, signed inclusion--exclusion sums,
factorials, automorphism orders, and sum over
\(O(\alpha^n/\sqrt n)\) types
preserve the \(O(n\log n)\) bit bound.
\end{proof}

\begin{corollary}
\label{cor:generalquotientpowers-cost}
For fixed positive integers \(a,b\), the construction in
\namedref{Corollary}{cor:generalquotientpowers} computes the numbers of
ordered pairs \((f,g)\in\MM_m^2\) satisfying \(f^a=fg^b\),
for \(0\leq m\leq n\), with
\(O(n^2(24/5)^n)\) integer arithmetic operations on integers of
\(O(n\log n)\) bits. The constants may depend on \(a\) and \(b\).
\end{corollary}

\begin{proof}
For \(F\in\mathcal Q_k\), the tuple transition in
\namedref{Corollary}{cor:generalquotientpowers} acts on \(X_k^b\),
which has \(k^b\) states. The same coefficient recurrences and cycle
products therefore cost \(O(n^2k^{b+2}M_F)\) per subset
\(J\subseteq L_F\). The condition
\(k+\lambda_F\leq n\) and the coefficient-array sizes are
unchanged. \namedref{Lemma}{lem:weighted-functional-types} with
\(p=b+2\) bounds the total by \(O(n^2(24/5)^n)\).
For fixed \(b\), the assigned tuples and vertex selections still
give coefficients of magnitude \(n^{O(n)}\), so the bit bound
follows as in \namedref{Corollary}{cor:766cost}.
\end{proof}

For \((a,b)=(3,1)\), the identity becomes \(fg=f^3\);
\namedref{Corollary}{cor:764polynomial} gives the sharper bound.

\subsection{The identity \texorpdfstring{\(f^2=gfg\)}{f squared equals gfg}}
\label{app:coupled-operation-bounds}

To bound the total over all maps, we combine the sum over isomorphism
types with inclusion--exclusion. Two vertex colors record membership
in the allowed subset.

\begin{lemma}
\label{lem:two-colored-functional-types}
For each integer \(m\geq0\), let \(\mathfrak T_m^{(2)}\) contain
one representative of each isomorphism type of pairs
\((f,S)\), where \(f:X_m\to X_m\) and \(S\subseteq X_m\),
and isomorphisms preserve both the map and the subset. Then
\[
 \sum_{m=0}^{n}|\mathfrak T_m^{(2)}|
   =O\!\left((17/3)^n\right).
\]
All these types, representatives, and automorphism orders can be
generated with \(O(n^3(17/3)^n)\) work on integers of
\(O(n\log n)\) bits.
\end{lemma}

\begin{proof}
Distinguish membership in \(S\) by two vertex colors.
Let \(u_m\) be the number of rooted-tree types on \(m\geq1\)
vertices with these colors, and put
\(U(x):=\sum_{m\geq1}u_mx^m\). A root has two color choices
and its children form a multiset, so
\[
 U(x)=2x\prod_{m\geq1}(1-x^m)^{-u_m}
     =2x\exp\left(\sum_{j\geq1}\frac{U(x^j)}j\right).
\]
Ordering the children gives the coefficientwise majorant
\[
 P(x):=\frac{1-\sqrt{1-8x}}2,
 \qquad P(x)=\frac{2x}{1-P(x)}.
\]
For \(0\leq z\leq1/32\), we have
\(P(z)\leq P(1/32)=(2-\sqrt3)/4<1/14\), and hence
\(P(z)\leq28z/13\).

Choose \(r:=3/17\) so that the higher arguments \(r^j\),
\(j\geq2\), lie within the plane-tree bound and the rooted-tree
iteration stays below one. Since \(r^2<1/32\),
\[
 \sum_{j\geq2}\frac{U(r^j)}j
 \leq\frac{28}{13}
       \left(\frac{r^2}{2}+\frac{r^3}{3(1-r)}\right)
 =\frac{144}{3757}.
\]
The same bound holds for trees of bounded height. Iterating the
rooted-tree equation from zero stays below one, because
\[
 2r\exp\left(1+\frac{144}{3757}\right)
 <\frac6{17}\frac{87}{32}\frac{3757}{3613}
 =\frac{57681}{57808}<1.
\]
Here \(\exp(1)<87/32\) and
\(\exp(t)\leq(1-t)^{-1}\) for \(0\leq t<1\), by the
power series. Passing to the increasing limit proves
\(U(r)\leq57681/57808<1\).

Let \(c_m\) be the number of connected functional-digraph types
on \(m\geq1\) vertices with these colors. Choosing a starting vertex of the directed cycle gives
a nonempty ordered list of rooted-tree types, so
\[
 C(r):=\sum_{m\geq1}c_mr^m
   \leq\frac{U(r)}{1-U(r)}<\infty.
\]
Arbitrary types are multisets of connected types. Thus
\[
 \mathcal H(x):=\sum_{m\geq0}|\mathfrak T_m^{(2)}|x^m
      =\prod_{m\geq1}(1-x^m)^{-c_m},
 \qquad
 \log \mathcal H(r)\leq\frac{C(r)}{1-r}<\infty.
\]
The last inequality can first be applied to finite products.
It follows that
\[
 \sum_{m=0}^{n}|\mathfrak T_m^{(2)}|
 \leq r^{-n}\mathcal H(r)
 =O\!\left((17/3)^n\right).
\]

For generation, use canonical codes that include the root color.
Generate rooted-tree types by increasing size, listing the child
types in nondecreasing order. During the search for size \(m\),
each partial child list is itself the child list of a rooted-tree
type of size at most \(m\). Sort the available child types by size
and then code, so only extensions fitting the remaining size are
considered. Codes have length \(O(n)\); generation and sorting
through size \(n\) therefore use
\(O(n^2\sum_{m\leq n}u_m)\) work.

Next generate ordered lists of rooted-tree types of total size at
most \(n\), retaining only the least cyclic rotation of each.
A list of size \(m\) represents a connected type, with at most
\(m\) different lists per type. Partial lists are already lists
of size at most \(n\). There are therefore at most
\(n\sum_{m\leq n}c_m\) lists to process, and comparing all
rotations of a length-\(O(n)\) code takes \(O(n^2)\) work.
Finally, generate nondecreasing lists of connected types; each
partial list is a functional-digraph type of size at most \(n\).
This costs \(O(n\sum_{m\leq n}|\mathfrak T_m^{(2)}|)\) work.

Automorphism orders use the child-multiplicity, cycle-rotation,
and component-multiplicity formulas from \namedref{Section}{sec:hom},
with the colors included in each code. Representatives are recovered
by assigning vertex labels to these codes. Both computations fit
within the preceding bounds. Since
\(u_m\leq c_m\leq|\mathfrak T_m^{(2)}|\), the total is
\(O(n^3(17/3)^n)\). Automorphism orders and factorials are at
most \(n!\), giving the stated bit bound.
\end{proof}

\Needspace{12\baselineskip}
\begin{corollary}
\label{cor:coupled-exponential-cost}
For a fixed map \(f\) on \(X_n\), use \(D\) and
\(m_{c,j}\) from \namedref{Lemma}{lem:coupled-color-grouping} with
\(A=f\), \(B=f^2\). Then \(K(f)\) can be computed with
\[
 O\!\left(n^4 2^{|D|}\prod_{c,j}(m_{c,j}+1)\right)
\]
integer arithmetic operations and \(O(n^3)\) stored integers.
The unrestricted totals \(\NN_{ff,gfg}(m)\), \({0\leq m\leq n}\),
can be computed with \(O(n^4(17/3)^n)\) integer arithmetic operations.
All integers involved have \(O(n\log n)\) bits.
\end{corollary}

\begin{proof}
Fix \(S\subseteq X_n\) and put \(p:=f^2\). Use the exponential
coefficients \(t_{x,y,m}\), \(a_{c,x,m}\), and \(e_{c,x,m}\)
from the proof of \namedref{Theorem}{thm:coupled-exponential}, for
\(x,y,c\in X_n\) and \(0\leq m\leq n\). Compute them by increasing
degree with~\eqref{eq:coupled-color-t-rec}--\eqref{eq:coupled-color-e-rec}.
The coefficients through degree \(n\) of the \(n^2\) series
\(E_{c,x}\) require \(O(n^4)\) arithmetic operations in total.
The coefficients \(t_{x,y,m}\) and sums \(a_{c,x,m}\) can be
formed in \(O(n^3)\) operations by grouping indices under \(f\)
and \(p\).

Products use binomial convolution of exponential coefficients,
and reciprocals use~\eqref{eq:766invRec}.
The cycles of \(\Phi_f\) contain at most \(n^2\) ordered pairs.
Their products, reciprocals, and the final product in
\eqref{eq:coupled-color-G}, truncated at degree \(n\), therefore
require another \(O(n^4)\) operations. Cycles longer than \(n\)
may be omitted. \namedref{Lemma}{lem:coupled-color-grouping} requires
\(2^{|D|}\prod_{c,j}(m_{c,j}+1)\) evaluations of these series.
Since \(m+1\leq2^m\) and \(\sum_{c,j}m_{c,j}=n-|D|\),
this is at most \(2^n\). Process the choices one at a time,
using \(O(n^3)\) stored integers, and combine them with the
binomial multiplicities in that lemma. The final integer
inclusion--exclusion sum is divided exactly by \(n!\).

For the unrestricted total, take \(f\in\MM_n\) and
\(S\subseteq X_n\), write \(\mathcal G_{f,S}\) for the series
\(\mathcal G_S\) associated with \(f\), and put
\[
 \Aut(f,S):=\{\tau\in\mathfrak S_n:
                       \tau f=f\tau,\ \tau(S)=S\}.
\]
Simultaneous relabeling of \(f\) and \(S\) preserves this series.
Grouping the sum over labeled pairs \((f,S)\) by their
two-colored types therefore gives
\[
 \NN_{ff,gfg}(n)
 =\coef{z^n}
  \sum_{(f,S)\in\mathfrak T_n^{(2)}}
   (-1)^{n-|S|}\frac{n!}{|\Aut(f,S)|}\mathcal G_{f,S}(z).
\]
For \(0\leq m\leq n\), compute
\(g_m(f,S):=m!\coef{z^m}\mathcal G_{f,S}\) with the integer
recurrences above. Sum the signed integers
\((n!/|\Aut(f,S)|)g_n(f,S)\), and then divide the whole sum
by \(n!\). This final division is exact because its numerator
equals \(n!\NN_{ff,gfg}(n)\).

\namedref{Lemma}{lem:two-colored-functional-types} generates the types
and their automorphism orders within \(O(n^3(17/3)^n)\) work.
Each type on \(m\leq n\) points requires \(O(m^4)\)
operations for its coefficient. The cumulative type bound thus gives
\(O(n^4(17/3)^n)\) operations for all totals through size \(n\).

For coefficient degrees \(m\leq n\), the intermediate series count
colored trees, forests, or functional digraphs with at most \(n\)
vertices, and fiber weights at most \(n\) per leaf. For example,
\({g_m(f,S)\leq m^m n^{2m}}\). The analogous rooted-forest bound is
\((m+1)^m n^{2m}\); factorials and binomial factors are also
\(n^{O(n)}\). The grouped inclusion--exclusion multiplicities sum
to \(2^n\), the class weights are at most \(n!\), and there are
\(\exp(O(n))\) colored types. Hence all intermediate integers
and signed partial sums have \(O(n\log n)\) bits.
\end{proof}

\Needspace{14\baselineskip}
\begin{corollary}
\label{cor:coupled-prescribed-sandwich-cost}
For prescribed maps \(A,B\) on \(X_n\), with
\(D,m_{c,j}\) as in \namedref{Lemma}{lem:coupled-color-grouping}, the
number of maps \(g\) on \(X_n\) satisfying \(B=gAg\)
can be computed with
\[
 O\!\left(n^4 2^{|D|}\prod_{c,j}(m_{c,j}+1)\right)
\]
integer arithmetic operations and \(O(n^3)\) stored integers of
\(O(n\log n)\) bits. For fixed positive integers \(a,b\), all
counts of pairs \((f,g)\in\MM_m^2\) satisfying \(f^a=gf^bg\),
for \(0\leq m\leq n\), can be computed with \(O(n^4(17/3)^n)\) arithmetic operations.
The constants may depend on \(a,b\).
\end{corollary}

\begin{proof}
Apply the substitutions of \namedref{Corollary}{cor:coupled-prescribed-sandwich}.
They preserve the number of indexed series and the form of the tree
and cycle recurrences. Their coefficients still count colored
structures with fiber weights at most \(n\).
\namedref{Lemma}{lem:coupled-color-grouping} gives the number of
evaluations, and the operation, storage, and bit bounds follow
as in \namedref{Corollary}{cor:coupled-exponential-cost}.

For fixed \(a,b\), forming \(A=f^b\) and \(B=f^a\) takes
\(O(n)\) operations per \(f\). Simultaneous relabeling of
\(f,S\) also preserves the modified series, so the same summation
over two-colored types applies. Its evaluation has the stated bound.
\end{proof}

For the idempotent restriction, retain \(D(t)\) and
\(d_m:=m!\coef{t^m}D(t)\), for \(m\geq0\), from
\namedref{Corollary}{cor:coupled-idempotent-egf}. These coefficients satisfy
\(d_0=1\), \(d_1=0\), \(d_2=2\), \(d_3=4\), and
\[
\begin{aligned}
 d_m={}&2(m-1)^2d_{m-2}
       +2(m-1)(m-2)d_{m-3}\\
      &+(m-1)(m-2)(m-3)(5-m)d_{m-4}
       \qquad(m\geq4).
\end{aligned}
\]
Indeed, logarithmic differentiation gives
\[
 \frac{D'(t)}{D(t)}
 =\frac{t(t^2+2t+2)}{(1-t^2)^2},
\]
and comparison of exponential coefficients in
\((1-t^2)^2D'=t(t^2+2t+2)D\) gives the recurrence.
It computes \(d_0,\ldots,d_n\) with \(O(n)\) integer arithmetic
operations. Generating the binomial row and the powers
\(n^0,\ldots,n^n\) successively then evaluates the finite sum for
\(J_n\) with \(O(n)\) operations and \(O(n)\) stored integers.
For a nonnegative integer \(N\), sharing the \(d_k\) gives
\(J_0,\ldots,J_N\) with
\(O(N^2)\) operations and \(O(N)\) stored integers.

For the bit bound,
\[
 \log D(t)=\sum_{m\geq1}
    \left(1-\frac1m+\frac12\mathbf1_{\{2\mid m\}}\right)t^m
\]
has nonnegative coefficients, so \(D\) does also. Consequently,
for \(0\leq k\leq n\),
\[
 0\leq d_k\leq D(1/2)\,2^k k!
       =\frac{e^{7/6}}2\,2^k k!.
\]
The recurrence multipliers are polynomials of degree at most four,
so their products and signed sums have \(O(n\log n)\) bits for
\(n\geq2\). The binomial coefficients, powers, and final sum have
the same bound. Thus the prefix computation uses integers of
\(O(N\log N)\) bits for \(N\geq2\).

\Needspace{18\baselineskip}
\begin{corollary}
\label{cor:coupled-height-two-cost}
For a nonnegative integer \(N\), the formulas of
\namedref{Theorem}{thm:coupled-height-two-total} compute
all terms through \(H_N\) with \(O(N^4)\) arithmetic operations and
\(O(N^3)\) stored integers.
\end{corollary}

\begin{proof}
Precompute Stirling numbers, binomial coefficients, falling factorials,
and powers of \(x\) and \(x-1\) through degree \(N\) for integers
\(1\leq x\leq N\). For each such \(x\) and each integer
\(0\leq r<N\), precompute the powers of \(x(r+1)-1\)
through degree \(N\). This takes
\(O(N^3)\) operations and stored integers.
Equation~\eqref{eq:coupled-marked-diagonal} then computes all
\(\Delta_{k,k}(x)\), for \(1\leq k,x\leq N\), with another
\(O(N^3)\) operations.

There are \(O(N^3)\) remaining values
\(\Delta_{m,a}(x)\) with \(1\leq a<m\leq N\) and
\(1\leq x\leq N\); each requires \(O(N)\) operations by
\eqref{eq:coupled-marked-stirling}.
The two-component sums over \(0\leq n\leq N\) likewise have
\(O(N^4)\) summands, and the final recurrence needs \(O(N^2)\)
operations. The values \(\Delta_{m,a}(x)\) and the arithmetic
tables use \(O(N^3)\) stored integers.
Here arithmetic operations are counted independently of operand bit lengths.
\end{proof}

\end{document}